\documentclass[preprint,12pt]{elsarticle}
\usepackage[utf8]{inputenc}
\usepackage{amssymb}
\usepackage{amsmath}
\journal{Computer Methods in Applied Mechanics and Engineering}

\usepackage{xspace}
\usepackage{graphicx} % Required for inserting images
\usepackage{algorithm} %HB: required for pseudocode
\usepackage{amsfonts} 
\usepackage{amsmath}
\usepackage{siunitx}
\usepackage{subcaption}
\usepackage{amsthm}
\theoremstyle{remark}
\newtheorem{remark}{Remark}[section]
\usepackage{bm}
\usepackage{breakurl} 
\usepackage[breaklinks]{hyperref}
\usepackage{diagbox}   % for diagonal cell headers
\usepackage{booktabs} % for professional looking tables
\usepackage{doi}
\usepackage{textcomp}

\providecommand{\norm}[1]{\lVert#1\rVert}

\newcommand{\normnew}[1]{\left\lVert#1\right\rVert}
\newcommand{\Rmark}

\usepackage{algorithm}
\usepackage{algorithmic}
\usepackage{xcolor}
\definecolor{ao(english)}{rgb}{0.0, 0.5, 0.0}

\begin{document}

\begin{frontmatter}



\title{Data-driven model order reduction for wave propagation in materials with damage and nonlinearities}


\author[TUBS]{Saddam Hijazi\corref{cor1}}
\ead{saddam.hijazi@tu-braunschweig.de}
\author[HSU]{Nikiema Fulgence}
\ead{nikiemaw@hsu-hh.de}
\author[TUBS]{Hannah Burmester}
\ead{hannah.burmester@tu-braunschweig.de}
\author[HSU]{Natalie Rauter}
\ead{natalie.rauter@hsu-hh.de}
\author[TUBS]{Carmen Gr\"a\ss{}le}
\ead{c.graessle@tu-braunschweig.de}

\affiliation[TUBS]{organization={Technische Universit\unexpanded{ä}t Braunschweig},%Department and Organization
            addressline={Institute for Partial Differential Equations, Universit\unexpanded{ä}tsplatz 2}, 
            city={Braunschweig},
            postcode={38106}, 
            %state={},
            country={Germany}}
\affiliation[HSU]{organization={Helmut-Schmidt-University / University of the Federal Armed Forces Hamburg},%Department and Organization
            addressline={Chair of Solid Mechanics, Holstenhofweg 85}, 
            city={Hamburg},
            postcode={22043}, 
            %state={},
            country={Germany}}    
                    
\cortext[cor1]{Corresponding author}

%% Abstract
\begin{abstract}
In this work, we consider wave propagation in materials with damage and/or nonlinearities. To accelerate the simulations of the resulting high-dimensional problems, we apply model order reduction methods. Depending on the knowledge of the underlying equations and the availability of their discrete operators, intrusive methods (here projection-based approaches based on proper orthogonal decomposition (POD)) or non-instrusive methods (here data-driven approaches including several variants of dynamic mode decomposition (DMD) and operator inference (OpInf)) can be used. We have transferred the existing reduction approaches to the second order mechanical systems. In three different numerical examples, we evaluate the performance of the reduction techniques.
\end{abstract}

\begin{highlights}
\item Model order reduction methods for second order mechanical systems.
\item Targeted second order mechanical systems exhibit a damage/nonlinearity. 
\item Comparison between intrusive and non-intrusive methods for second order systems.
\item Reduced method learns reduced models preserving the structure of mechanical systems.
\end{highlights}

%% Keywords

\end{frontmatter}



\section*{Novelty statement}
We apply the methods proper orthogonal decomposition (POD), dynamic mode decomposition (DMD), multiresolution DMD (mrDMD) and variants of operator inference (OpInf) to different problems of guided ultrasonic wave (GUW) propagation  which are the classical  wave equation, a mechanical system resulting from GUW propagation and a nonlinear hyperelastic model. The operator inference approach is modified by following aspects: We present a time discrete operator inference approach for second order mechanical systems. In addition, the re-projection method which recovers the intrusive reduced operators has been extended for second order mechanical systems. Further, we continue and extend our study on the nonlinear aluminium model started in \cite{GH24}, by introducing a damage in the nonlinear model.

\section{Introduction}
\noindent Partial differential equations (PDEs) can be used to describe  physical phenomena based on conservation laws (e.g.\ mass and momentum) \cite{quarteroni2009numerical}. In structural mechanics, the governing equation of motion takes the form of a PDE derived from the balance of momentum \cite{hughes2012finite}. This PDE depends on parameters that influence the structural response, such as material properties, geometry, or boundary and initial conditions. In this work, motivated within the context of structural health monitoring (SHM) \cite{Farrar2006}, we consider a fixed but arbitrary damage configuration in a wave propagation model. Our overall research goal is to estimate the damage characteristics implemented by model parameters for subsequent localization and characterization. Since solving this inverse problem requires a large number of forward simulations, efficient and accurate surrogate models are needed. For this reason, this work undertakes a comparative analysis of intrusive and non-intrusive model reduction techniques - depending on the different levels of availability of model information - and uses both established approaches and extensions of reduction methods to the context of wave propagation.

%In this article, we focus on the problem of wave propagation in damaged materials within the context of SHM \cite{Farrar2006}. 
Non-destructive evaluation (NDE) techniques are commonly employed to detect and localize defects in structures. The objective of SHM is to leverage NDE methods to continuously monitor structural integrity in an automated manner, using appropriately placed sensors. Among the various approaches, a widely used technique is guided ultrasonic waves (GUW) \cite{amafabia2017review}. In this method, guided waves are propagated through the structure, and their interaction with damage alters the signals recorded by the sensors \cite{lammering2017lamb,su2006guided,guy2003guided,Fulgence2024,FulgenceNikiema2023}. The objective is to process these measured signals using numerical methods to characterize the presence, location, and severity of damage.
In this work, we consider different wave propagation models each including a material degradation defect:  a linear wave equation (Section~\ref{sec:waveEq}), a linear mechanical system resulting for GUW propagation (Section~\ref{sec:mechGUW}) and a nonlinear elastic material (Section~\ref{sec:mechSys2}). %We are interested in modeling the wave propagation when a damage is present in the physical domain. 
%In addition, we consider wave propagation in nonlinear elastic material without defects (Section~\ref{sec:mechSys2}) as a first step to later include damage in that setting. The type of damage considered in this work is always a material degradation damage.
%In this work, we focus on parameters associated with material degradation, often referred to as \emph{damage parameters}. The motivation stems from structural health monitoring (SHM), where the identification and localization of damage is essential to prevent catastrophic failures and to reduce maintenance costs \cite{Farrar2006}.
% Solving parameterized PDE problems for multiple parameter values is often computationally prohibitive. The FEM discretization typically yields systems of very large dimension, leading to high computational demands in terms of both time and memory. Repeated simulations, as required in real-time control, uncertainty quantification, inverse problems, or optimization, can therefore become infeasible. To address this issue, the scientific community has developed a wide range of \emph{model order reduction} (MOR) techniques over the last decades 
For a general overview of model reduction techniques we refer to e.g.\  \cite{hesthaven2015certified,Bader2016,prud2002reliable, Barrault2004, Grepl2007, Benner2015, quarteroni2014reduced, quarteroniRB2016}. These methods aim to construct low-dimensional surrogate models that retain the essential dynamics of the high-fidelity system while significantly reducing computational costs. Particularly effective reduction is achieved for problems in computational mechanics, encompassing both structural and fluid dynamics applications \cite{Rozza2008,Guo2018,Kunisch2002,Carlberg2013,Ballarin2016}.\\
 %For spatial discretization, we use the finite element method for both the wave equation and the equation of motion. In the case of the GUW simulations, the discretization in space and time is fine due to the high frequency of the simulation. This results in a high number of degrees of freedom of the discretized model raising the computational cost needed to perform the simulation. This motivates the use of model order reduction techniques in order to mitigate the computational costs associated with running the GUW simulations such that the reduced model can later be used to estimate the damage parameters solving an inverse problem. 
Model order reduction techniques can broadly be classified into intrusive and non-intrusive approaches. Intrusive methods (including e.g.\ proper orthogonal decomposition (POD) \cite{volkwein2013proper,Bergmann2009516,Baiges2014189,Burkardt2006,Ballarin2016,Graessle2021}, proper generalized decomposition (PGD) \cite{Dumon20111387,Chinesta2011} and greedy reduced basis method \cite{hesthaven2015certified}) require full knowledge of the governing equations and direct access to the discretized system matrices. In contrast, non-intrusive methods construct reduced models either solely from solution data without explicit access to the full order operators (e.g.\ dynamic mode decomposition (DMD) \cite{Schmid2010} or methods using regression techniques or artificial neural networks (ANNs) \cite{Loiseau2018,Guo2018,Hesthaven2018,Guo2019}) or require additional knowledge about the structure of the underlying system (e.g.\ operator inference (OpInf)). For the mechanical examples arising from the GUW simulations, we focus on non-intrusive techniques due to restricted access to the system's operators caused by using a commercial software.
In particular, we investigate in different variants of OpInf. Originally introduced in \cite{Peherstorfer2016Opinf}, the OpInf approach has since been extended in various directions \cite{Filanova2023,Benner2020,Qian2020,Yldz2021,SHARMA2024116865,Uy2023,Freitag2025}.  In \cite{Qian2022}, the authors extend the OpInf framework to nonlinear PDEs by formulating the methodology directly in a Hilbert space setting. Their approach incorporates the underlying inner product into the reduced modelling process and enables the inference of reduced operators that remain consistent with the continuous PDE formulation. Furthermore, the framework can exploit lifting transformations to expose polynomial structure in systems with non-polynomial nonlinearities.\par

The wave-type problems considered in this work arise from the spatial discretization of PDEs, such as the equation of motion for elastic bodies or the classical wave equation. The resulting second order mechanical system (often referred to simply as the \textit{mechanical model} in the literature) consists of a system of second order ordinary differential equations governing the (generalized) displacement vector $\mathbf{u}\in\mathbb{R}^{N}$, where $N$ denotes both the dimension of the state space and the number of ordinary differential equations.

For linear mechanical systems, the governing equations generally contain the following terms:

\begin{itemize}
\item \textit{Inertial forces}: $\mathbf{M}\ddot{\mathbf{u}}$, where $\mathbf{M}\in\mathbb{R}^{N\times N}$ is the mass matrix.
\item \textit{Damping forces}: $\mathbf{D}\dot{\mathbf{u}}$, where $\mathbf{D}\in\mathbb{R}^{N\times N}$ is the damping matrix.
\item \textit{Internal forces}: $\mathbf{K}\mathbf{u}$, where $\mathbf{K}\in\mathbb{R}^{N\times N}$ is the stiffness matrix.
\item \textit{External forces}: $\mathbf{f}\in\mathbb{R}^{N}$.
\end{itemize}

Typically, the mass and stiffness matrices, $\mathbf{M}$ and $\mathbf{K}$, are symmetric, with $\mathbf{M}$ being positive definite and $\mathbf{K}$ positive (semi-)definite. Under these standard assumptions, the undamped mechanical system is Lyapunov stable. Moreover, if the damping matrix $\mathbf{D}$ is positive definite, the system becomes asymptotically stable \cite{kailath1980linear,Muller1972Stability,Thomson2021Vibrations,Salimbahrami2006SecondOrderKrylov}.

In this work, our objective is to construct a reduced order model (ROM) that serves as an accurate surrogate for the full order mechanical model (FOM), providing a faithful approximation of the displacement field. An intrusive approach based on POD and Galerkin projection can be employed to achieve this objective. A key advantage of the POD--Galerkin methodology is that it preserves important structural properties of the original mechanical system, including its second order formulation. In particular, the reduced operators $\mathbf{M}_r$ and $\mathbf{K}_r$ (and $\mathbf{D}_r$ in case of presence of damping) inherit the symmetry and positive-definiteness properties of their full order counterparts. Consequently, the reduced mechanical model retains the stability characteristics of the original system while preserving its second order structure.

Alternatively, the second order system may first be rewritten as an equivalent first order system by introducing the momentum as an additional state variable. In this case, the dimension of the state vector, and consequently the number of ordinary differential equations, doubles, i.e.,
\[
N_{\mathrm{1st}} = 2N.
\]
Standard projection-based reduction techniques, such as POD--Galerkin, may then be applied directly to the reformulated first order system. However, this reformulation no longer explicitly preserves the second order mechanical structure and may introduce undesirable numerical behavior.

For conservative mechanical systems, the first order formulation possesses a Hamiltonian structure governed by the skew-symmetric symplectic matrix $\mathbb{J}$ \cite{leimkuhler2004simulating,hairer2006geometric}. Since the standard POD--Galerkin projection does not, in general, preserve the underlying symplectic geometry, the resulting reduced order model may fail to reproduce the qualitative dynamics of the original system. This geometric inconsistency can lead to artificial energy dissipation or spurious numerical energy growth, which may ultimately result in unstable long-time simulations. Consequently, special care must be taken when applying projection-based model reduction techniques to first order systems obtained from second order mechanical models.

To address this issue, several structure-preserving reduction techniques have been proposed, among which symplectic POD is one of the most prominent. The central idea is to restrict the reduced trial space to a symplectic subspace so that the reduced basis satisfies the required symplectic structure-preserving conditions \cite{lall2003structure,peng2016symplectic}. Throughout this work, whenever reduced order models are constructed from the governing second order mechanical equations, the original second order formulation is retained rather than being reformulated as an equivalent first order system. This choice preserves the intrinsic structure of the governing equations and avoids the potential difficulties associated with first order formulations that require structure-preserving reduction techniques. Purely data-driven methods, such as DMD, do not rely on the governing equations in the first place and are therefore outside the scope of this discussion.

It is well-known that hyperbolic problems and systems involving transport phenomena are challenging for model order reduction due to the Kolmogorov barrier, see e.g.\ \cite{OR16,BSU19,GU19}. In order to overcome this barrier, nonlinear model reduction methods can be employed which are not restricted to fixed linear subspaces, see e.g.\ \cite{Peherst2022,HestPehUng26}.

Finally, we also highlight hybrid reduced order models which leverage both physical equations and solution data. Examples include calibration-based methods \cite{noack2005,GALLETTI2004,Couplet2005}, data-driven filtering and correction techniques \cite{Xie2018,Ivagnes2023}, POD with interpolation \cite{Demo2019,Nguyen2022,Hijazi2020,Hijazi2020JCP}, and the use of physics-informed neural networks (PINNs) at the reduced level \cite{Chen2021,Hijazi2023}.
%{\color{red} integrate the first paragraph of section 2.2 here}

\section{Model description}

%In the following the models considered in this work are introduced. Apart from the classical wave equation, mechanical systems arising from GUW simulations are used.

\subsection{Wave equation}\label{sec:waveEq}
%The first model we consider is based on the classical wave equation, in which damage is incorporated via parameterization. Following  \cite{NandaLorenzWave}, damage is characterized by its spatial coordinates and by a local reduction of the wave propagation speed. Here, we study a modification of this model as proposed in \cite{Bur2024} by additionally introducing the damage size and considering a smoothing modification of the initial condition and the damage parameterization.
\noindent Let $(t_0,t_{end})$ be a given  time interval, $\Omega \subset \mathbb{R}^2$ denotes the spatial domain with sufficiently smooth boundary $\partial \Omega$. The scalar initial boundary value problem of the wave equation reads as follows:
\begin{equation}\label{eq:ibvp_wave}
\left\{\begin{aligned}
 u_{tt}(t,\mathbf{x}) - c^2(\mathbf{x}) \Delta u(t,\mathbf{x}) &= f(t,\mathbf{x}), && \quad \text{in } (t_0,t_{end}) \times \Omega,  \\
  u(t,\mathbf{x}) & = b(t,\mathbf{x}), && \quad \text{on } (t_0,t_{end}) \times \partial\Omega,\\
 u(t_0,\mathbf{x}) &= u_0(\mathbf{x}) , &&\quad \text{in } \Omega, \\
 u_t(t_0,\mathbf{x}) &= u_1(\mathbf{x}) , && \quad \text{in } \Omega.
\end{aligned}\right.
\end{equation}
% note: we take here (t0,tend) since in the numerical tests, we start at t=0.1 and end at t=5
Here, $u$ represents the displacement, $c$ is the wave speed coefficient and $f$ is the external forcing term. The initial conditions for the solution and the first time derivative at $t=t_0$, respectively, are given by $u_0$ and $u_1$. We follow the ideas in \cite{NandaLorenzWave}, where a damage is modeled by a local reduction of the wave propagation speed $c$ and in \cite{Bur2024} where an additional parameter is added describing the size of the damage. We previously used this model in the overview article \cite{GH24} and in the context of parameter selection procedures \cite{pamm_proceeding}. For the numerical solution of \eqref{eq:ibvp_wave}, we use piecewise linear finite elements for discretization in space, resulting in a system of ordinary differential equations, see Section~\ref{sec:mechSys2}. For the discretization in time, we use the a $\theta$-scheme, see Section~\ref{sec:waveEqNum} for more details.

\subsection{Mechanical system resulting from GUW propagation}\label{sec:mechGUW}
\noindent The second model is a mechanical system resulting from guided ultrasonic wave (GUW) propagation in the presence of damage. 
%{\color{red} In the field of structural analysis, it is possible to generate Lamb wavages in thin-walled structures making use of the physical phenomena of wave propagation interacting with the defects present in the structure.}\par 
The governing equation which describes the virtual displacement resulted from the GUW reads as:
\begin{equation}\label{eq:eqofmotion}
Div \, \mathbf{P} + \rho_0 \mathbf{b} = \rho_0 \mathbf{\ddot{u}}.
\end{equation}
Here, $Div$ is the material divergence operator, $\mathbf{P}$ is the first Piola-Kirchhoff stress tensor, $\rho_0$ is the density of the material, $\mathbf{b}$ is the density of the volumetric forces and $\mathbf{u}$ is the displacement. The equation of motion  \eqref{eq:eqofmotion} is solved together with the following boundary conditions for both the displacement $\mathbf{u}$ and the stress tensor  $\mathbf{P}$:
\begin{equation}\label{eq:BC_CurrentConfig}
\mathbf{u} = \mathbf{\bar{u}} \quad \textrm{at} \quad S_u, \quad \textrm{and} \quad \mathbf{P} \mathbf{N} = \mathbf{\bar{t}^{I}} \quad \textrm{at} \quad S_{\sigma}.
\end{equation} 
Here, the boundary of the domain is given by $\partial \Omega = S_u \cup S_\sigma$, thus $\mathbf{\bar{u}}$ is the value of the displacement field at the Dirichlet boundary $S_u$, $\mathbf{N}$ is the normal vector and $\mathbf{\bar{t}^{I}}$ is the strain vector at the Neumann boundary $S_\sigma$. In addition we have initial conditions which describe the displacement and first time derivative as follows:
\begin{equation}\label{eq:IC_CurrentConfig}
\mathbf{u} = \mathbf{u}_0 \quad \textrm{at} \quad t=0, \quad \textrm{and} \quad \mathbf{\dot{u}} = \mathbf{\dot{u}}_0 \quad \textrm{at} \quad t=0.
\end{equation} 
The equation of motion is discretized using the finite element method, employing  quadrilateral plane strain elements and quadratic Lagrange shape functions for spatial discretization. The discretization is conducted within the commercial software of COMSOL Multiphysics{{\textsuperscript{\textregistered}\xspace}}.
In the numerical experiments, we consider the linear two-dimensional model of fiber metal laminate (FML) with damage considered in \cite{BellamMuralidhar2021,BellamMuralidhar2023,Abschlussband2025}, see Section~\ref{sec:mechSys} for more details.

\subsection{Nonlinear hyperelastic model}\label{sec:nonlmodel}

%{\color{red} @Natalie: is this nonlinear Aluminium model already published somewhere? we refer here to the COMSOL model given to us by Niki\\
%The model is well-known and published in all major books dealing with nonlinear continuum mechanics. You could cite for example 978-1-118-63270-3}
\noindent The third model involves hyperelastic material with nonlinearity, see e.g.\ \cite{Belytschko2000}. For this purpose we give a brief introduction of how the nonlinearity enters the mathematical formulation of the problem. As in the last section, we have the governing equation being the equation of motion in \eqref{eq:eqofmotion}. We introduce the deformation gradient $\mathbf{F}$ which relates to the displacement through the following relation
\begin{equation}
    \mathbf{F} = \mathbf{I} + \nabla \mathbf{u},
\end{equation}
where $\mathbf{I}$ is the unity tensor. The determinant of $\mathbf{F}$ is called the Jacobian of the motion and is denoted by $J$. We further introduce rotation-independent deformation tensor denoted by $\mathbf{C} = \mathbf{F}^T \mathbf{F}$ which is called the right Cauchy-Green deformation tensor. In addition, we denote by $\mathbf{S}$ the second Piola-Kirchoff tensor which relates to $\mathbf{P}$ by the relation $\mathbf{P} = \mathbf{S} \mathbf{F}^T$. The equation of motion could be re-written using the second Piola-Kirchoff tensor $\mathbf{S}$.\par 

If we assume hyperelastic material behavior, these materials are characterized by the existence of a strain energy density function $\Psi$ as the potential of the stresses \cite{Belytschko2000,Hahn1985}. Based on this, the second Piola-Kirchoff stress tensor is given as:
\begin{equation}
\mathbf{S} = 2 \frac{\partial \Psi}{\partial \mathbf{C}}.
\end{equation}

We address one example of nonlinear material models which is the well known Neo-Hookean model. The Neo-Hookean material model serves as a simple extension of Hooke's law to large deformations. For compressible isotropic material behavior, the strain energy function as a function of the right Cauchy-Green deformation tensor and the determinant of the material deformation tensor is given by:
\begin{equation}
\mathbf{\Psi}_\textrm{Neo}(\mathbf{C}) = \frac{1}{2} \lambda {(\ln J)}^2 - \mu \ln \left(J\right) + \frac{1}{2}\mu  (tr(\mathbf{C})-3),
\end{equation}
where $\lambda$ is the first Lam\'e parameter and $\mu$ is the shear modulus or the second Lam\'e parameter \cite{Pence2014}. Then the resulting stress tensor is given by:
\begin{equation}
\mathbf{S} = \mu \ln(J) \mathbf{C}^{-1} + \mu (\mathbf{I} - \mathbf{C}^{-1}).
\end{equation}
The equation of motion is coupled with boundary and initial conditions as in \eqref{eq:BC_CurrentConfig} and \eqref{eq:IC_CurrentConfig}.\par

\subsection{Second order mechanical system}\label{sec:mechSys2}

Spatial discretization of the wave equation \eqref{eq:ibvp_wave} or the equation of motion \eqref{eq:eqofmotion} leads to a linear second order system of ordinary differential equations of the form:
\begin{equation}\label{eq:mechanicalSystem}
\mathbf{M} \ddot{\mathbf{u}}(t) + \mathbf{K} \mathbf{u}(t) = \mathbf{f}(t), \quad t \in [0,T].
\end{equation}
In case of the nonlinear hyperelastic model in Section~\ref{sec:mechSys2}, the spatial discretization leads to 
\begin{equation}\label{eq:mechanicalSystemNonl}
\mathbf{M}\ddot{\mathbf{u}}(t) + \mathbf{g}(\mathbf{u}(t)) = \mathbf{f}(t), \quad t \in [0,T].
\end{equation}
By $t$ we denote time and $[0,T]$ is the time interval of interest. The vector $\mathbf{u} \in \mathbb{R}^N$ collects the nodal  displacements\footnote{Note that we use the notation $\mathbf{u}$ to refer to both the continuous displacement field in \eqref{eq:eqofmotion} and its discrete counterpart in \eqref{eq:mechanicalSystem}; the intended interpretation is clear from the context.}. The matrices $\mathbf{M} \in \mathbb{R}^{N\times N}$ and $\mathbf{K} \in \mathbb{R}^{N\times N}$ are the mass and stiffness matrices, respectively,   $\mathbf{f}(t) \in \mathbb{R}^N$ is the external force vector, $\mathbf{g}$ is the nonlinear part and $N$ is the number of degrees of freedom. Often, the right hand side $\mathbf{f}$ in \eqref{eq:mechanicalSystem} is expressed in terms of input signals $\mathbf{z}(t) \in \mathbb{R}^{N_I}$, with $N_I$ input signals, and an input-mapping matrix $\mathbf{B} \in \mathbb{R}^{N \times N_I}$ as follows:
\begin{equation}\label{eq:mechanicalSystem_Input}
\mathbf{M} \ddot{\mathbf{u}}(t) + \mathbf{K} \mathbf{u}(t) = \mathbf{B} \mathbf{z}(t).
\end{equation}

\begin{remark}[Damping] In this work, the effect of model order reduction on the numerical simulation of wave propagation and its interaction with damage is analyzed, where the interaction of the propagating waves with damage is limited to the generation of additional wave modes and mode conversion. Since these effects are not influenced by damping, it is not included in the numerical model. 
\end{remark}
\section{Model order reduction techniques}\label{sec:morsec}
The aim of this work is to provide a broad comparison between different variants of model reduction applied to second order systems including intrusive and non-intrusive techniques depending on the availability of information on the structure of the system. For this we either utilize well-known approaches, propose modifications to them or extend existing methods to our application examples. In particular, we consider POD-Galerkin ROM in Section~\ref{sec:intrusivePOD}, variants of operator inference (OpInf) in Sections~\ref{sec:Opinf1},\ref{sec:Opinf2},\ref{sec:Opinf3},\ref{sec:Opinf4},\ref{sec:Opinf5} and different versions of DMD in Section~\ref{sec:DMD}.

\subsection{Snapshot data}\label{sec:snapdata}
Let us introduce the notation for the snapshot data to which we refer in the subsequent chapters.
At given equidistant time instants $0 = t_1 < \dots < t_M = T$ the data is gathered which  consists of the input signals $\mathbf{z}(t)$ and the nodal displacements of the mechanical system collected in the vector $\mathbf{u}(t)$ at the given time points $t=t_i$ for $i=1,\dots,M$. The snapshot matrices of the displacement and input signals are formed as follows:
\begin{align}\label{eq:data1}
\mathbf{U} &= 
\begin{bmatrix} 
\vert & \cdots & \vert \\
\mathbf{u}(t_1) & \cdots &  \mathbf{u}(t_M)   \\
\vert & \cdots & \vert
\end{bmatrix}  \in \mathbb{R}^{N \times M}, \quad 
\mathbf{Z} = 
\begin{bmatrix} 
\vert & \cdots & \vert \\
\mathbf{z}(t_1) & \cdots &  \mathbf{z}(t_M)   \\
\vert & \cdots & \vert
\end{bmatrix} \in \mathbb{R}^{N_I \times M}. 
\end{align}

\subsection{The intrusive POD-Galerkin ROM for second order mechanical systems}\label{sec:intrusivePOD}

In this subsection, we recall the intrusive POD-Galerkin ROM for a mechanical system of second order without damping. For the derivation of a reduced order approximation of the full order model (FOM) \eqref{eq:mechanicalSystem}  and  \eqref{eq:mechanicalSystemNonl}, respectively, the FOM displacement field $\mathbf{u}(t)$ is approximated using a Galerkin ansatz as follows:
\begin{equation}\label{eq:2}
\mathbf{u}(t) \approx \mathbf{\Phi} \mathbf{u}_r(t),
\end{equation}
with coefficient vector $\mathbf{u}_r(t) \in \mathbb{R}^r$ of dimension $r \ll N$ and reduced basis matrix $\mathbf{\Phi} = [\bm{\phi}_1,  \dots, \bm{\phi}_r ] \in \mathbb{R}^{N\times r}$. Here, the reduced basis vectors $\{\bm{\phi}_i\}_{i=1}^r$ are generated using the POD method. For this a SVD of the snapshot matrix of the displacement can be used:
\begin{equation}\label{eq:svdeq}
\mathbf{U} = \mathbf{V} \mathbf{\Sigma} \mathbf{W}^T,
\end{equation}
where $\mathbf{\Sigma}$ is a matrix containing the singular values in descending order, $\mathbf{V} \in \mathbb{R}^{N \times N}$ is the matrix of left singular vectors and $\mathbf{W} \in \mathbb{R}^{M \times M}$ is the matrix of right singular vectors. Then, a reduced basis $\mathbf{\Phi}$ is defined as the first $r$ columns of the matrix $\mathbf{V}$.\par  

By inserting the reduction assumption \eqref{eq:2} in \eqref{eq:mechanicalSystem}  and \eqref{eq:mechanicalSystemNonl}, respectively and performing a Galerkin projection onto the POD space, one may obtain the following linear reduced order dynamical system:
\begin{equation}
\mathbf{\Phi}^T \mathbf{M} \mathbf{\Phi} \; \ddot{\mathbf{u}}_r(t) + \mathbf{\Phi}^T \mathbf{K} \mathbf{\Phi} \; \mathbf{u}_r(t) = \mathbf{\Phi}^T \mathbf{f}(t), \quad \textrm{for} \quad t  \in  [0,T],
\end{equation}
and nonlinear system, respectively:
\begin{equation}
     \mathbf{\Phi}^T \mathbf{M} \mathbf{\Phi} \; \ddot{\mathbf{u}}_r(t) + \mathbf{g}(\mathbf{\Phi}  \mathbf{u}_r(t)) = \mathbf{\Phi}^T \mathbf{f}(t) , \quad \textrm{for} \quad t  \in  [0,T].
\end{equation}
We define the reduced operators as follows:
\begin{equation}
\mathbf{M}_r = \mathbf{\Phi}^T \mathbf{M}\mathbf{\Phi}, \quad \mathbf{K}_r = \mathbf{\Phi}^T \mathbf{K} \mathbf{\Phi}, \quad \mathbf{f}_r(t) = \mathbf{\Phi}^T \mathbf{f}(t).
\end{equation}
Thus, the linear reduced order system can be written as:
\begin{equation}\label{eq:3}
\mathbf{M}_r \ddot{\mathbf{u}}_r(t) + \mathbf{K}_r \mathbf{u}_r(t) = \mathbf{f}_r(t), \quad \textrm{for} \quad  t \in  [0,T].
\end{equation}
In case model \eqref{eq:mechanicalSystem_Input} is considered, we obtain the right-hand side
\begin{equation}\label{eq:Brz}
\mathbf{B}_r \mathbf{z}(t),  \quad \text{with} \quad  \mathbf{z}(t) \in \mathbb{R}^{N_I},
\end{equation}
where $\mathbf{B}_r = \mathbf{\Phi}^T \mathbf{B} \in \mathbb{R}^{r \times N_I}$ is the reduced control operator.

\subsection{The operator inference approach: motivation and data preparation}\label{sec:Opinf1}
The basic idea of operator inference is to recover a numerical model directly from snapshot data by approximating the system operators. Usually, the structure of the system that resulted in the numerical data is known to a certain extent. This partial knowledge is used to formulate a mathematical problem for the approximation of the system operators. Since the estimation of the full order operators would be computationally immense, one directly estimates reduced system operators which are of much smaller dimension. In \cite{Peherstorfer2016Opinf}, this approach has been proposed and is termed as operator inference. The operator inference approach has also been extended to specific physical problems and to non-polynomial nonlinear cases, see, e.g.\ \cite{Benner2020,Qian2020,Yldz2021}. An operator inference method tailored for second order mechanical systems has been proposed in \cite{Benner2021,Filanova2023}, which we closely follow here to apply it to our mechanical systems of interest.

In the operator inference approach, derivative information of the solution data is required. In many occasions, the numerical solver may provide this data. If that is not the case, then one may approximate the derivatives using finite difference schemes, e.g.\ using an eighth order finite difference scheme, see \cite{Sharma2024}. The second derivatives of the displacement are assembled in the following snapshot matrix:
\begin{align}\label{eq:data2}
\ddot{\mathbf{U}} &= 
\begin{bmatrix} 
\vert & \cdots & \vert \\
\ddot{\mathbf{u}}(t_1) & \cdots &  \ddot{\mathbf{u}}(t_M)   \\
\vert & \cdots & \vert
\end{bmatrix} \in \mathbb{R}^{N \times M}.
\end{align}
In order to obtain a reduced version of the snapshot data, POD can be used as follows:
\begin{equation}\label{eq:reduced_data}
\mathbf{\widehat{U}} = \mathbf{\Phi}^T \mathbf{U}, %\quad \dot{\mathbf{\widehat{U}}} = \mathbf{\Phi}^T \dot{\mathbf{U}}, 
\quad {\widehat{\ddot{\mathbf{U}}}} = \mathbf{\Phi}^T \ddot{\mathbf{U}},
\end{equation}
where $\mathbf{\Phi}$ is the matrix of POD modes computed using e.g.\ a singular value decomposition (SVD) of $\mathbf{U}$, see Section \ref{sec:intrusivePOD}.

\subsection{The operator inference approach for second order mechanical systems}\label{sec:Opinf2}

In this subsection, we recall the main ideas of the operator inference approach applied to the linear second order mechanical system and closely follow \cite{Peherstorfer2016Opinf,Filanova2023}.
 First, we write the reduced order system \eqref{eq:3} (but for the case of an input signal at the right hand side) for all time snapshots at once as follows:
\begin{equation}\label{eq:rom_all_snapshots}
\mathbf{M}_r \ddot{\mathbf{U}}_r + \mathbf{K}_r \mathbf{U}_r = \mathbf{B}_r \mathbf{Z}, 
\end{equation}
where $\mathbf{U_r}$ and $\ddot{\mathbf{U}}_r$ are the matrices of the ROM displacement and acceleration (second time derivative of the solution), respectively, and are defined as follows:

\begin{align}\label{eq:data_rom}
\mathbf{U}_r = 
\begin{bmatrix} 
\vert & \cdots & \vert \\
\mathbf{u}_r(t_1) & \cdots &  \mathbf{u}_r(t_M)   \\
\vert & \cdots & \vert
\end{bmatrix} \in \mathbb{R}^{r \times M}, \quad
\ddot{\mathbf{U}}_r &= 
\begin{bmatrix} 
\vert & \cdots & \vert \\
\ddot{\mathbf{u}}_r(t_1) & \cdots &  \ddot{\mathbf{u}}_r(t_M)   \\
\vert & \cdots & \vert
\end{bmatrix} \in \mathbb{R}^{r \times M}.
\end{align}
Multiplying both sides of \eqref{eq:rom_all_snapshots} by $\mathbf{M}_r^{-1}$ (assuming that $\mathbf{M}_r$ is invertible) from the left yields the following dynamical system:
\begin{equation}\label{eq:rom_modified_strucutre}
\ddot{\mathbf{U}}_r + \mathbf{M}_r^{-1} \mathbf{K}_r \mathbf{U_r} = \mathbf{M}_r^{-1} \mathbf{B}_r \mathbf{Z}.
\end{equation}
The reduced dynamical system \eqref{eq:rom_modified_strucutre} is used as a model for inferring the reduced operators by putting forward a least squares problem formulation. Namely, we seek to identify the second order mechanical system which takes the following form:
\begin{equation}\label{eq:opinf_rom1}
\widehat{\ddot{\mathbf{u}}}(t) + 
\mathbf{\widehat{K}}_M \mathbf{\widehat{u}}(t) =\mathbf{\widehat{B}}_M \mathbf{z}(t), \quad \textrm{for} \quad t  \in  [0,T],
\end{equation}
where $\mathbf{\widehat{K}}_M \in\mathbb{R}^{r\times r}, \mathbf{\widehat{B}}_M\in \mathbb{R}^{r\times N_I}$ and $ \mathbf{\widehat{u}}(t) \in \mathbb{R}^r$. This is done using the reduced data in \eqref{eq:reduced_data} and the input signals $\mathbf{z}(t)$. The system matrices are identified by solving the following minimization problem:
\begin{equation}\label{eq:opinf_forces_informed_unconstrained_pb}
\min_{\widehat{\mathbf{K}}_M, \widehat{\mathbf{B}}_M} 
\| \widehat{\ddot{\mathbf{U}}} + 
\widehat{\mathbf{K}}_M \widehat{\mathbf{U}} - \widehat{\mathbf{B}}_M \mathbf{Z} \|_F^2,
\end{equation}
where for a matrix $\mathbf{A} \in \mathbb{R}^{n \times m}$ the Frobenius norm is defined by $\normnew{\mathbf{A}}_F  = \sqrt{\sum^n_{i=1} \sum^m_{j=1} {|a_{ij}|}^2}$ with $a_{ij}$ denoting the entries of $\mathbf{A}$.\par
The optimization problem~\eqref{eq:opinf_forces_informed_unconstrained_pb} can be written in a more compact way by combining the data matrices as follows
\begin{equation}
\widehat{\mathbf{D}} = \begin{bmatrix}
  \mathbf{\widehat{U}}\\
  \mathbf{Z}
\end{bmatrix},
\end{equation}
and by defining the following combined unknown operators:
\begin{equation}
\widehat{\mathbf{P}} = \left[ - \widehat{\mathbf{K}}_M, \, \widehat{\mathbf{B}}_M \right],
\end{equation}
which allows to write the optimization problem as follows:
\begin{equation}\label{eq:opinf_min}
\min_{\widehat{\mathbf{P}} \in \mathbb{R}^{r \times (r+N_I)}} 
\norm{ \widehat{\mathbf{P}} \widehat{\mathbf{D}} -  \ddot{\widehat{\mathbf{U}}}}_F^2.
\end{equation}
The solution to the minimization problem above satisfies the normal equations:
\begin{equation}
\widehat{\mathbf{P}} \widehat{\mathbf{D}} \widehat{\mathbf{D}}^T  = \ddot{\widehat{\mathbf{U}}} \widehat{\mathbf{D}}^T .
\end{equation}
The equation above has a unique solution if and only if $\widehat{\mathbf{D}} \widehat{\mathbf{D}}^T$ is invertible which is fulfilled when $\widehat{\mathbf{D}}$ has full rank. The resulted operator $\widehat{\mathbf{P}}$ can then be used to integrate the second order mechanical system in \eqref{eq:opinf_rom1} and obtain the operator inference reduced solution. 
Then by multiplying the reduced solution from the left by the matrix of POD modes $\mathbf{\Phi}$, one can approximate the FOM displacement. We note that this approach requires no access to the full order model code or algorithm. Therefore, it falls under the category of non-intrusive reduced order models or data-driven reduced order models. The approach so far assumes partial knowledge about the external forces represented in the input signals $\mathbf{z}$ and is termed as \textit{unconstrained operator inference} \cite{Filanova2023}, the name is indicating that no constraints are imposed on the OpInf matrices. In the remaining part of this section we will recall another version of operator inference which is based on full knowledge of the forces. This version is termed \textit{forces-informed operator inference} \cite{Filanova2023}.\par 

Assume that the external forces or load vector $\mathbf{f}$ is available at each time instant of the numerical simulation:
\begin{align}
\mathbf{F} &= 
\begin{bmatrix} 
\vert & \cdots & \vert \\
\mathbf{f}(t_1) & \cdots &  \mathbf{f}(t_M)   \\
\vert & \cdots & \vert
\end{bmatrix} \in \mathbb{R}^{N \times M}.
\end{align}
In this case, the reduced solution of the intrusive ROM satisfies the following system:
\begin{equation}\label{eq:rom_all_snapshots_fullforces}
\mathbf{M}_r \ddot{\mathbf{U}}_r + \mathbf{K}_r \mathbf{U}_r = \mathbf{\Phi}^T \mathbf{F}.
\end{equation}
The goal is to infer the reduced operators using the data for the solution and the second time derivative and the data of the forces. This implies that we are seeking a reduced order system which takes the following form:
\begin{equation}
\mathbf{\widehat{M}} \widehat{\ddot{\mathbf{u}}}(t) + 
\mathbf{\widehat{K}} \mathbf{\widehat{u}}(t) = \mathbf{\Phi}^T \widehat{\mathbf{f}}(t),  \quad \text{for } t \in [0,T],
\end{equation}
with the condition that $\mathbf{\widehat{M}}$ and $\mathbf{\widehat{K}}$ are symmetric positive definite matrices:
\begin{equation}
\mathbf{\widehat{M}} \succ 0, \quad \mathbf{\widehat{K}} \succ 0.
\end{equation}
The reduced forces data can be computed in a similar way to the reduced displacement as $\mathbf{\widehat{F}} = \mathbf{\Phi}^T \mathbf{F}$.\par 

\noindent The forces-informed operator inference problem is given by:
\begin{equation}\label{eq:opinf_forces_informed_pb}
\min_{\widehat{\mathbf{M}} \succ 0, \widehat{\mathbf{K}} \succ 0} 
\normnew{  \widehat{\mathbf{M}} \widehat{\ddot{\mathbf{U}}} +
\widehat{\mathbf{K}} \mathbf{\widehat{U}} - \mathbf{\widehat{F}} }_F^2.
\end{equation}
In the current case, we define the combined data matrix:
\begin{equation}
\widehat{\mathbf{D}} = \begin{bmatrix}
  \widehat{\ddot{\mathbf{U}}}\\
  \mathbf{\widehat{U}}
\end{bmatrix}.
\end{equation}
Then, the operator inference forces-informed problem is given by:

\begin{equation}\label{eq:opinf_forces_informed_pb}
\min_{\widehat{\mathbf{M}} \succ 0, \widehat{\mathbf{K}} \succ 0} 
\normnew{ \left[ \widehat{\mathbf{M}} \quad 
\widehat{\mathbf{K}} \right] \widehat{\mathbf{D}}  - \mathbf{\widehat{F}} }_F^2.
\end{equation}

We remark that it is possible for certain optimization algorithms to transform the optimization problem in \eqref{eq:opinf_forces_informed_pb} into an unconstrained one. For this, we utilize a parameterization of stable matrices originally proposed in \cite{GILLIS2017113} and used in \cite{GoyalPontesBenner2023,Goyal2025GuaranteedStable}. This is done by assuming that the mass and the stiffness matrices admit the following decompositions:
\begin{equation}
\widehat{\mathbf{M}} = \mathbf{L}^T \mathbf{L}, \quad \widehat{\mathbf{K}} = \mathbf{E}^T \mathbf{E},
\end{equation}
which implies that the unconstrained optimization problem for the forces-informed operator inference is:
\begin{equation}\label{eq:opinf_forces_informed_pb_symDecompo}
\min_{\mathbf{L}, \mathbf{E} } 
\normnew{ \mathbf{L}^T \mathbf{L} \widehat{\ddot{\mathbf{U}}} + 
\mathbf{E}^T \mathbf{E} \mathbf{\widehat{U}}  - \mathbf{\widehat{F}} }_F^2.
\end{equation}
\begin{remark}[Operator inference for a nonlinear model]
We note that we consider operator inference here only for the linear mechanical model of Section~\ref{sec:mechGUW} due to the complex structure of the nonlinear model. In the work \cite{Qian2020}, the nonlinearity requires to be of polynomial structure or can be lifted to polynomial structure using a variable transformation. In \cite{Benner2020} an operator inference approach is introduced for systems with non-polynomial nonlinear terms which requires non-polynomial terms analytically.
\end{remark}

\subsection{Scaling the data for the operator inference approach}\label{sec:Opinf3}

In several scientific applications, one may encounter the issue of having data which belong to different numerical scales. This issue could pose a difficulty for the operator inference approach and in particular for solving the optimization problem both in the unconstrained and in the forces-informed case. Therefore, we propose to deal with this problem by introducing scaling of the data snapshot matrices as follows:
\begin{equation}
\tilde{\ddot{\mathbf{U}}} = \frac{\ddot{\mathbf{U}}}{\norm{\ddot{\mathbf{U}}}_{F}} , \quad \tilde{\mathbf{U}} = \frac{\mathbf{U}}{\norm{\mathbf{U}}_{F}}, \quad \tilde{\mathbf{Z}} = \frac{\mathbf{Z}}{\norm{\mathbf{Z}}_{F}}.
\end{equation}
In case of employing the forces-informed approach the matrix $\mathbf{F}$ will substitute the matrix $\mathbf{Z}$. The reduced order basis matrix $\mathbf{\Phi}$ is obtained by doing a SVD on the newly scaled matrix of the displacement snapshots:
\begin{equation}
\tilde{\mathbf{U}} = \tilde{\mathbf{V}} \tilde{\mathbf{\Sigma}} \tilde{\mathbf{W}}^T,
\end{equation} 
and by retaining the first $r$ columns of the matrix $\tilde{\mathbf{V}}$. The approach continues with computing the reduced data based on the scaled snapshots:
\begin{equation}\label{eq:reduced_data_scaling}
\mathbf{\widehat{U}} = \mathbf{\Phi}^T \mathbf{\tilde{U}}, \quad \widehat{\ddot{\mathbf{U}}} = \mathbf{\Phi}^T \tilde{\ddot{\mathbf{U}}}.
\end{equation}
Then, we solve the operator inference problem in \eqref{eq:opinf_forces_informed_unconstrained_pb} but with $\tilde{\mathbf{Z}}$ instead of $\mathbf{Z}$. The effects of the data scaling have to be considered when solving the reduced order system which takes the following form:
\begin{equation}\label{eq:scaledeq}
\frac{1}{\norm{\ddot{\mathbf{U}}}_{F}} \widehat{\ddot{\mathbf{u}}}(t) + 
\frac{1}{\norm{\mathbf{U}}_{F}} \widehat{\mathbf{K}}_M \mathbf{\widehat{u}}(t) = \widehat{\mathbf{B}}_M \tilde{\mathbf{z}}(t), \quad \text{with} \quad \tilde{\mathbf{z}}(t) = \frac{\mathbf{z}(t)}{\norm{\mathbf{Z}}_{F}}.
\end{equation}
After solving \eqref{eq:scaledeq}, an approximation of the full displacement can be computed by multiplying the reduced solution from the left by the matrix of POD modes $\mathbf{\Phi}$.

\subsection{The discrete operator inference approach for second order mechanical systems}\label{sec:Opinf4}

In Section~\ref{sec:Opinf2}, we have presented the operator inference approach for the second order mechanical system of interest written in a semi-discretized setting. However, the operator inference method could be also applied in a fully discretized setup. We refer the reader to the works \cite{Farcas2023,Peherstorfer2020} for an application of discrete OpInf to first order systems.\par

The continuous time formulation of OpInf requires the evaluation of time derivatives of the displacement field. These derivative data are typically obtained through finite difference approximations when they are not provided directly by the numerical solver. In practice, however, estimating derivatives from data may introduce significant numerical errors, particularly when the data are noisy or when the temporal resolution is limited. Such errors propagate directly into the regression problem used to identify the reduced operators and may deteriorate the stability and predictive accuracy of the inferred model.\par

An alternative strategy consists in formulating the operator inference problem directly at the level of the time discretized dynamics. Depending on the time integration scheme used for the full order model, the discrete evolution of the displacement field can be expressed as a recurrence relation linking several consecutive time levels. For instance, explicit or semi-implicit schemes such as central-difference or theta-type methods lead to relations of the following form:
\begin{equation}\label{eq:DiscreteFOM}
\mathbf{G} \mathbf{u}_{k+2} = \mathbf{H} \mathbf{u}_{k+1} + \mathbf{R} \mathbf{u}_k + \mathbf{C} \mathbf{z}_{k+2} ,
\end{equation}
where $\mathbf{G} \in \mathbb{R}^{N \times N}$, $\mathbf{H} \in \mathbb{R}^{N \times N}$, $\mathbf{R} \in \mathbb{R}^{N \times N}$ and $\mathbf{C} \in \mathbb{R}^{N \times N_I}$. Furthermore, $\mathbf{u}_k = \mathbf{u}(t_k)$ denotes the displacement vector at time $t_k$, while $\mathbf{z}_k = \mathbf{z}(t_k)$ represents the input signal. Equation \eqref{eq:DiscreteFOM} thus defines an update relation for computing the displacement vector $\mathbf{u}_{k+2}$ based on the two previous displacement states $\mathbf{u}_k$, $\mathbf{u}_{k+1}$ and the input vector $\mathbf{z}_{k+2}$.

It should be noted that not all full order time integration schemes can be written in the form of \eqref{eq:DiscreteFOM}. For instance, certain methods (such as the Newmark scheme) evolve the displacement and velocity variables simultaneously and therefore do not admit a recurrence relation involving displacement variables alone. Consequently, the fully discrete formulation of OpInf should be viewed as an additional modeling option that may be advantageous when the underlying time integration scheme naturally leads to a displacement-based recurrence relation.

Similar to \eqref{eq:2}, a Galerkin ansatz could be introduced for the displacement vector $\mathbf{u}_k$ as follows:
\begin{equation}\label{eq:DiscreteGalerkin}
\mathbf{u}_k \approx \mathbf{\Phi} \mathbf{u}_{k,r},
\end{equation}
where the reduced displacement at time step $t_k$ is $\mathbf{u}_{k,r} \in \mathbb{R}^r$. The substitution of \eqref{eq:DiscreteGalerkin} into \eqref{eq:DiscreteFOM} and the projection onto the reduced space spanned by the columns of $\mathbf{\Phi}$ yields:
\begin{equation}\label{eq:DiscreteGalerkinROM}
\mathbf{G}_r \mathbf{u}_{r,k+2} = \mathbf{H}_r \mathbf{u}_{r,k+1} + \mathbf{R}_r \mathbf{u}_{r,k} + \mathbf{C}_r \mathbf{z}_{k+2},
\end{equation} 
where
\begin{equation}
\mathbf{G}_r = \mathbf{\Phi}^T \mathbf{G}\mathbf{\Phi}, \quad \mathbf{H}_r = \mathbf{\Phi}^T \mathbf{H} \mathbf{\Phi} \quad \mathbf{R}_r = \mathbf{\Phi}^T \mathbf{R} \mathbf{\Phi}, \quad \mathbf{C}_r = \mathbf{\Phi}^T \mathbf{C}.
\end{equation}

The solution of the FOM and the ROM in \eqref{eq:DiscreteFOM} and \eqref{eq:DiscreteGalerkinROM}, respectively, is carried out by inverting the matrix $\mathbf{G}$ or $\mathbf{G}_r$ (assuming that they are invertible). In this section, we assume $\mathbf{G}_r = \mathbf{I}$. From the perspective of operator inference, this assumption does not restrict generality since the discrete system can be left-multiplied by $\mathbf{G}_r^{-1}$. Consequently, the operators identified by the inference procedure correspond to the effective matrices $\mathbf{G}_r^{-1}\mathbf{H}_r$, $\mathbf{G}_r^{-1}\mathbf{R}_r$, and $\mathbf{G}_r^{-1}\mathbf{C}_r$.\par

Unlike the time continuous formulation, the operator inference approach developed for the discrete time setting does not require time derivative data. The first step of the discrete OpInf approach consists of decomposing the displacement data into three matrices: the snapshot data, the one-step--shifted snapshots, and the two-step--shifted snapshots, defined as follows:
\begin{align}\label{eq:dataDiscreteOpInf}
\mathbf{U}_0 &= 
\begin{bmatrix} 
\vert & \cdots & \vert \\
\mathbf{u}_1 & \cdots & \mathbf{u}_{M-2}   \\
\vert & \cdots & \vert
\end{bmatrix} \in \mathbb{R}^{N \times (M-2)}, \quad 
\mathbf{U}_1 = 
\begin{bmatrix} 
\vert & \cdots & \vert \\
\mathbf{u}_2 & \cdots & \mathbf{u}_{M-1}   \\
\vert & \cdots & \vert
\end{bmatrix} \in \mathbb{R}^{N \times (M-2)}, \quad \\
\mathbf{U}_2 &= 
\begin{bmatrix} 
\vert & \cdots & \vert \\
\mathbf{u}_3 & \cdots & \mathbf{u}_M   \\
\vert & \cdots & \vert
\end{bmatrix} \in \mathbb{R}^{N \times (M-2)}, \quad
\mathbf{Z}_2 = 
\begin{bmatrix} 
\vert & \cdots & \vert \\
\mathbf{z}_3 & \cdots & \mathbf{z}_M   \\
\vert & \cdots & \vert
\end{bmatrix} \in \mathbb{R}^{N_I \times (M-2)}, 
\end{align}
The corresponding reduced data is computed as follows
\begin{equation}
\mathbf{\widehat{U}}_0 = \mathbf{\Phi}^T \mathbf{U}_0, \quad \mathbf{\widehat{U}}_1 = \mathbf{\Phi}^T \mathbf{U}_1, \quad 
\mathbf{\widehat{U}}_2 = \mathbf{\Phi}^T \mathbf{U}_2.
\end{equation}
Then, the least squares regression problem for operator inference reads as follows:
\begin{equation}
\min_{\widehat{\mathbf{P}} \in \mathbb{R}^{r \times (2r+N_I)}} 
\norm{ \widehat{\mathbf{P}} \widehat{\mathbf{D}} -  \mathbf{\widehat{U}}_2 }_F^2,
\end{equation}
where
\begin{equation}\label{eq:dataDiscreteOpInf_PandD}
\widehat{\mathbf{D}} = \begin{bmatrix}
  \mathbf{\widehat{U}}_1\\
  \mathbf{\widehat{U}}_0\\
  \mathbf{Z}_2
  \end{bmatrix}, 
  \quad \widehat{\mathbf{P}} = \left[  \widehat{\mathbf{H}}, \,  \widehat{\mathbf{R}}, \, \widehat{\mathbf{C}} \right].
\end{equation}

\subsection{The re-projection method for second order mechanical systems}\label{sec:Opinf5}

In this section we propose an extension of the re-projection method presented in \cite{Peherstorfer2020} to the case of second order mechanical systems. The re-projection method aims at generating reduced trajectories (which we denote with a bar $\mathbf{\bar{U}}$) which are identical to the ones produced by the intrusive reduced operators. This results in a zero closure error $ \norm{\mathbf{U}_r - \mathbf{\bar{U}}}_F = 0$. It has been shown in \cite{Peherstorfer2020} that OpInf applied on the re-projected trajectories recovers the intrusive reduced operators under certain conditions.\par 

The re-projection method aims at generating snapshot trajectories which enforce a zero closure error. The closure error is caused by the non-Markovian dynamics of the projected trajectories. The non-Markovian dynamics are not produced by the reduced order model and by eliminating them one could recover the trajectories generated by the reduced order model. This is done by a data sampling procedure which cancels out the non-Markovian dynamics in the reduced space. We introduce the re-projection method in \autoref{alg:reproj_cont} and \autoref{alg:reproj_disc} for the time continuous case and the time discrete case, respectively. It can be seen from \autoref{alg:reproj_cont} that we need to have both the first and second time derivative of the full order solve at each iteration in order to perform the method. The OpInf problem in the continuous time case will read as follows:
\begin{equation}\label{eq:opinf_reproj_cont}
\min_{\mathbf{\bar{P}} \in \mathbb{R}^{r \times (r+N_I)}} 
\norm{ \mathbf{\bar{P}} \mathbf{\bar{D}} -  \bar{\ddot{\mathbf{U}}}}_F^2, \quad \text{where} \quad \mathbf{\bar{D}} = \begin{bmatrix}
  \mathbf{\bar{U}}\\
  \mathbf{Z}
\end{bmatrix}, \quad \text{and} \quad \mathbf{\bar{P}} = \left[ - \mathbf{\bar{K}}_M, \, \mathbf{\bar{B}}_M \right].
\end{equation}
On the other hand, the discrete OpInf problem is:
\begin{equation}\label{eq:opinf_reproj_disc}
\min_{\mathbf{\bar{P}} \in \mathbb{R}^{r \times (2r+N_I)}} 
\norm{ \mathbf{\bar{P}} \mathbf{\bar{D}} - \mathbf{\bar{U}}_2 }_F^2, \quad \text{where} \quad \mathbf{\bar{D}} = \begin{bmatrix}
  \mathbf{\bar{U}}_1\\
  \mathbf{\bar{U}}_0\\
  \mathbf{Z}_2
  \end{bmatrix}, \quad \text{and} \quad \mathbf{\bar{P}} = \left[  \mathbf{\bar{H}}, \,  \mathbf{\bar{R}}, \, \mathbf{\bar{C}} \right].
\end{equation}

\begin{algorithm}
\caption{The re-projection method for second order mechanical systems (time continuous case)}\label{alg:reproj_cont}
\begin{algorithmic}[1]
  \scriptsize
  \STATE Start with the initial displacement $\mathbf{u}_1$ and the initial time derivative of the displacement $\mathbf{\dot{u}}_1$ with the input signal $\mathbf{z}_1$.
  \scriptsize \\
  \STATE Solve the FOM to obtain $\mathbf{\ddot{u}}_1$.
  \STATE Compute $\mathbf{\bar{u}}_1 =  \mathbf{\Phi}^T \mathbf{u}_1 \, \bar{\dot{\mathbf{u}}}_1 =  \mathbf{\Phi}^T \mathbf{\dot{u}}_1 \, \bar{\ddot{\mathbf{u}}}_1 =  \mathbf{\Phi}^T \mathbf{\ddot{u}}_1$
  \scriptsize \\
  \FOR{$i = 2, \ldots, M$}
    \STATE Solve the FOM for one time step using $\mathbf{\Phi} \mathbf{\bar{u}}_{i-1}, \, \mathbf{\Phi} \bar{\dot{\mathbf{u}}}_{i-1}, \, \mathbf{\Phi} \bar{\ddot{\mathbf{u}}}_{i-1}, \, \mathbf{z}_{i-1} \rightarrow \mathbf{u}_i^\star, \, \mathbf{\dot{u}}_i^\star, \, \mathbf{\ddot{u}}_i^\star$.
    \STATE Set $\mathbf{\bar{u}}_i =  \mathbf{\Phi}^T \mathbf{u}_i^\star, \, \bar{\dot{\mathbf{u}}}_i = \mathbf{\Phi}^T \mathbf{\dot{u}}_1^\star, \, \bar{\ddot{\mathbf{u}}}_i = \mathbf{\Phi}^T \mathbf{\ddot{u}}_i^\star$.
\ENDFOR
  \STATE Return the re-projected reduced trajectories $\mathbf{\bar{U}} = [\mathbf{\bar{u}}_1, \mathbf{\bar{u}}_2, \ldots, \mathbf{\bar{u}}_M]$ and $\bar{\ddot{\mathbf{U}}} = [\bar{\ddot{\mathbf{u}}}_1, \bar{\ddot{\mathbf{u}}}_2, \ldots, \bar{\ddot{\mathbf{u}}}_M]$ 
  \end{algorithmic}
\end{algorithm}

  \begin{algorithm}
\caption{The re-projection method for second order mechanical systems (time discrete case)}\label{alg:reproj_disc}
\begin{algorithmic}[1]
  \scriptsize
  \STATE Start with the displacement at the first two steps $\mathbf{u}_1$ and $\mathbf{u}_2$ with the input signal discretized vector $\mathbf{z}_1$.
  \scriptsize \\
  \STATE Compute $\mathbf{\bar{u}}_1 =  \mathbf{\Phi}^T \mathbf{u}_1 \, \mathbf{\bar{u}}_2 =  \mathbf{\Phi}^T \mathbf{u}_2$.
  \scriptsize \\
  \FOR{$i = 2, \ldots, M-1$}
    \STATE Solve the FOM for one time step using $\mathbf{\Phi} \mathbf{\bar{u}}_{i-1}, \, \mathbf{\Phi} \mathbf{\bar{u}}_i, \, \mathbf{z}_{i+1} \rightarrow \mathbf{u}_{i+1}^\star.$
    \STATE Set $\mathbf{\bar{u}}_{i+1} =  \mathbf{\Phi}^T \mathbf{u}_{i+1}^\star$.
\ENDFOR
  \STATE Return the re-projected reduced trajectories $\mathbf{\bar{U}}_0 = [\mathbf{\bar{u}}_1, \mathbf{\bar{u}}_2, \ldots, \mathbf{\bar{u}}_{M-2}]$, $\mathbf{\bar{U}}_1 = [\mathbf{\bar{u}}_2, \mathbf{\bar{u}}_3, \ldots, \mathbf{\bar{u}}_{M-1}]$ and $\mathbf{\bar{U}}_2 = [\mathbf{\bar{u}}_3, \mathbf{\bar{u}}_4, \ldots, \mathbf{\bar{u}}_M]$.
\end{algorithmic}
\end{algorithm}

As already mentioned in \cite{Peherstorfer2016,Peherstorfer2020}, we have observed that in several cases one may need to augment the data matrix $\mathbf{\bar{D}}$ or $\widehat{\mathbf{D}}$ with trajectories coming from different input signals in order to reduce the condition number of the matrix $\mathbf{\bar{D}} \mathbf{\bar{D}}^T$ or $\widehat{\mathbf{D}} \widehat{\mathbf{D}}^T$ which could pollute the learned operators and introduce numerical errors. For instance in the time continuous case, we assume that we have $q$ different input signals ${[\mathbf{z}^j(t)]}_{j=1}^q$ which give the matrices ${[\mathbf{Z}^j]}_{j=1}^q$, then we apply \autoref{alg:reproj_cont} $q$ times to obtain $\bar{\mathbf{U}}^1,\dots,\bar{\mathbf{U}}^q$ and $\bar{\ddot{\mathbf{U}}}^1,\dots,\bar{\ddot{\mathbf{U}}}^q$. We concatenate the matrices to obtain:
\begin{equation}
\mathbf{\bar{U}} = [\bar{\mathbf{U}}^1, \bar{\mathbf{U}}^2, \ldots, \bar{\mathbf{U}}^q], \quad \bar{\ddot{\mathbf{U}}} = [\bar{\ddot{\mathbf{U}}}^1, \bar{\ddot{\mathbf{U}}}^2, \ldots, \bar{\ddot{\mathbf{U}}}^q], \quad \mathbf{\bar{Z}} = [\mathbf{Z^1}, \mathbf{Z^2}, \ldots, \mathbf{Z^q}],
\end{equation} 
then we solve the OpInf \eqref{eq:opinf_reproj_cont} with the matrices above.

\subsection{Dynamic mode decomposition}\label{sec:DMD}

Dynamic mode decomposition (DMD) is a purely data-driven approach proposed in \cite{Schmid2010} to extract spatio-temporal coherent structures from high-dimensional dynamical systems. For a comprehensive study of the method, we refer to \cite{BN19,HTu2014,KBBP16}, for example. The basic idea is to treat a sequence of measurement snapshots as observations of an underlying linear operator (often interpreted as a finite-dimensional approximation of the Koopman operator, see \cite{ROWLEY2009}). In order to characterize the system's dynamics, DMD computes eigenvalues and eigenvectors of this linear operator. The resulting modes evolve exponentially in time allowing to identify growth rates or oscillatory patterns without explicitly knowing the underlying governing equations.

Let us recall the main steps of the DMD algorithm. We denote the discrete approximation at time $t_k$ by $\mathbf{u}_k \approx \mathbf{u}(t_k)$. The DMD technique seeks a best-fit linear operator $\mathbf{A} \in \mathbb{R}^{N \times N}$ such that
\begin{equation}\label{eq:dmdadvance}
\mathbf{u}_{k+1} \approx \mathbf{A} \mathbf{u}_k, \quad k = 1, \dots, M-1.
\end{equation}
Using the notations
\begin{equation}\label{eq:notationdmd}
\mathbf{U} = 
\begin{bmatrix} 
\vert &    & \vert \\
\mathbf{u}_1 &    \cdots &  \mathbf{u}_M   \\
\vert  & & \vert
\end{bmatrix} \in \mathbb{R}^{N \times M}
\end{equation}
and
\begin{equation}
\mathbf{U}_j^k =
\begin{bmatrix} 
\vert &   & \vert \\
\mathbf{u}_j  & \cdots &  \mathbf{u}_k   \\
\vert  & & \vert
\end{bmatrix}
\end{equation}
the DMD problem \eqref{eq:dmdadvance} can be written utilizing two time-shifted snapshot matrices as
\begin{equation*}
\mathbf{U}_2^M \approx  \mathbf{A} \mathbf{U}_1^{M-1}.
\end{equation*}
The solution is given by:
\begin{equation}
\mathbf{A} = \mathbf{U}_2^M {(\mathbf{U}_1^{M-1})}^\dagger,
\end{equation}
where $\dagger$ denotes the Moore--Penrose pseudoinverse. However, explicitly forming and storing the full operator $\mathbf{A} \in \mathbb{R}^{N \times N}$ and consequently calculating its spectral decomposition quickly becomes infeasible when $N$ is large. In order to circumvent this, a reduced operator $\tilde{\mathbf{A}} \in \mathbb{R}^{r \times r}$ with $r \ll N$ is constructed via 
\begin{equation*}
\tilde{\mathbf{A}} = \mathbf{V}^T \mathbf{A} \mathbf{V} = \mathbf{V}^T \mathbf{U}_2^M \mathbf{W} \mathbf{\Sigma}^{-1} \in \mathbb{R}^{r \times r},
\end{equation*}
where $\mathbf{V} \in \mathbb{R}^{N \times r}$ denotes the POD modes, $\mathbf{\Sigma} \in \mathbb{R}^{r\times r}$ contains the first $r$ leading singular values on the diagonal and $\mathbf{W}\in \mathbb{R}^{(M-1) \times r}$ has as its columns the associated right singular vectors from a truncated SVD of $\mathbf{U}_1^{M-1}$ given by
\begin{equation*}
\mathbf{U}_1^{M-1} \approx \mathbf{V} \mathbf{\Sigma} \mathbf{W}^T.
\end{equation*}
This reduced operator $\tilde{\mathbf{A}}$ allows to construct a model for the computation of the reduced coefficients $\tilde{\mathbf{u}}_k$ as
\begin{equation}\label{eq:dmdadvance_red}
\tilde{\mathbf{u}}_{k+1} = \tilde{\mathbf{A}} \tilde{\mathbf{u}}_k,
\end{equation}
where the full order state $\mathbf{u}$ can be reconstructed via $\mathbf{u}_k \approx \mathbf{V}\tilde{\mathbf{u}}_k$.
We compute the eigenvalues and eigenvectors of the reduced operator  $\tilde{\mathbf{A}}$ as
$$\tilde{\mathbf{A}} \mathbf{\Psi} = \mathbf{\Psi} \mathbf{\Lambda}.$$
It can be shown that the reduced operator $\tilde{\mathbf{A}}$ has the same non-zero eigenvalues as the full order operator $\mathbf{A}$, see \cite{HTu2014}. We use the notation $\mathbf{\Lambda}=\text{diag}(\lambda_1,\dots,\lambda_r)$ and refer to $(\lambda_i)_{i=1}^r$ as the DMD eigenvalues. The eigenvectors of the reduced operator $\tilde{\mathbf{A}}$ can be utilized in order to construct the DMD modes\footnote{Note that we use the notation $\mathbf{\Phi}$ to refer to the DMD modes in this Section and the POD modes in Section~\ref{sec:intrusivePOD}; the intended interpretation is clear from the context.} $\mathbf{\Phi} \in \mathbb{R}^{N \times r}$ as
$$\mathbf{\Phi} = \mathbf{U}_2^M \mathbf{W} \mathbf{\Sigma}^{-1} \mathbf{\Psi}.$$
These DMD modes are the eigenvectors of the full order operator $\mathbf{A}$ such that we can compute the spectral decomposition of $\mathbf{A}$ without the need to carry out explicit computations on the matrix $\mathbf{A}$ directly.
Utilizing the spectral decomposition $(\mathbf{\Phi},\mathbf{\Lambda})$ it is possible to reconstruct the system state, here written in a formulation which is continuous in time (see e.g.\ \cite[(7.27)]{BN19}):
\begin{equation}\label{eq:DMD_expression}
\mathbf{u}(t) = \sum_{j=1}^r b_j \bm{\phi}_j e^{\omega_j t} = \mathbf{\Phi} \, \exp(\Omega t) \, \mathbf{b}, \quad \omega_i = \log(\lambda_i)/\Delta t. 
\end{equation}
For the computation of the initial mode amplitudes $\mathbf{b}$, we consider two choices in this work:  the first option (see e.g.\ \cite[(7.29)]{BN19}) is given by
\begin{equation}\label{eq:amplitudes_using_projectedPOD}
    \mathbf{b}  = {\left( \mathbf{W} \mathbf{\Lambda} \right)}^{-1} \mathbf{\tilde{u}}_1,
\end{equation}
where the last expression computes the amplitudes using POD projected data. As second option, the optimal vector of amplitudes $\mathbf{b}$ could be computed by solving an optimization problem as mentioned in \cite{Jovanovi2014} (Section II.B). The final formula for the optimal $\mathbf{b}$ is:
\begin{equation}\label{eq:amplitudes_optimal}
\mathbf{b} = \left( \left( \mathbf{\Phi}^* \mathbf{\Phi} \right) \circ 
\left( \overline{\mathbf{V}_{\mathrm{and}} \mathbf{V}_{\mathrm{and}}^*} \right) \right)^{-1} 
\overline{\mathrm{diag} \left( \mathbf{V}_{\mathrm{and}} \mathbf{V} \mathbf{\Sigma}^* \mathbf{\Phi} \right)},
\end{equation}
where $\circ$ is the elementwise multiplication of two matrices, the overline denotes the complex-conjugate of a vector (matrix) and $\mathbf{V}_{\mathrm{and}}$ is the Vandermonde matrix defined as:
\begin{equation}
\mathbf{V}_{\mathrm{and}} = \begin{bmatrix} 
1 & \lambda_1 & \cdots & \lambda_1^{M} \\
1 & \lambda_2 & \cdots & \lambda_2^{M}   \\
\vdots & \vdots& \ddots & \vdots \\
1 & \lambda_r & \ddots & \lambda_r^{M}   \\
\end{bmatrix} \in \mathbb{R}^{r \times M}.
\end{equation} 
For the rest of the work, we refer to this approach as optimal DMD. The DMD approximation of the snapshot matrix $\mathbf{U}$ can be also computed as follows:
\begin{equation}\label{eq:Opt_DMD}
\mathbf{U}_\mathrm{DMD} = \mathbf{\Phi} \,\mathrm{diag}(\mathbf{b}) \, \mathbf{V}_{\mathrm{and}}.
\end{equation} 

$\mathbf{U}_{j,w}=[\mathbf{u}_{w_1},\mathbf{u}_{w_2},\dots,\mathbf{u}_{w_k}]$,
here $w_1,w_2,…,w_k$ represent the time indices in the sub-window. At each level $j$, for sub-window $w$, standard DMD is applied to $\mathbf{U}_{j,w}$, yielding:
\begin{equation}
\mathbf{u}(t) = \sum_{k=1}^{r_{max}} b^{j,w}_k \bm{\phi}^{j,w}_k e^{\omega^{j,w}_k t},
\end{equation}
where in the last summation the terms $b^{j,w}_k$, $\bm{\phi}^{j,w}_k$ and $\omega^{j,w}_k$ have similar meaning to the ones in \eqref{eq:DMD_expression}.

\section{Numerical examples}\label{sec:numres}

This section presents the results of applying the model order reduction techniques of Section \ref{sec:morsec} to three different numerical examples of wave propagation in materials with damage. The first example considers the wave equation \eqref{eq:ibvp_wave} in a two-dimensional domain, where the damage is modeled as a localized reduction in the wave speed coefficient. For this benchmark problem, we assume full access to the governing equations and simulation data, enabling a comprehensive comparison of both intrusive and non-intrusive MOR techniques. The second example involves a complex linear elasticity model of a composite material, specifically fiber metal laminate (FML), in which the damage is represented by a localized reduction in the Young’s modulus. Here, we apply the data-driven, non-intrusive MOR methods DMD and OpInf, reflecting a practical scenario where knowledge about the structure of the system is available as well as input-output data. The third example concerns GUW propagation in a nonlinear elastic model of aluminium. Here, we focus exclusively on data-driven DMD assuming only displacement data is available and the full order model is not accessible.\par 

In order to assess the accuracy of the reduced order models, we consider the relative $L^2$-error at the sensor location defined as:
\begin{equation}\label{eq:l2_error}
\epsilon_s = \frac{{\left\lVert \mathbf{u}_s - \mathbf{u}_s^* \right\rVert}_{L^2(0,T)}}{{\left\lVert \mathbf{u}_s \right\rVert}_{L^2(0,T)}} \times 100 \%,
\end{equation}
where $\mathbf{u}_s$ and $\mathbf{u}_s^*$ are time-dependent signals at the sensor for the full order model (FOM) and the reduced order model (ROM), respectively. We also compute the time-dependent relative $L^2$-error over the entire spatial domain:
\begin{equation}\label{eq:l2_error_t}
\epsilon_u (t) = \frac{{\left\lVert \mathbf{u}(t) - \mathbf{u}^*(t) \right\rVert}_{L^2(\Omega)}}{{\left\lVert \mathbf{u}(t) \right\rVert}_{L^2(\Omega)}} \times 100 \%,
\end{equation}
where $\mathbf{u}(t)$ and $\mathbf{u}^*(t)$ denote the spatial solutions of the FOM and ROM at time $t$, respectively. Moreover, we define the cumulative energy $c_r$ (for a reduced dimension $r$) of the singular values as
\begin{equation}\label{eq:cumenerg}
c_r = \frac{\sum_{k=1}^r \sigma_k}{\sum_{k=1}^M \sigma_k},
\end{equation}
where $\sigma_k$ is $k$-th element in the diagonal of the matrix of the singular values $\mathbf{\Sigma}$, compare \eqref{eq:svdeq}.

In the numerical experiments presented in this section, our objective is to accurately reconstruct the system dynamics over the same time interval used during training. This choice is motivated by the fact that, in the context of damage identification and characterization, time extrapolation is generally of limited practical relevance. In particular, extending the simulation to significantly longer time horizons typically provides little additional information regarding the nature of the damage, while increasing the complexity of the inference procedure and making the interpretation of damage-induced wave responses more challenging. Nevertheless, since the accuracy and stability of long-time extrapolation are topics of considerable interest in the reduced order modelling community, we include a dedicated numerical experiment addressing this aspect in Section~\eqref{sec:waveEqNum}.

\subsection{Numerical results for the wave equation}\label{sec:waveEqNum}

\begin{figure}[htbp]
  \centering
  \begin{subfigure}[b]{0.25\textwidth}
    \centering
    \includegraphics[width=\linewidth]{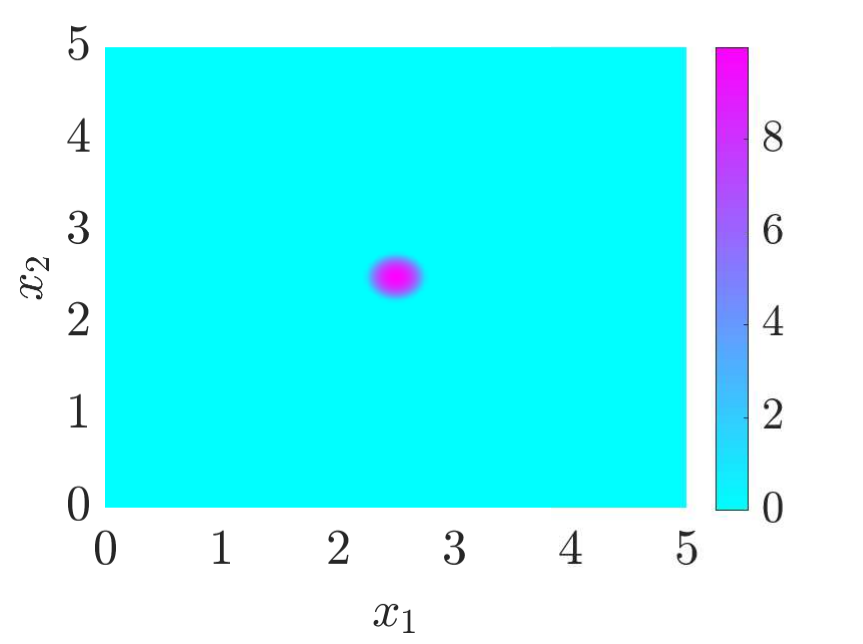}
    \caption{$t=0$}
    \label{fig:sub1}
  \end{subfigure}\hfill
  \begin{subfigure}[b]{0.25\textwidth}
    \centering
    \includegraphics[width=\linewidth]{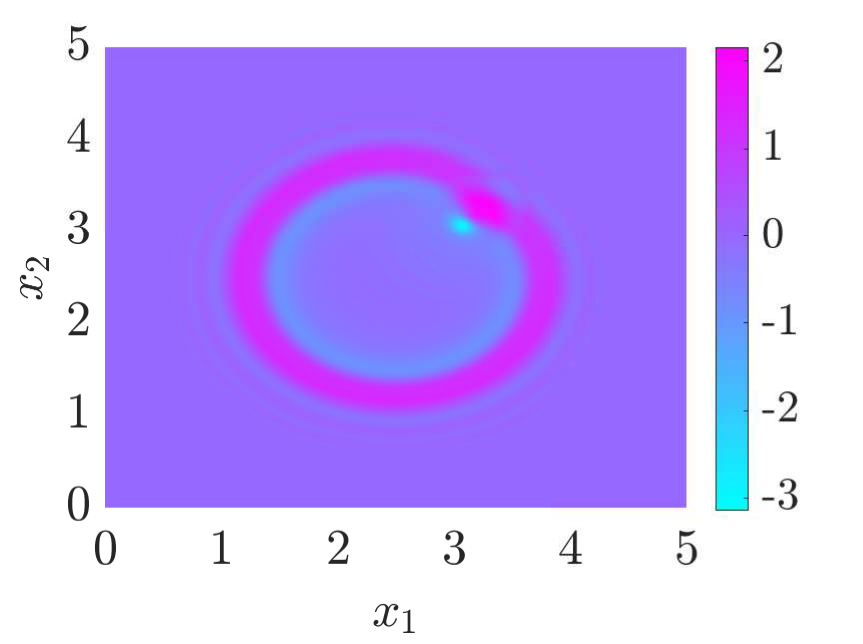}
    \caption{$t=2.5$}
    \label{fig:sub2}
  \end{subfigure}\hfill
  \begin{subfigure}[b]{0.25\textwidth}
    \centering
    \includegraphics[width=\linewidth]{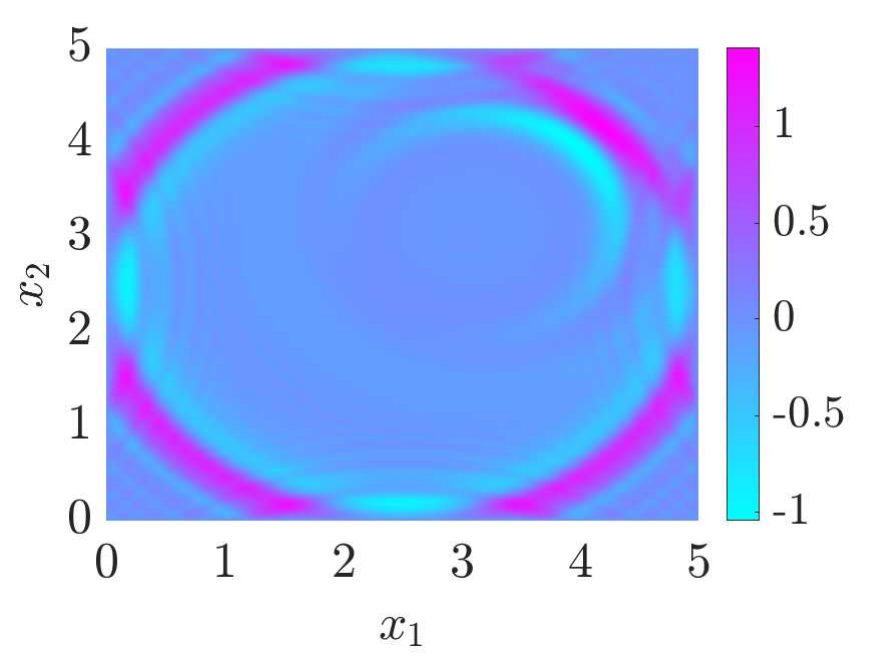}
    \caption{$t=5$}
    \label{fig:sub3}
  \end{subfigure}\hfill
  \begin{subfigure}[b]{0.25\textwidth}
    \centering
    \includegraphics[width=\linewidth]{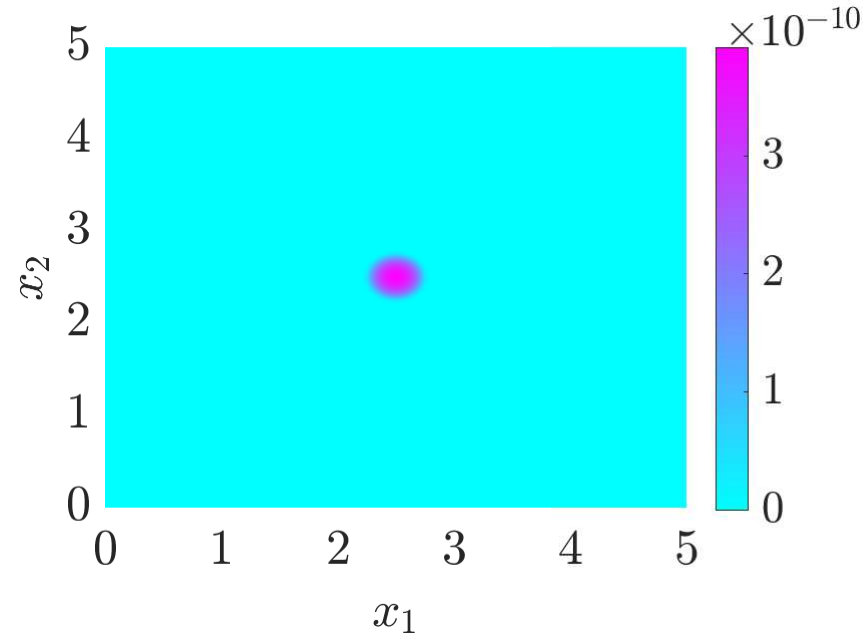}
    \caption{$\dot{\mathbf{u}}_0$}
    \label{fig:sub4}
  \end{subfigure}
  \caption{Displacement solution in [\si{m}] of the FOM for the wave equation at times $t=0, 2.5, 5$ and initial time derivative of the displacement field in [\si{m/s}]}
  \label{fig:2DwaveModel:highDimSolutionModel}
\end{figure}

In this subsection, we consider the problem of the wave equation \eqref{eq:ibvp_wave} with a parameterized damage and recall that this model is taken from \cite{Bur2024}, where the problem setting is based on \cite{NandaLorenzWave}. For the initial boundary value problem of the wave equation \eqref{eq:ibvp_wave} we make the following choices: the spatial domain is $\Omega = (0,5)\times (0,5) \subset \mathbb{R}^2$ and the time interval is set to $(t_0,t_{end})=(0.1,5)$. We set the right-hand side $f(t,\mathbf{x})=0$ and the boundary condition $b(t,\mathbf{x})=0$. The initial conditions are chosen as $u_0(\mathbf{x}) = u_{\mathbf{x}_M}(0.1) g(\mathbf{x})$ and $u_1(\mathbf{x}) = \frac{d u_{\mathbf{x}_M}}{dt}(0.1) g(\mathbf{x})$, where $g$ is defined as follows:
\begin{equation}\label{eq:g_wave_damaged}
g(\mathbf{x}) = \left\{\begin{aligned}
 & \exp \left( 1 + \frac{0.09}{ {(x_1 - x_c)}^2 + {(x_2 - y_c)}^2 - 0.09 } \right), && \quad \text{for } \sqrt{{(x_1 - x_c)}^2 + {(x_2 - y_c)}^2} < 0.3,  \\
 & 0, && \quad \text{for } \sqrt{{(x_1 - x_c)}^2 + {(x_2 - y_c)}^2} \geq 0.3,\\
\end{aligned}\right.
\end{equation}
with $\mathbf{x}=(x_1,x_2) \in \Omega, \quad \mathbf{x}_c = (x_c,y_c)=(2.5,2.5) \in \Omega$ and $u_{\mathbf{x}_M}$ is given by
\begin{equation}
u_{\mathbf{x}_M}(t) = 10 + \sin(2\pi A t), \quad A = 450000, \quad \text{for } t\leq 0.1.
\end{equation}
For $t \in [0,0.1]$, the displacement function $u(t,\mathbf{x})$ is defined as the product of $u_{\mathbf{x}_M}(t)$ and $g(\mathbf{x})$ and for $t > 0.1$ as the solution to \eqref{eq:ibvp_wave}. 
The reason behind splitting the time interval is to model a first phase where a wave motion is introduced in the material and a second phase where the interaction of the propagating wave with the material is observed.
The wave speed coefficient $c^2$ satisfies the following relationship:
\begin{equation}
c^2(\mathbf{x}) =  0.25 - (0.25 - c_D) \exp \left( - \frac{ {(x_1 - x_{1D})}^2 + {(x_2 - x_{2D})}^2 }{d_D} \right),
\end{equation} 
where $\mathbf{x}_D = (x_{1D} , x_{2D})$ is the location of the damage at which the value of the wave speed is reduced to $c_D < 0.25$ and $d_D$ is a parameter which describes the size of the damage. Piecewise linear finite elements with $h_{max}=0.05$ are used to discretize~\eqref{eq:ibvp_wave} in space leading to $N=11317$ spatial degrees of freedom and the $\theta$-method with $\theta=\frac{1}{4}$  is used to solve the discretized system in time starting from $t_0=0.1$ with a time step value of $\tau = 0.02$ which gives $M=245$ snapshots. % The initial fields of the displacement and its time derivative are shown in \eqref{fig:2DwaveModel:IC}.
%\begin{figure}[htbp]
%\centering
%\begin{subfigure}[t]{0.48\textwidth} 
%\includegraphics[width=0.9\linewidth]{Figures/2DWaveModel/u0.pdf}
%\includegraphics[width=0.7\textwidth]{fig/ex1/u0.pdf}
%\caption{$\mathbf{u}_0$}
%\end{subfigure}
%\begin{subfigure}[t]{0.48\textwidth} 
%\includegraphics[width=0.9\linewidth]{Figures/2DWaveModel/v0.pdf}
%\includegraphics[width=0.7\linewidth]{fig/ex1/v0.pdf}
%\caption{$\dot{\mathbf{u}}_0$}
%\end{subfigure}
%\vspace{10pt}
%  \caption{Initial displacement and initial time derivative of the displacement fields for the problem of the damaged wave equation}\label{fig:2DwaveModel:IC} 
%\end{figure}
Figure \ref{fig:2DwaveModel:highDimSolutionModel} shows the resulting numerical solution for the displacement of the finite element model at the times $t = 0, 2.5,5$ together with the initial time derivative of the displacement field on a domain with the damage at position $x_D = (3,3) \in \Omega$ with $d_D = 0.06$ and a reduction of the wave speed to $c_D = 0.05$. The wave motion is excited around the center $\mathbf{x}_c= (2.5,2.5)\in \Omega$ and propagates outwards from there. The excited wave is reflected at the damage in position $\mathbf{x}_D $ resulting in a second overlaid wave motion centering around the damage.

%\begin{figure}[htbp]
%\centering
%\begin{subfigure}[t]{0.325\textwidth} 
%\includegraphics[width=0.98\linewidth]{Figures/2DWaveModel/uFOMt0.pdf} 
%\includegraphics[width=0.6\textwidth]{fig/ex1/u0.pdf}
%\caption{$t= 0$}
%\end{subfigure}
%\begin{subfigure}[t]{0.325\textwidth} 
%\includegraphics[width=0.98\linewidth]{Figures/2DWaveModel/uFOMt25.pdf}
%\includegraphics[width=0.6\textwidth]{fig/ex1/u25.pdf}
%\caption{$t=2.5$}
%\end{subfigure}
%\begin{subfigure}[t]{0.325\textwidth} 
%\includegraphics[width=0.98\linewidth]{Figures/2DWaveModel/uFOMt5.pdf}
%\includegraphics[width=0.6\textwidth]{fig/ex1/u50.pdf}
%\caption{$t=5$}
%\end{subfigure}
%\begin{subfigure}[t]{0.325\textwidth} 
%\includegraphics[width=0.98\linewidth]{Figures/2DWaveModel/uFOMt5.pdf}
%\includegraphics[width=0.75\textwidth]{fig/ex1/v0.pdf}
%\caption{$\dot{\mathbf{u}}_0$}
%\end{subfigure}
%\vspace{10pt}
%  \caption{Solution of the FOM for the wave equation at times $t = 0, 2.5,5.$}\label{fig:2DwaveModel:highDimSolutionModel} 
%\end{figure} 

\begin{figure}[htbp]
  \centering
   \includegraphics[width=0.33\linewidth]{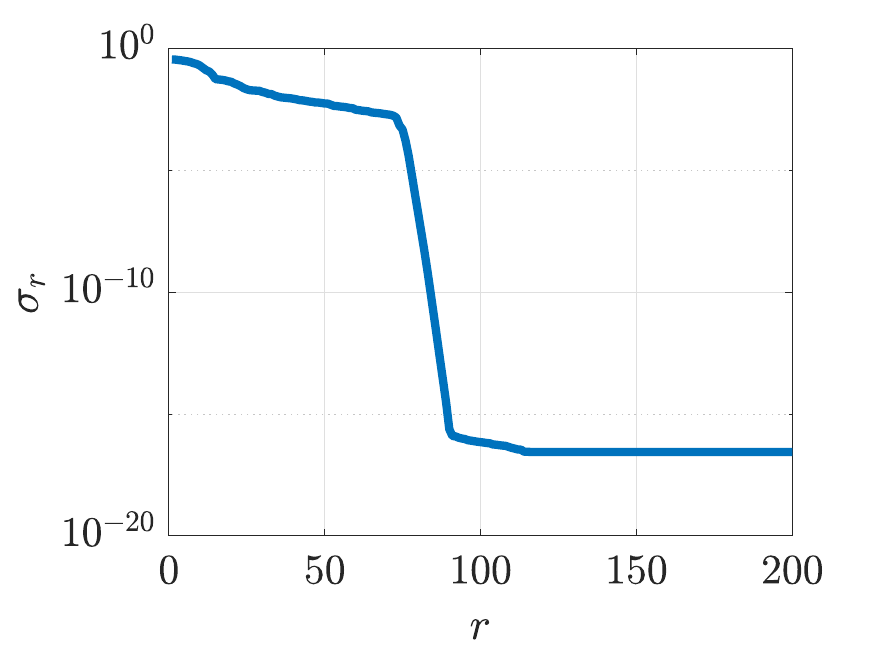}
   \includegraphics[width=0.33\linewidth]{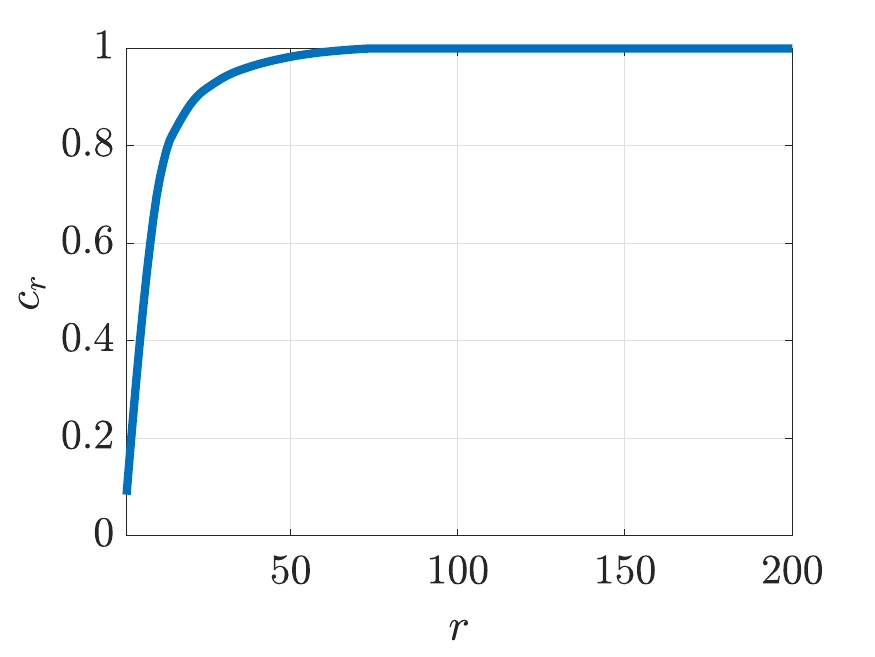}
  \caption{Decay of the singular values of the snapshot matrix $\mathbf{U}$ (left) and cumulative energy $c_r$ (right) for the wave equation.}
\label{fig:sv_decay_damagedWaveEq} 
\end{figure}

The decay of the singular values associated with the snapshot matrix $\mathbf{U}$ and the cumulative energy $c_r$ \eqref{eq:cumenerg} of the singular values are shown in Fig.~\ref{fig:sv_decay_damagedWaveEq}. The cumulative energy $c_r$ surpasses $95\%$ with $r = 33$, while the $99\%$ mark is reached for $r = 57.$  
In the following, we compare the solutions to the reduced order models constructed using POD, DMD and OpInf in terms of accuracy and efficiency. The implementation details are given as follows: The DMD initial amplitudes are computed using two options: (i) projection of the first snapshot onto the POD basis (see \eqref{eq:amplitudes_using_projectedPOD}), and (ii) an optimized amplitude formulation based on a least-squares fit over the time series (see \eqref{eq:amplitudes_optimal}). We refer to these two versions of DMD by 'DMD' and 'DMD Optimal', respectively. We mention here that the DMD solution is obtained using \eqref{eq:dmdadvance_red}.\par

Regarding operator inference, we apply both continuous and discrete versions of OpInf to this problem. We describe first the procedure for the time continuous setting. In this problem we have full knowledge about the external forces which are assumed to be identically zero for all time instants. This corresponds to the forces-informed variant of OpInf. Along the solution data, the approach requires the data for the second order derivative of the displacement mentioned in \eqref{eq:data2}. The required acceleration (second time derivative) data is computed using an eighth order finite difference scheme, see \cite{Sharma2024}, chosen to ensure high accuracy in the OpInf optimization problem. In order to improve solution stability during the optimization, a Tikhonov regularization term, scaled by the Frobenius norm of the inferred operators, is added to the OpInf optimization problem \eqref{eq:opinf_forces_informed_pb_symDecompo}. The optimization problem solves for the auxiliary matrices $\mathbf{L}$ and $\mathbf{E}$ which are subsequently used to compute the reduced mass and stiffness matrices. The optimization is implemented using the TensorFlow library~\cite{TensorFlow}, with the Adam algorithm~\cite{Adam} used as the optimizer. To account for scale differences between the displacement and its second derivative, the scaling method described in Section~\ref{sec:Opinf3} is applied. It has to be noted that in this example the zero operator $\mathbf{\widehat{P}} = \mathbf{0}$ is a trivial solution to the minimization problem because of the fact that the right hand side is identically zero. Therefore, a suitable initialization of the auxiliary matrices has to be done in order to avoid the trivial solution.\par

We shift the discussion to the time discrete version of OpInf which we denote by dOpInf. We aim at learning the reduced order operators $\mathbf{\widehat{H}}$ and $\mathbf{\widehat{R}}$ in \eqref{eq:dataDiscreteOpInf_PandD}. We remark that $\mathbf{\widehat{C}} = 0$ because of the absence of external input signals or loads. In this numerical example we have applied dOpInf without and also with the use of the re-projection method in \eqref{sec:Opinf5}. For both cases of dOpInf, we have generated snapshots which correspond to $q=6$ different initial data. In particular we sampled uniformly $x_c$ in \eqref{eq:g_wave_damaged} in the range of $[2.5,3.5]$ and we fixed $y_c = 2.5$. The generated snapshots are then concatenated to give the following matrices:
\begin{equation}
\mathbf{U}_0 = [\mathbf{U}_0^1, \mathbf{U}_0^2, \ldots, \mathbf{U}_0^q], \quad \mathbf{U}_1 = [\mathbf{U}_1^1, \mathbf{U}_1^2, \ldots, \mathbf{U}_1^q], \quad \mathbf{U}_2 = [\mathbf{U}_2^1, \mathbf{U}_2^2, \ldots, \mathbf{U}_2^q].
\end{equation} 
We then apply the re-projection method in \autoref{alg:reproj_disc} to each initial data which gives the re-projected trajectories $\mathbf{\bar{U}}_0^1, \dots \mathbf{\bar{U}}_0^q$, $\mathbf{\bar{U}}_1^1, \dots \mathbf{\bar{U}}_1^q$ and $\mathbf{\bar{U}}_2^1, \dots \mathbf{\bar{U}}_2^q$, concatenating the re-projected trajectories gives:
\begin{equation}
\mathbf{\bar{U}}_0 = [\mathbf{\bar{U}}_0^1, \mathbf{\bar{U}}_0^2, \ldots, \mathbf{\bar{U}}_0^q], \quad \mathbf{\bar{U}}_1 = [\mathbf{\bar{U}}_1^1, \mathbf{\bar{U}}_1^2, \ldots, \mathbf{\bar{U}}_1^q], \quad \mathbf{\bar{U}}_2 = [\mathbf{\bar{U}}_2^1, \mathbf{\bar{U}}_2^2, \ldots, \mathbf{\bar{U}}_2^q].
\end{equation} 
Then, the data matrices $\mathbf{\widehat{D}}$ and $\mathbf{\bar{D}}$ are assembled as in \eqref{eq:dataDiscreteOpInf_PandD} and \eqref{eq:opinf_reproj_disc} for the dOpInf without and with the re-projection method, respectively. In this setting, the matrices $\mathbf{\widehat{D}}$ and $\mathbf{\bar{D}}$ are both full rank which implies that the minimization problems correspond to solving two sets of matrix equations. The accuracy of the learned operators depends highly on the condition number of the matrices $\mathbf{\widehat{D}} \mathbf{\widehat{D}}^T$ and $\mathbf{\bar{D}} \mathbf{\bar{D}}^T$. We solve the system of equations to obtain the learned operators $\mathbf{\widehat{P}}$ and $\mathbf{\bar{{P}}}$ by inverting the matrices $\mathbf{\widehat{D}} \mathbf{\widehat{D}}^T$ and $\mathbf{\bar{D}} \mathbf{\bar{D}}^T$, respectively. These computations were carried out using Matlab 2023b. Both dOpInf without and with re-projection are applied without regularization in this comparison in order to isolate the effect of the re-projection procedure.\par 

\begin{figure}[htbp]
\centering
\includegraphics[width=0.49\linewidth]{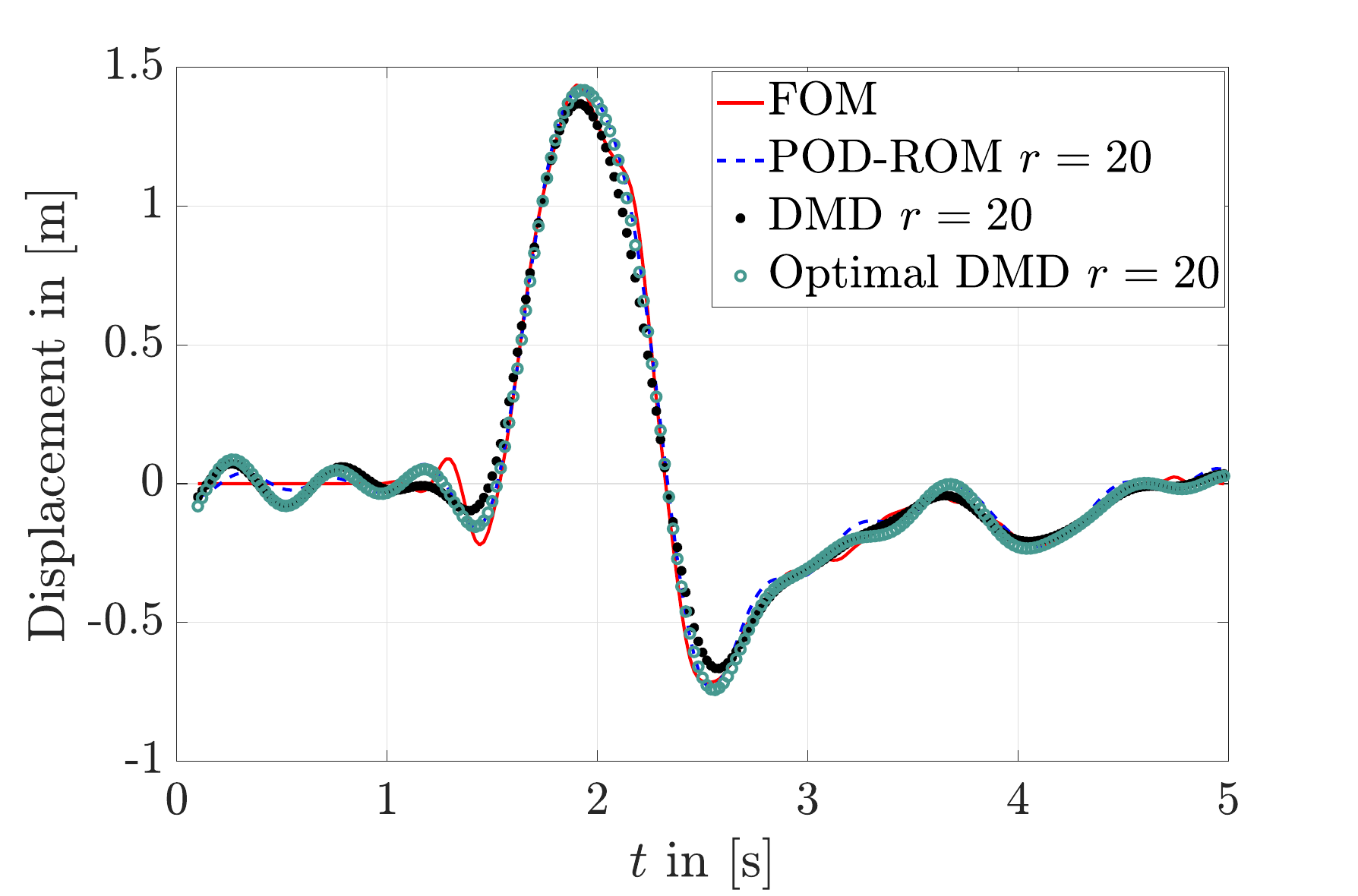}
\includegraphics[width=0.49\linewidth]{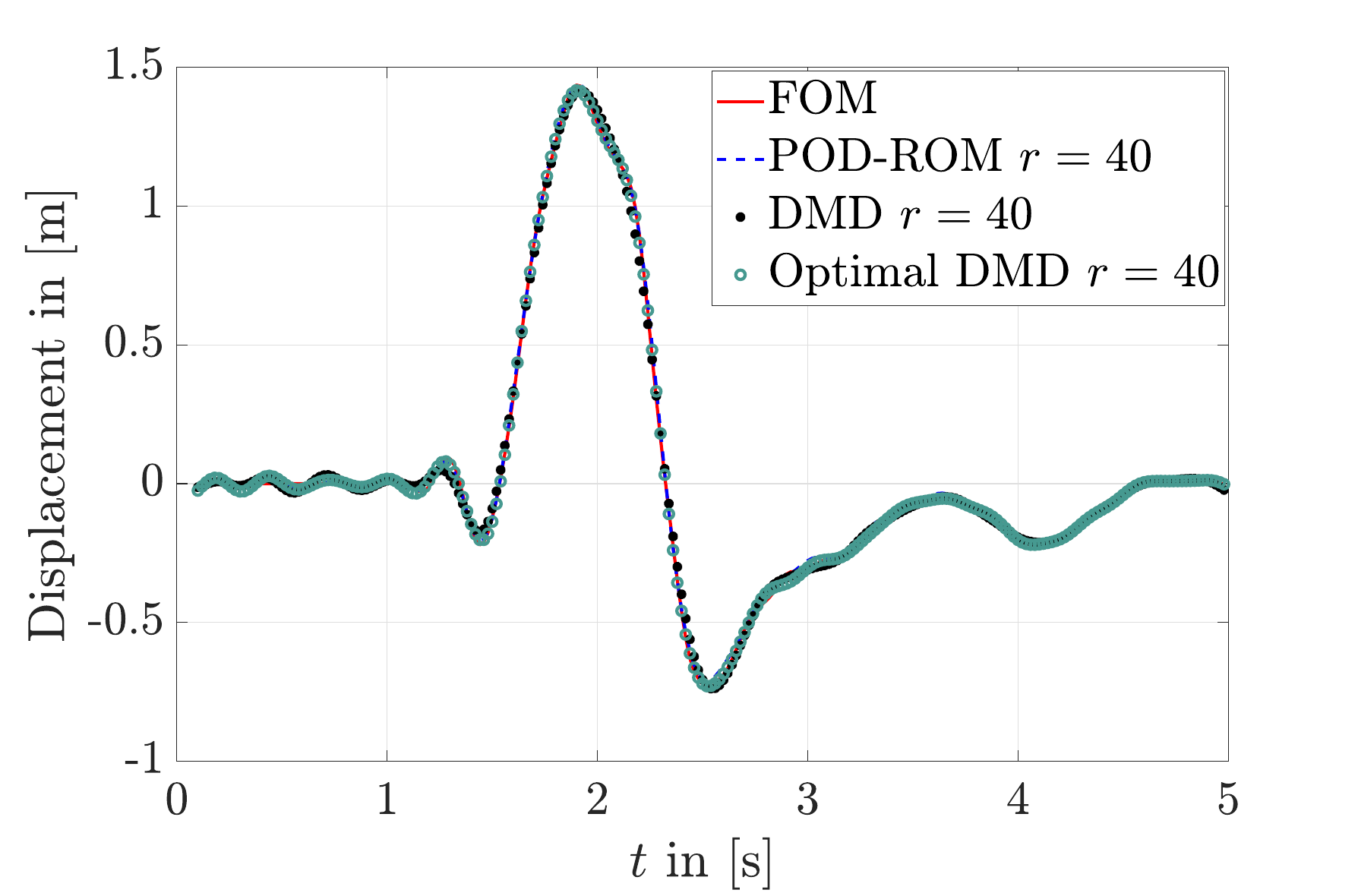}
\caption{Full order displacement versus the POD-ROM and DMD solutions at the sensor location for the reduced dimensions $r = 20$ (left) and $r = 40$ (right) for the wave equation, respectively.} 
\label{fig:2DwaveModel:displacement1}
\end{figure}

At this point, we present the results of the reduction performed by POD-ROM, DMD and OpInf. We start by visualizing the results obtained by POD-ROM and DMD versus the FOM. As an example, we consider the reduced dimensions of $20$ and $40$. The displacement solutions at the sensor location ($(3.5,2.5) \in \Omega$) for POD, DMD and optimal DMD are plotted in Fig.~\ref{fig:2DwaveModel:displacement1}. We notice that the POD-ROM is able to reproduce qualitatively the full order model signal. As for the quantitative assessment of the POD results, it yields a relative sensor error \eqref{eq:l2_error} $\epsilon_s = 8.872$ $\%$ and $3.021$ $\%$ for the reduced dimensions of $20$ and $40$, respectively.  As for DMD and optimal DMD, we observe that both versions have reconstructed the FOM solution successfully. The relative error values of DMD and optimal DMD for the case of $r=20$ are $13.371$ $\%$ and $9.671$ $\%$, respectively. As for the case of $r=40$, the relative error values of both DMD versions are $5.137$ $\%$ and $2.702$ $\%$, respectively. The results here indicate a more accurate reconstruction of the FOM signal is obtained via the optimal DMD variant with respect to the classical one. In addition we notice that optimal DMD and POD-ROM are close to each other in terms of accuracy.\par 

\begin{figure}[htbp]
\centering
\includegraphics[width=0.49\textwidth]{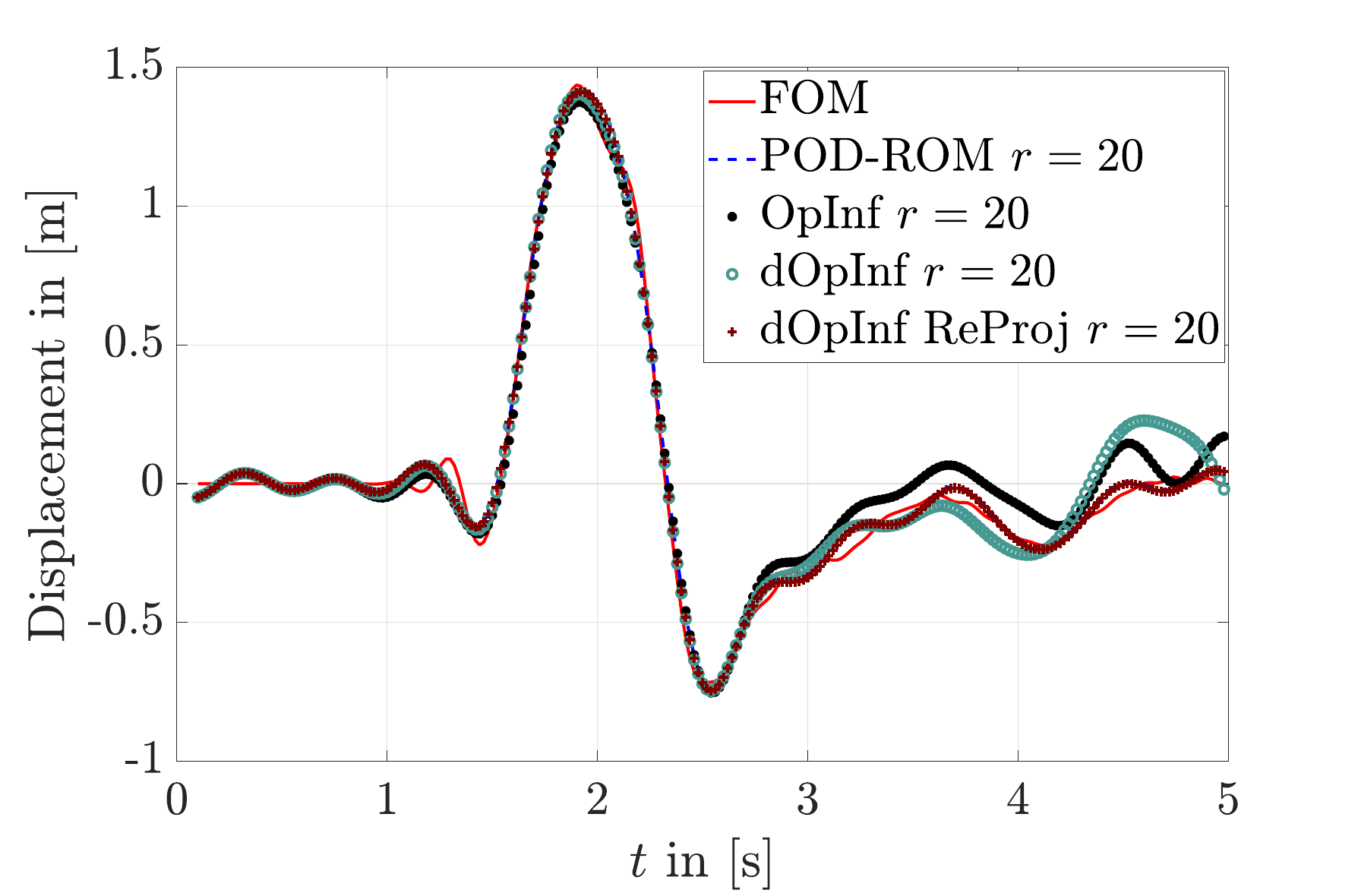}
\includegraphics[width=0.49\textwidth]{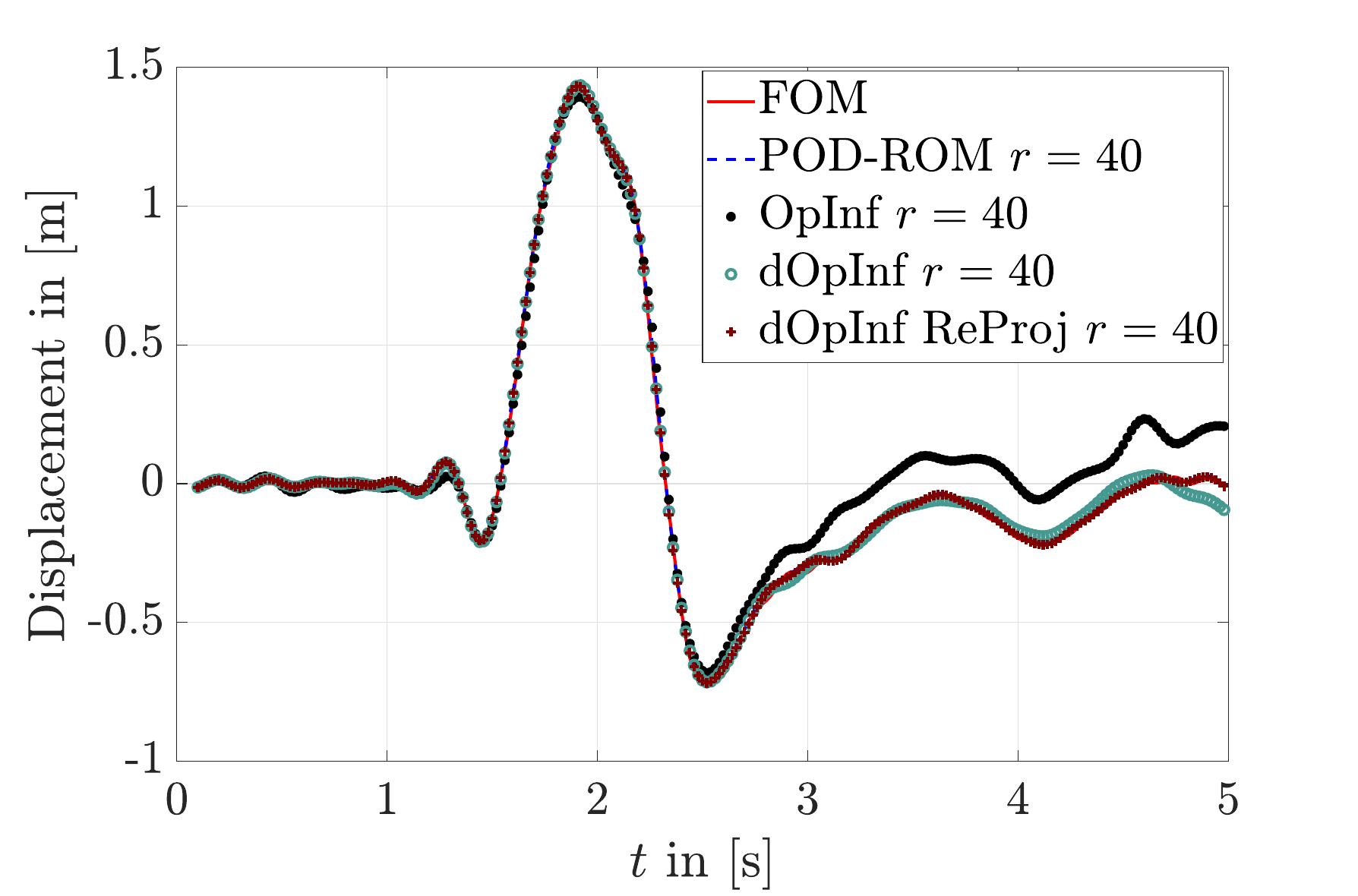}
\caption{Full order displacement versus the POD-ROM, OpInf and dOpInf solutions at the sensor location for the reduced dimensions $r = 20$ (left) and $r = 40$ (right) for the wave equation, respectively.} 
\label{fig:2DwaveModel:displacement2}
\end{figure}

We now present the results obtained using operator inference. In particular, we compare the performance of time continuous and time discrete operator inference against the POD-ROM and the FOM. For consistency, the same displacement signal shown in Fig.~\ref{fig:2DwaveModel:displacement1} is considered for all approaches, and the corresponding comparisons are shown in Fig.~\ref{fig:2DwaveModel:displacement2}. The results indicate that the time continuous operator inference approach is not sufficiently accurate in reconstructing the full order model signal, as a clear deviation is observed during the second half of the time history for both reduced dimensions $r=20$ and $r=40$. In contrast, the final plot shows that the time discrete operator inference approach achieves a more accurate quantitative reconstruction of the FOM signal, both with and without the use of re-projection.\par

To provide a quantitative assessment, time continuous OpInf yields relative sensor errors of $\epsilon_s = 16.812$ $\%$ and $24.943$ $\%$ for $r=20$ and $r=40$, respectively. For time discrete OpInf without re-projection, the corresponding sensor errors are $16.884$ $\%$ and $4.832$ $\%$. Finally, time discrete OpInf with re-projection achieves the highest accuracy among all tested reduced order models, with sensor errors of $8.532$ $\%$ and $2.4$ $\%$ for $r=20$ and $r=40$, respectively.\par

These results demonstrate that time discrete operator inference with re-projection is a highly accurate reduction methodology for second order mechanical systems. The improved performance relative to the case without re-projection can be attributed to a more accurate recovery of POD-intrusive operators enabled by the re-projection procedure. In fact this can be seen by computing the relative error between the learned operator in case dOpInf with and without re-projection as follows:
\begin{equation}
\epsilon_H
=
\frac{
\left\lVert
\operatorname{vec}(\mathbf{H}_r)
-
\operatorname{vec}(\mathbf{H}^\star)
\right\rVert_{L^2(0,T)}
}{
\left\lVert
\operatorname{vec}(\mathbf{H}_r)
\right\rVert_{L^2(0,T)}
}
\times 100\%, \quad
\epsilon_R
=
\frac{
\left\lVert
\operatorname{vec}(\mathbf{R}_r)
-
\operatorname{vec}(\mathbf{R}^\star)
\right\rVert_{L^2(0,T)}
}{
\left\lVert
\operatorname{vec}(\mathbf{R}_r)
\right\rVert_{L^2(0,T)}
}
\times 100\%,
\end{equation}
where $\operatorname{vec}$ is the vectorization operator and the operators $\mathbf{H}^\star$ and  $\mathbf{R}^\star$ are the learned operators of dOpInf either  $\mathbf{\bar{H}}, \, \mathbf{\bar{R}}$ or $\widehat{\mathbf{H}}, \,  \widehat{\mathbf{R}}$. In the current example, the maximum value of $\epsilon_H$ for dOpInf with re-projection is $1.201$ $\%$ at $r=100$, while the corresponding $\epsilon_R$ remains below $1 \times 10^{-7}$ $\%$ for all reduced dimensions. In contrast, for dOpInf without re-projection, the lowest value of $\epsilon_H$ is $3.956$ $\%$ at $r=20$, whereas the highest value reaches $25.118$ $\%$ at $r=90$. For the same model, the minimum value of $\epsilon_R$ is $7.779$ $\%$ at $r=30$, while the maximum value is $39.102$ $\%$ at $r=90$.\par 

\begin{figure}[htbp]
\centering
\includegraphics[width=0.4\linewidth]{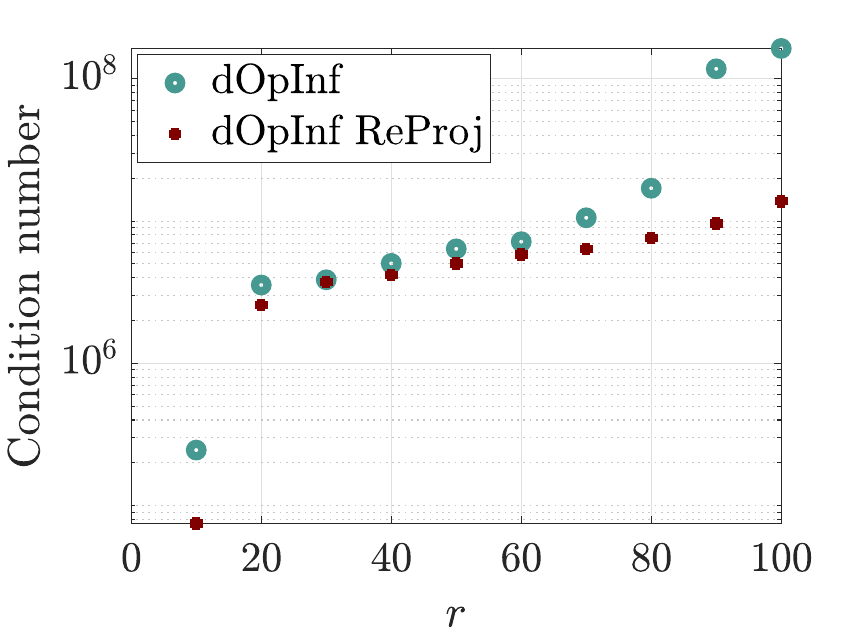}
\caption{Condition number for the matrices $\widehat{\mathbf{D}} \widehat{\mathbf{D}}^T$ and $\mathbf{\bar{D}} \mathbf{\bar{D}}^T$ which are inverted to solve the dOpInf problem without and with re-projection, respectively.} 
\label{fig:condition_number}
\end{figure}

We plot the condition number of the matrices $\widehat{\mathbf{D}} \widehat{\mathbf{D}}^T$ and $\mathbf{\bar{D}} \mathbf{\bar{D}}^T$ versus the reduced dimension in Fig.~\ref{fig:condition_number}. As it can be seen from the last plot, the condition number of the re-projection case ($\mathbf{\bar{D}} \mathbf{\bar{D}}^T$) is always smaller than the one  without re-projection ($\widehat{\mathbf{D}} \widehat{\mathbf{D}}^T$).

\begin{figure}[htbp]
\centering
\includegraphics[width=0.4\linewidth]{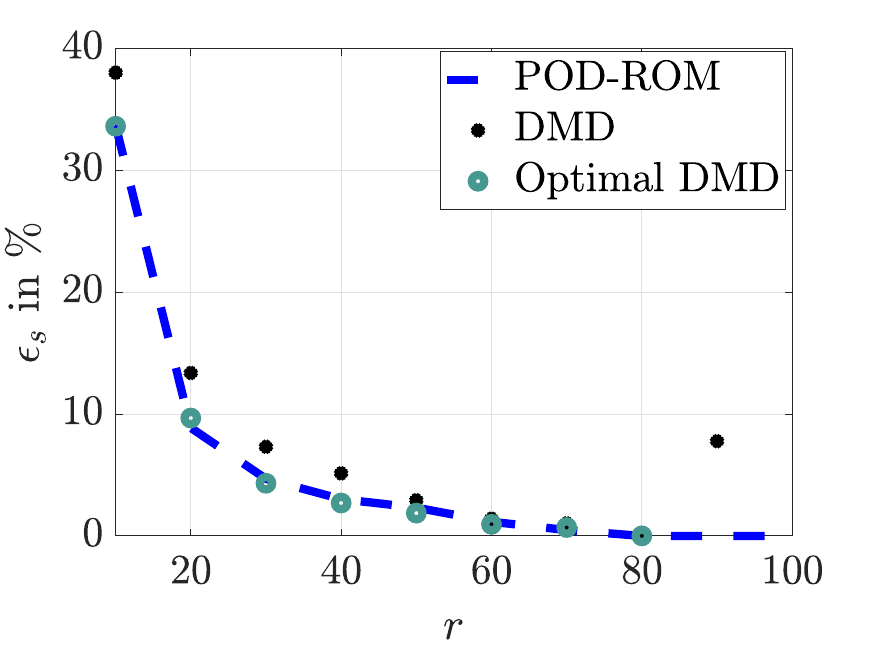}
\includegraphics[width=0.4\linewidth]{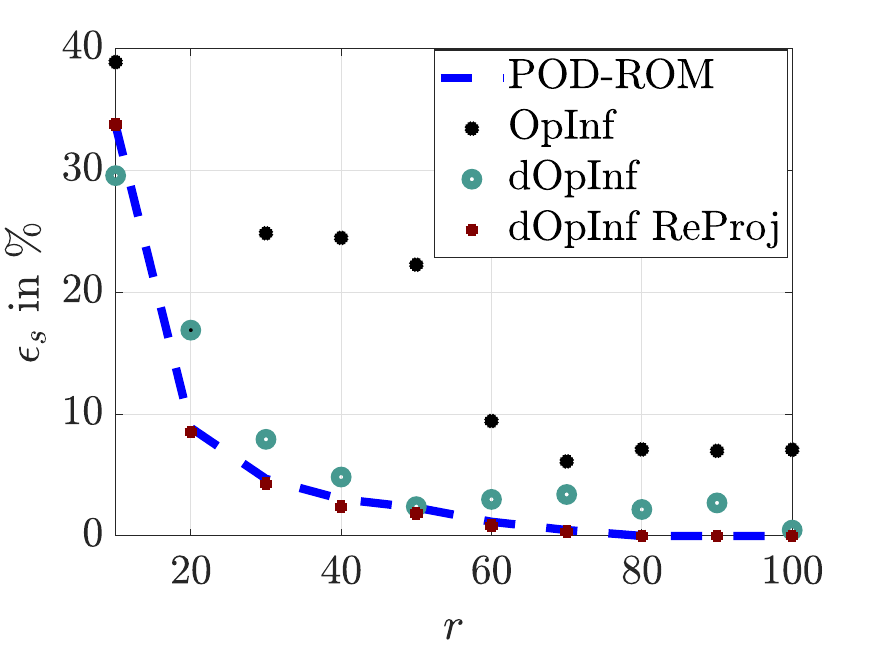}
\caption{Relative $L^2$-error $\epsilon_s$ at the sensor location as a function of the reduced dimension $r$ for POD-ROM compared with DMD variants (left) and POD-ROM compared with variants of OpInf (right).}
\label{fig:comparison_relerr_all}
\end{figure}

Next, we compare the decay of the relative error $\epsilon_s$ at the sensor location $(3.5,2.5) \in \Omega$ depending on the reduced dimension $r$ for the reduced order solutions generated by POD, DMD, Optimal DMD, OpInf and dOpInf in Fig.~\ref{fig:comparison_relerr_all}. For POD-ROM, we observe that the relative error $\epsilon_s$ decays rapidly and falls below $10 \%$ for $r = 20$ and below $5 \%$ for $r = 30$. DMD and Optimal DMD have similar convergence patterns with the latter being more accurate for all selected reduced dimensions, the Optimal DMD sensor error goes below the $3$ $\%$ threshold at $r=40$. As for OpInf ROM, the sensor error for $r=20$ and $r=40$ is significantly high, then the quality of the approximation of OpInf becomes higher with $r=60$ and then starts to saturate at $r=80$ with error value of $4.36$ $\%$. As for the dOpInf approaches, we notice that both perform better than their continuous counterpart. In particular, dOpInf without re-projection has achieved an error value of $2.4$ $\%$ for $r=50$, however, dOpInf with re-projection has resulted in low error values even lower than the ones obtained by POD-ROM for some choices of $r$. For example, for $r=50$ the dOpInf with re-projection has yielded an error of $1.84$ $\%$ (lower than the POD-ROM with $2.322$ $\%$) and it has given values lower than $1$ $\%$ for $r \geq 60$.

\begin{figure}[htbp]
\centering
\includegraphics[width=0.4\linewidth]{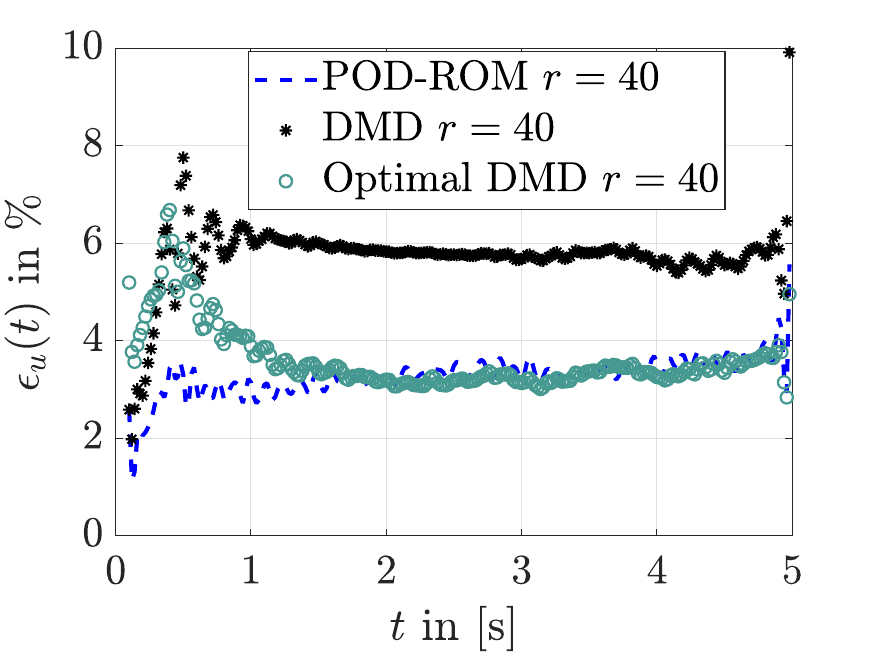}
\includegraphics[width=0.4\linewidth]{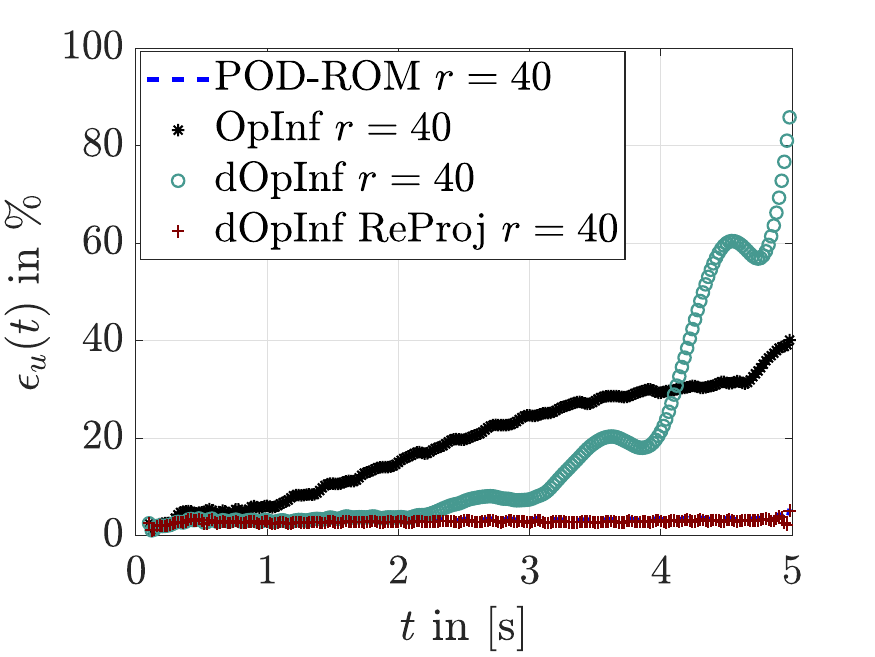}
\caption{Relative $L^2$-error $\epsilon_u(t)$ over the entire spatial domain as a function of time for the reduced dimension $r=40$ for reduced order solutions computed using POD and DMD variants (left), and OpInf and dOpInf variants (right) for the wave equation. The values are in percentages.}
\label{fig:comparison_epsilon_u_all}
\end{figure}

The time-dependent relative $L^2$-error $\epsilon_u(t)$ \eqref{eq:l2_error_t} over the entire spatial domain for the ROM of dimension $r = 40$ is plotted in Fig.~\ref{fig:comparison_epsilon_u_all}. For POD-ROM and Optimal DMD, most of the simulation the error stays between $3\%$ and $5\%$, while DMD gives slightly higher error values. In case of OpInf, the error value grows as time passes. This is attributed to the fact that the reduced inferred operators are not identical to the projected true ones which results in an error that accumulates throughout the simulation time. On the other hand, more accurate results are obtained using the discrete variants (dOpInf), where the one without re-projection is giving better results (than OpInf) for the first $4$ seconds of the simulation. The dOpInf with re-projection is the most accurate reduction methodology in this test with an error value even slightly better than the intrusive POD-ROM at certain time values. It is worth noting that dOpInf without re-projection exhibits a relatively large spatial relative error, reaching approximately $90$ $\%$ at the end of the simulation window. This behavior is associated with the unregularized formulation and can be mitigated by incorporating Tikhonov regularization. In the present comparison, however, regularization is not employed for the discrete time variants so that the effect of re-projection can be examined independently.

We compare the reduced order methods in terms of computational complexity. Table~\ref{tab:comptimes2dwave} summarizes the offline and online times for POD, DMD, OpInf and dOpInf variants for the reduced dimension of $r=40$. We would like to mention that the offline times include the FOM simulation time which is $4.086$ \si{s}. Clearly, OpInf has very large offline computational times compared to POD and DMD. However, the online computational times are comparable for all methods. The offline time needed for dOpInf with re-projection is obviously larger than dOpInf without re-projection because of the additional offline time required to query the FOM in the re-projection method.

\begin{table}[h!]
\centering
\small
\caption{Computational times of different reduced order methodologies for the wave equation with damage for $r=40$. The unit is [\si{s}].}
\label{tab:comptimes2dwave}

\begin{tabular}{m{7em} m{4.2em} m{3.2em} m{4.2em} m{4.2em} m{5.2em}}
\toprule
\diagbox[width=7em,height=4.8em]{Stage}{ROM Type}
& \textbf{POD-ROM}
& \textbf{DMD}
& \textbf{OpInf}
& \textbf{dOpInf}
& \shortstack{\textbf{dOpInf}\\\textbf{Re-projection}} \\
\midrule
\textbf{Offline} & 5.193 & 4.163 & 52.83 & 27.655 & 76.351 \\
\textbf{Online}  & 0.0149 & 0.013 & 0.0156 & 0.0108 & 0.0120 \\
\bottomrule
\end{tabular}
\end{table}

The final numerical experiment considered in this subsection addresses the time extrapolation task. The objective is to assess the ability of the reduced order models to accurately predict the system dynamics beyond the training interval. To this end, the reduced models are integrated over the time interval $[0,t_{\mathrm{final}}]$, where $t_{\mathrm{final}} > t_{\mathrm{end}}$, and their predictions are compared against the corresponding full order solution.

\begin{figure}[htbp]
\centering
\includegraphics[width=0.49\linewidth]{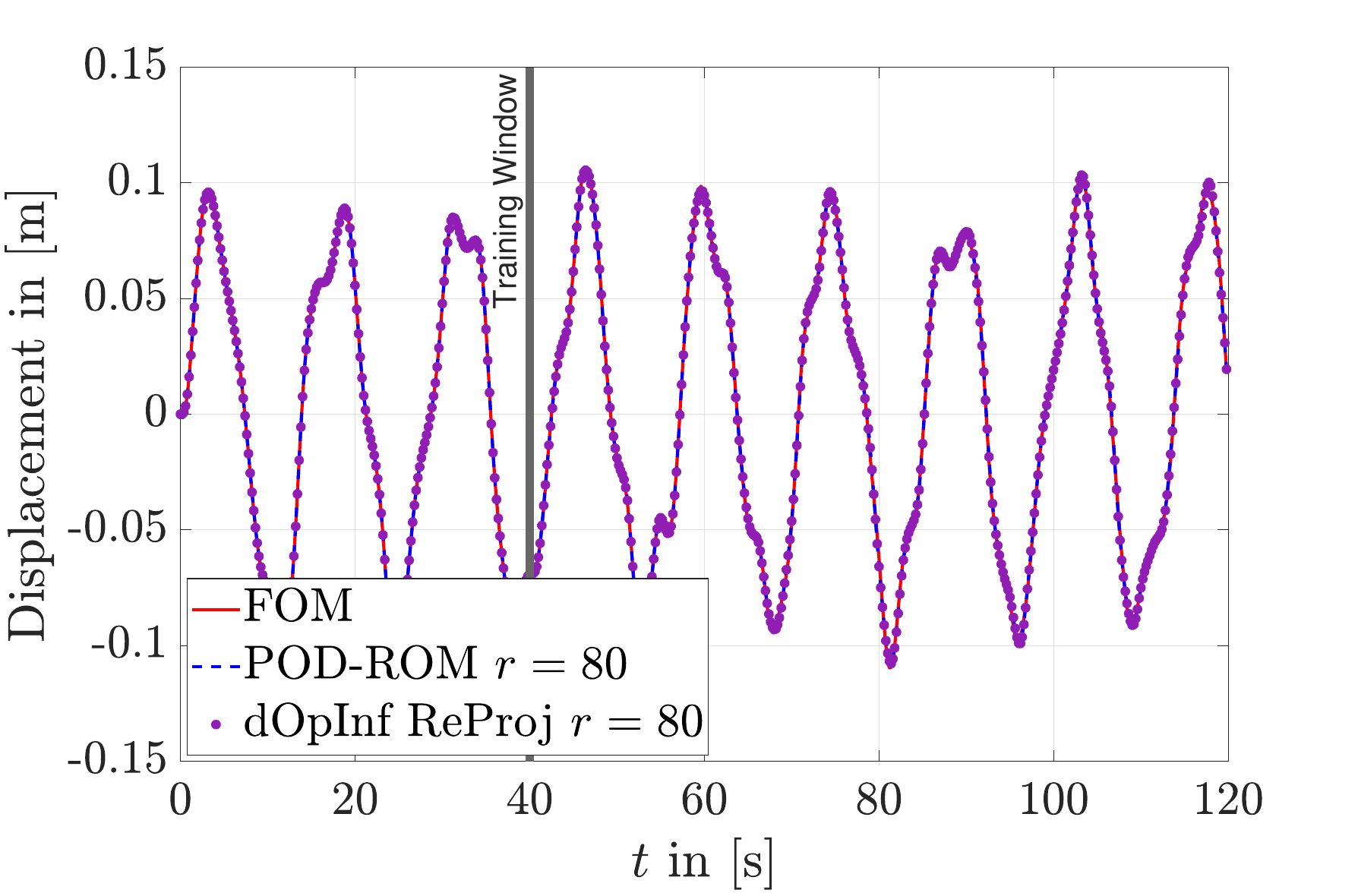}
\includegraphics[width=0.49\linewidth]{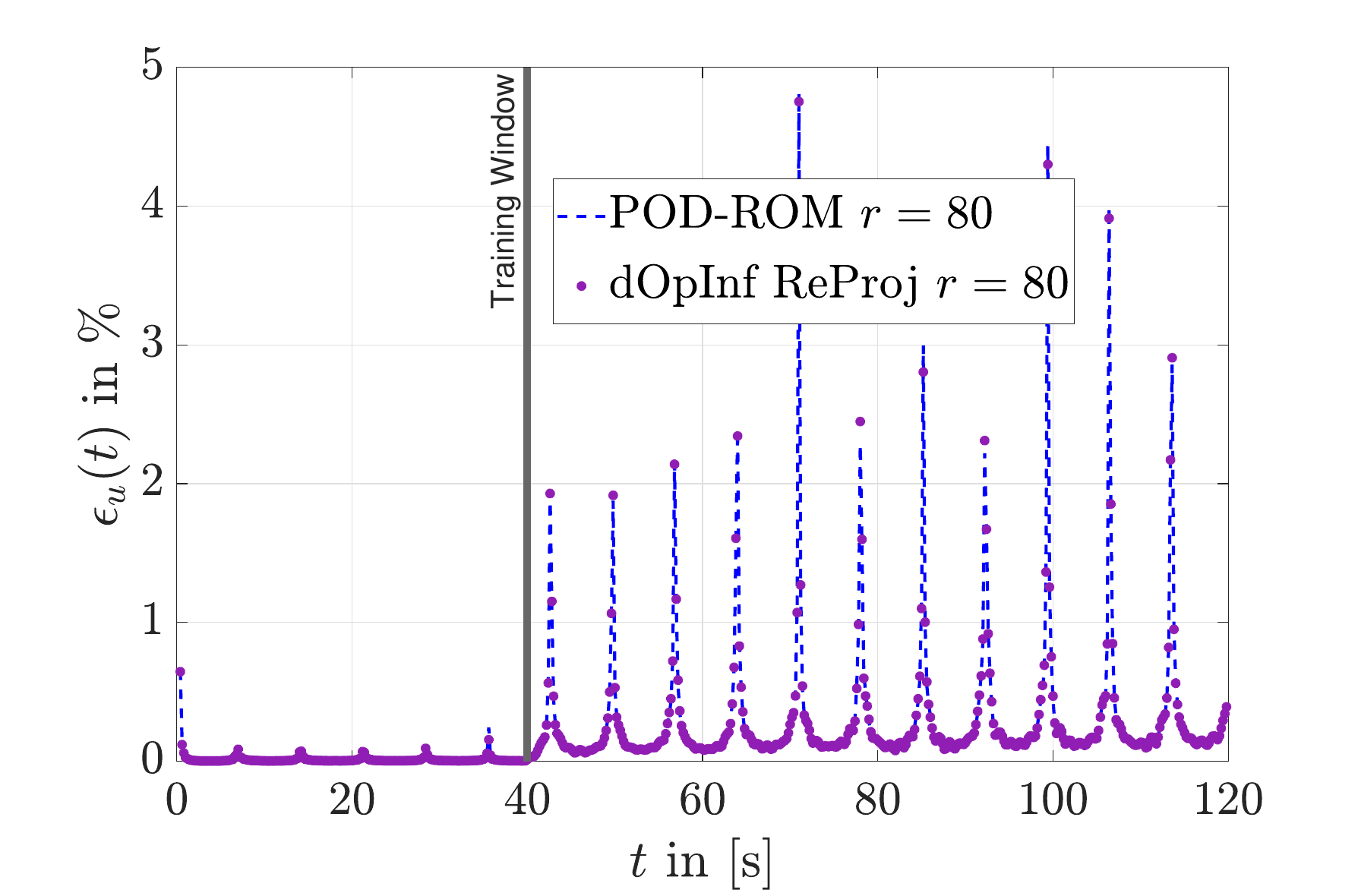}
\caption{Full order displacement versus the POD-ROM and dOpInf with reprojection solutions at the sensor location for the reduced dimensions $r = 80$, for the wave equation extrapolation test.} 
\label{fig:2DwaveModel:extrapolation}
\end{figure}

For this experiment, we consider the wave equation model introduced previously with the following modifications:

\begin{itemize}

\item The initial displacement and velocity fields are taken to be identically zero throughout the domain, i.e.,
\[
u_0(\mathbf{x}) = u_1(\mathbf{x}) = 0.
\]

\item The wave equation is driven by the following non-homogeneous forcing term:
\[
f(\mathbf{x},t)
=
\exp\left(
-\frac{(x_1-2.5)^2 + (x_2-2.5)^2}{4}
\right)
\exp(-0.1t)
\sin(0.4\pi t).
\]

\end{itemize}
The damage configuration remains unchanged throughout this experiment. We follow the same procedure as in the previous setting by collecting solution snapshots and projecting the governing equations onto the POD-generated reduced space. The training interval is defined by $t_{\mathrm{end}} = 40$~\si{s} with a time step of $0.2$~\si{s}. A total of $M=200$ snapshots are collected, and the reduced dimension is set to $r=80$. The POD-Galerkin reduced order model and the dOpInf model with reprojection are then simulated over the interval $[0,t_{\mathrm{final}}]$, where $t_{\mathrm{final}} = 120$~\si{s}, corresponding to a time horizon three times longer than the training interval. The results, evaluated at a representative spatial location in the domain are shown in Fig.~\ref{fig:2DwaveModel:extrapolation} (left) demonstrate that both the POD-ROM and the dOpInf model with reprojection maintain a high level of accuracy throughout the extrapolation period, without exhibiting noticeable instability or significant error growth. A more comprehensive assessment of the prediction accuracy is provided in the right subfigure of \eqref{fig:2DwaveModel:extrapolation}, which depicts the evolution of the spatial relative error as a function of time for both the POD-ROM and the dOpInf with reprojection. It can be observed that both reduced order models maintain low relative errors throughout the simulation, demonstrating their ability to accurately reproduce the full order dynamics over the entire spatial domain, including the extrapolation interval.

\subsection{Numerical results for the mechanical system resulting from GUW propagation}\label{sec:mechSys}

\begin{figure}[htbp]
    \centering
   \includegraphics[width=0.47\linewidth]{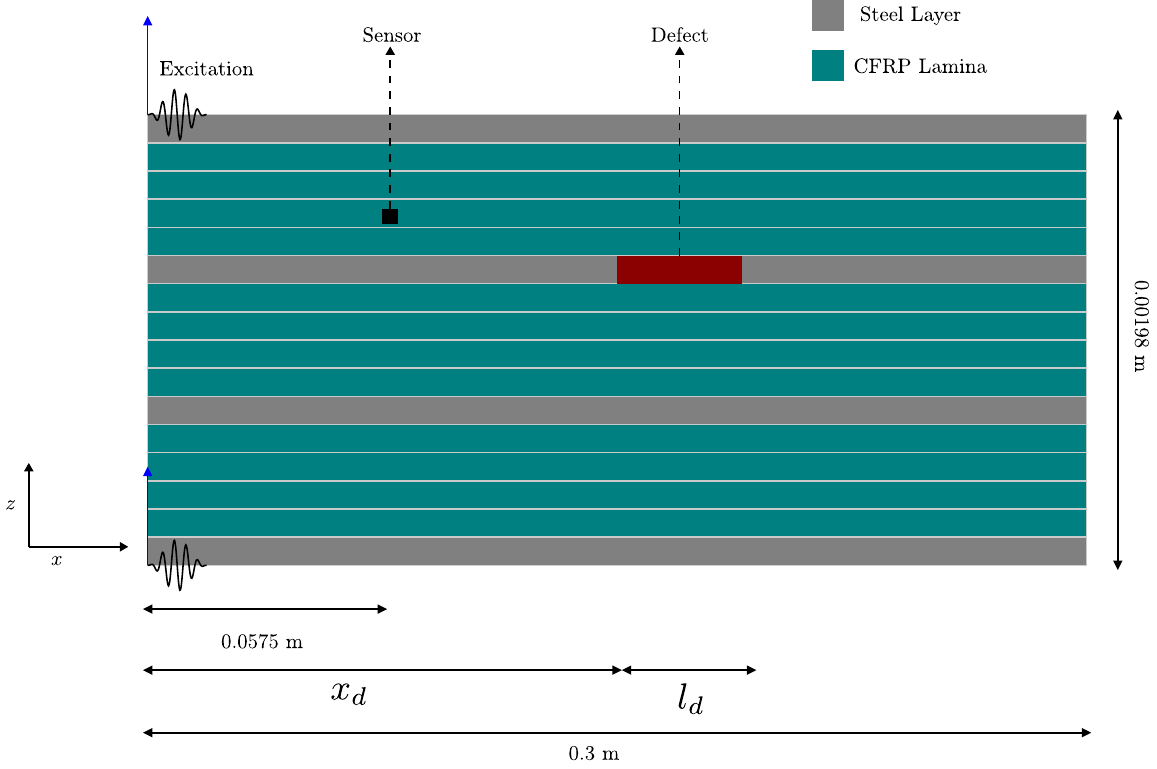} 
   \includegraphics[width=0.52\linewidth]{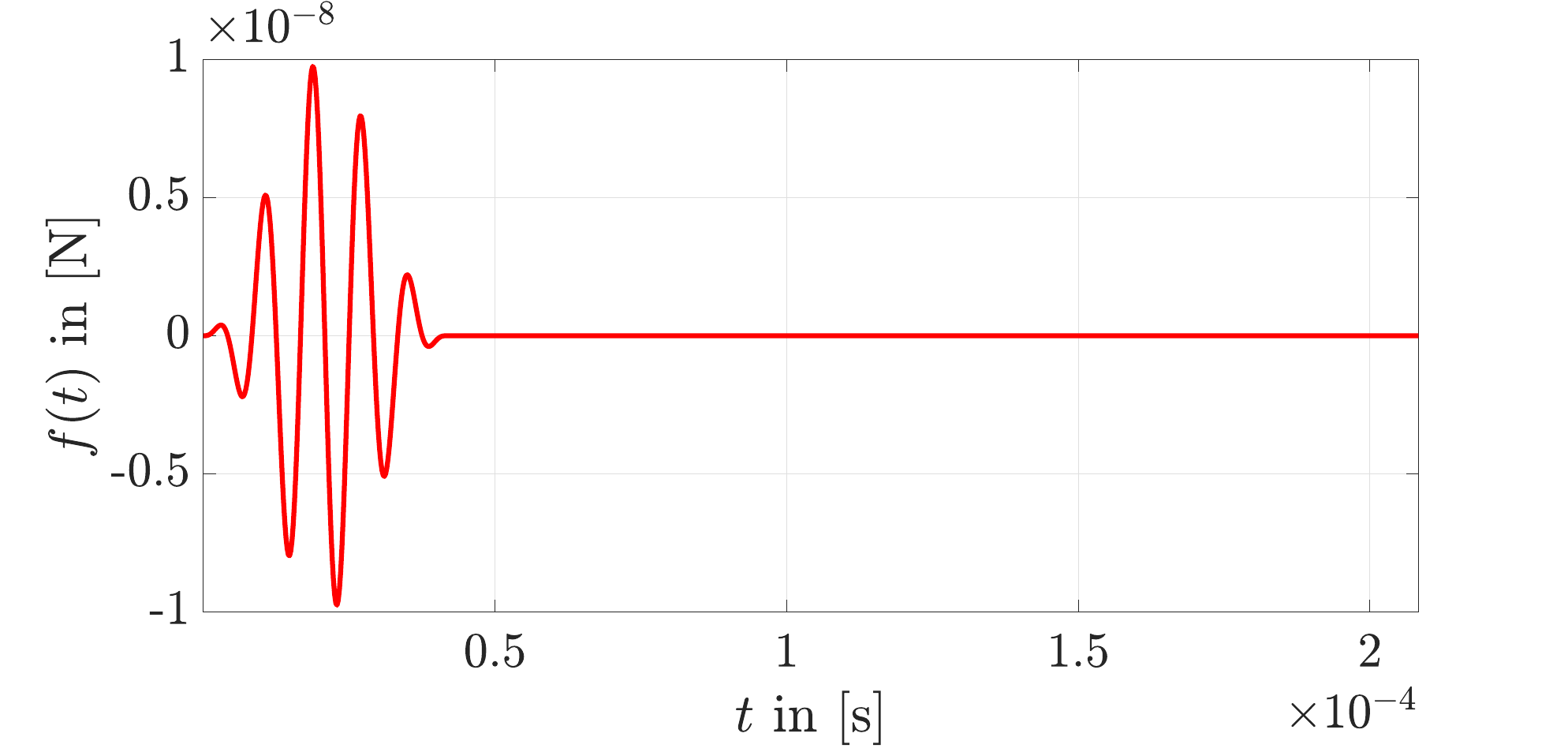} 
  \caption{FML model setup for the GUW simulation (left) and excitation force applied at top and bottom left nodes of the FML (right).}
\label{fig:FML_sketch} 
\end{figure}

In this subsection, we present the numerical results for the application of the non-intrusive model order reduction techniques operator inference (\ref{sec:Opinf1}, \ref{sec:Opinf2}, \ref{sec:Opinf3}) and dynamic mode decomposition (\ref{sec:DMD}) to the problem of the mechanical system \eqref{eq:mechanicalSystem_Input} obtained upon the discretization of the equation of motion \eqref{eq:eqofmotion} for elasticity problems. We consider the two-dimensional physical model of a fiber metal laminate (FML) of \cite{BellamMuralidhar2021}. The model contains one physical defect and describes the propagation of an anti-symmetric lamb wave mode ($A_0$) in the FML and the subsequent interaction with the defect. A sketch of the model is shown in Fig.~\ref{fig:FML_sketch} (left). The FML comprises sixteen alternating layers of carbon fiber reinforced plastic (CFRP) and steel. Damage is introduced by locally reducing the Young’s modulus $E_d$ in one of the steel layers. The simulation of the GUW propagation and its interaction with the damage is carried out using the commercial software COMSOL-Multiphysics software\textcopyright. The model is excited on the top and bottom left nodes of the model with a force excitation using a five-cycle Hanning window sinusoidal burst with a central frequency of $120$ kHz. The model is excited for $4.167 \times 10^{-5}$ \si{s} and the total time of the simulation is $2.083 \times 10^{-4}$ \si{s}. The thicknesses of steel and CFRP layers are $0.12$ \si{m} and $0.125$ \si{m}, respectively. The spatial discretization yields $N=79266$ degrees of freedom. The location of the damage along the FML is called $x_d$ and the length of the damage along the FML is denoted by $l_d$. The excitation force as function of time is shown in Fig.~\ref{fig:FML_sketch} (right). The time step size of the numerical simulation is $2.0833 \times 10^{-7}$, the solution is gathered at each time step leading to $M=1001$ snapshots. We present results for the damage configuration with the Young's modulus of elasticity $E_d = 2.55 \times 10^9$ \si{Pa}, the position of the damage $x_d = 0.06$ \si{m} and the length of the damage $l_d = 0.002$ \si{m}. A sensor is placed in the domain measuring the out-of-plane displacement (vertical displacement).

\begin{figure}[htp]
    \centering
       \includegraphics[width=0.55\linewidth]{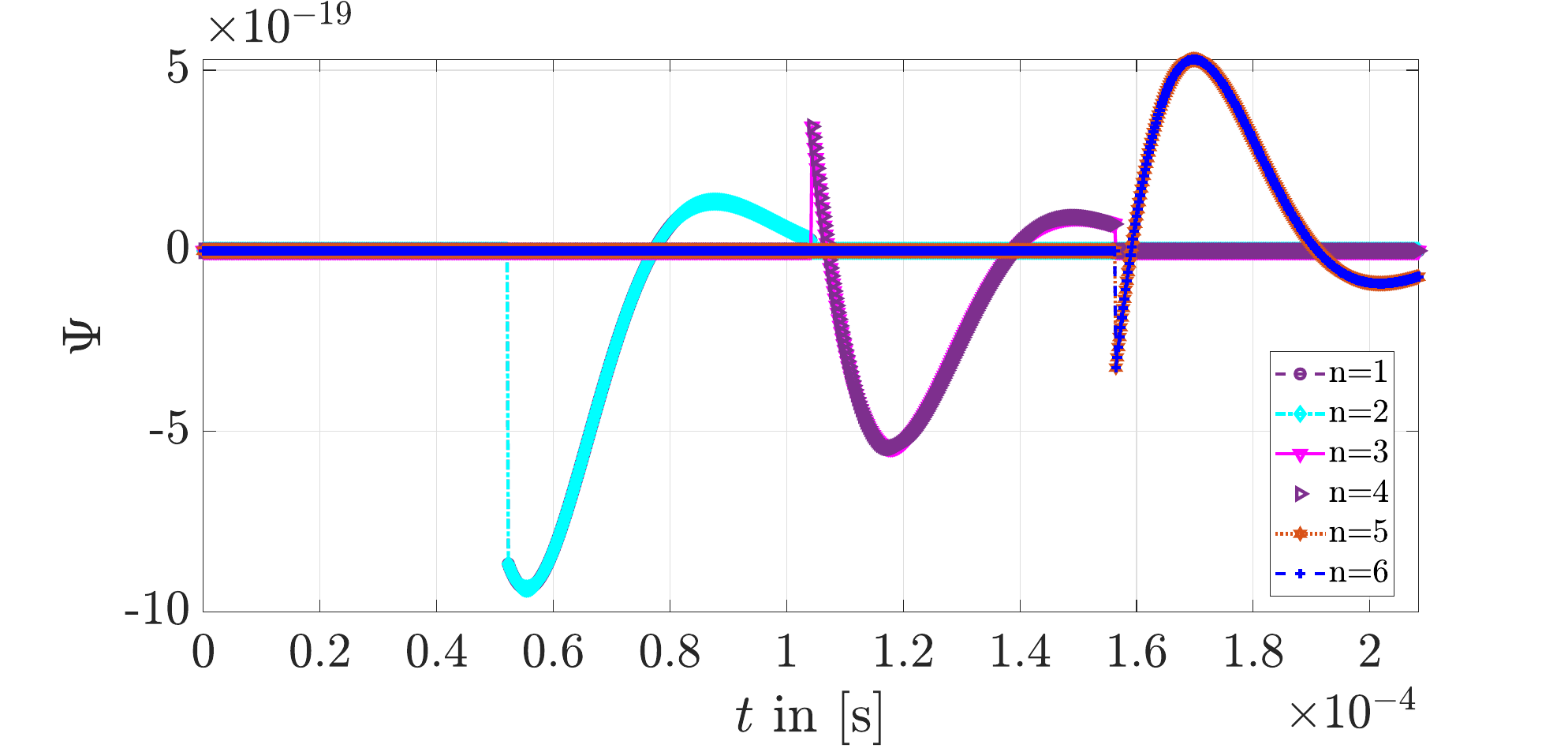} 
  \caption{Second level time dynamics obtained by mrDMD for the mechanical model of the damaged FML.}
\label{fig:mrdmd_2ndlevel_dynamics_damagedFML} 
\end{figure}

The access to the FOM is restricted in this numerical example. Therefore, we compare the DMD method and its multiresolution variant (mrDMD) to the operator inference approach. For DMD, the implementation details are given as follows. Like in the first example, we use two different strategies to compute the initial mode amplitudes: (i) projection of the first snapshot onto the POD basis (see \eqref{eq:amplitudes_using_projectedPOD}), and (ii) an optimized amplitude formulation based on a least-squares fit over the time series (see \eqref{eq:amplitudes_optimal}). In addition to the classical DMD, we use mrDMD which aggregates modal contributions from successive time-bins at different scales, ultimately producing a multiscale reconstruction of the original dynamics. For mrDMD, the maximum number of hierarchical time levels is set to $6$. This decomposition technique requires specifying a frequency threshold to distinguish between slow and fast dynamical modes. In this example, no slow modes are identified at the first level, which spans the full time window. At the second level, six significant modes are extracted. The temporal evolution of their coefficients, defined by $\mathbf{\Psi} = \mathbf{b} \, \exp(\Omega t)$, is illustrated in Fig.~\ref{fig:mrdmd_2ndlevel_dynamics_damagedFML} for the second level dynamics which contains $4$ different time-bins, the number of extracted modes here is $6$ corresponding to $3$ different pairs of frequencies.

\begin{figure}[htp]
  \centering
        \includegraphics[width=0.49\linewidth]{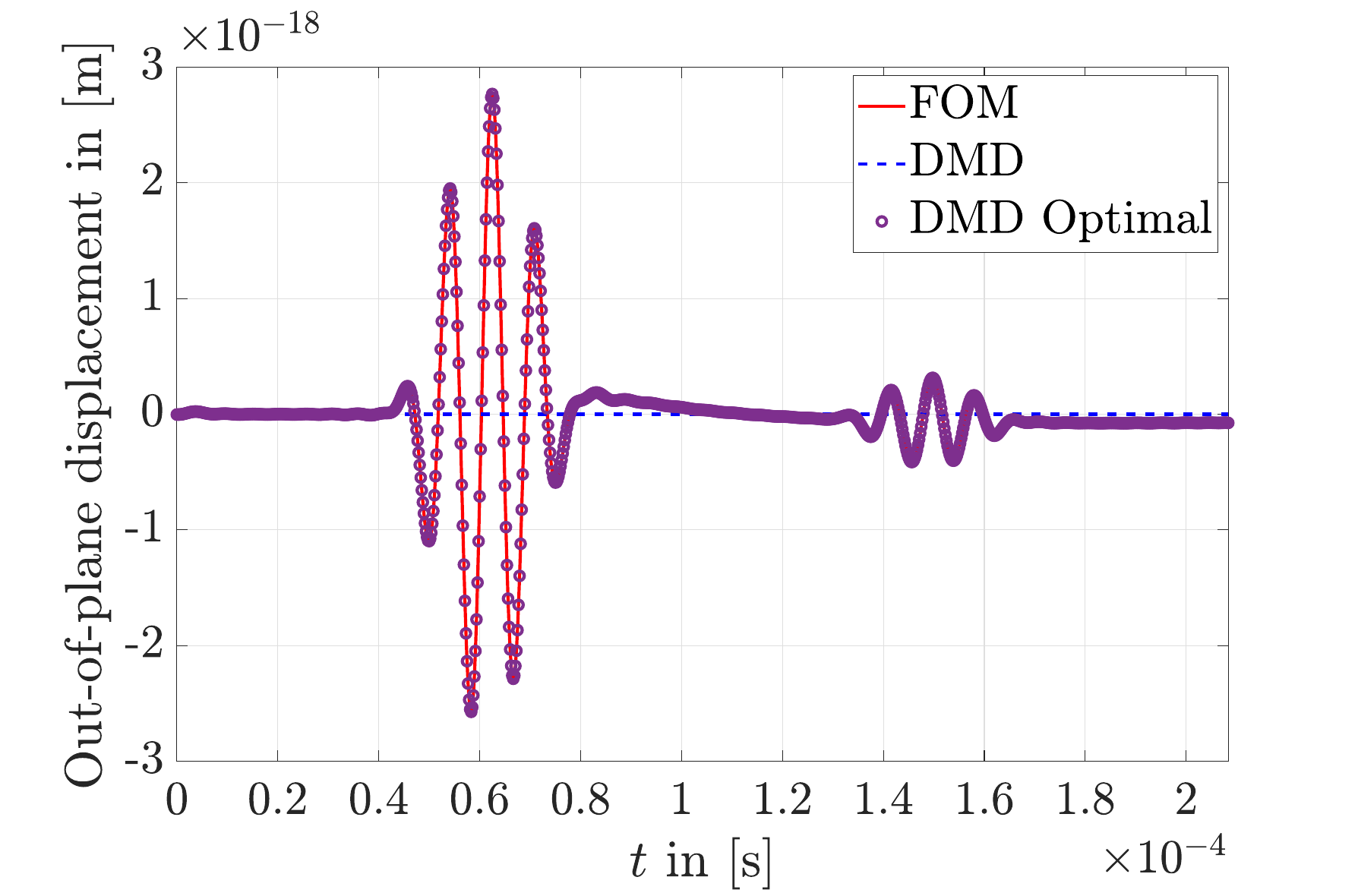} 
              \includegraphics[width=0.49\linewidth]{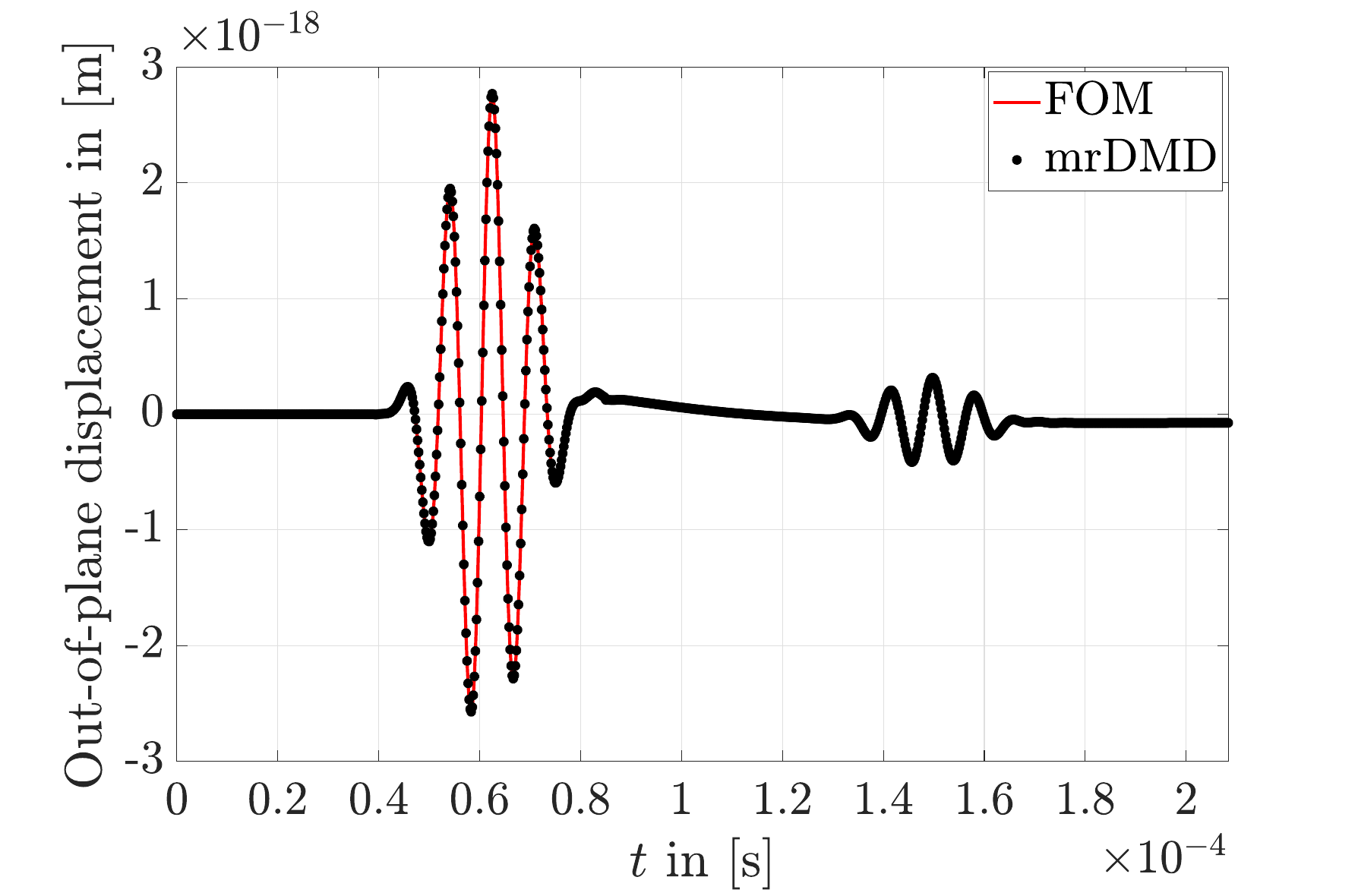} 
  \caption{Comparison between the FOM and DMD predicted displacement (left) and the mrDMD predicted displacement (right) at the sensor for the mechanical model of the damaged FML. The reduced dimension  is $r=80$ in the case of DMD and optimal DMD.}
\label{fig:FOM_vs_DMD_FMLModel_80modes} 
\end{figure}

Figure~\ref{fig:FOM_vs_DMD_FMLModel_80modes} (left) compares the full order model and DMD predictions at the sensor location. The DMD signal is obtained using \eqref{eq:dmdadvance}, whereas the optimal DMD is computed using \eqref{eq:Opt_DMD} and the reduced dimension is set to $r=80$. The DMD reconstruction using the first snapshot for amplitude estimation performed poorly, failing to recover the dynamics. This is attributed to the fact that the initial displacement field is nearly zero everywhere except at the excitation points, leading to negligible initial amplitudes. In contrast, DMD with optimized amplitudes accurately captured the signal, yielding a relative error \eqref{eq:l2_error} at the sensor of $\epsilon_s = 0.626 \%$.

The mrDMD reconstruction closely matched the FOM signal, see Fig.~\ref{fig:FOM_vs_DMD_FMLModel_80modes} (right), achieving a lower sensor error of $\epsilon_s = 0.137 \%$. This demonstrates that the multiresolution approach is effective at isolating time-localized dynamics in this guided wave simulation.

Finally, we compare the computational costs. DMD with optimized amplitudes and rank $r=80$ required $5.511$ \si{\second} for the computation of the DMD modes and eigenvalues, the construction of the Vandermode matrix and the resolution of the linear system for the optimal amplitudes vector $\mathbf{b}$. Then only $0.456$ \si{\second} are needed to compute the time dynamics $\mathbf{\Psi}$ and the reconstruction of the reduced solution. In the case of mrDMD, the method required $8.919$ \si{\second} for the computation of the DMD modes over the different time levels and time-bins inside each level, additional $3.218$ \si{\second} were taken for the extraction of the time dynamics and DMD modes from each time-bin and finally reconstructing the mrDMD solution. Although mrDMD offered improved accuracy, the optimal DMD was computationally more efficient.

\begin{figure}[htbp]
    \centering
    \includegraphics[width=0.33\linewidth]{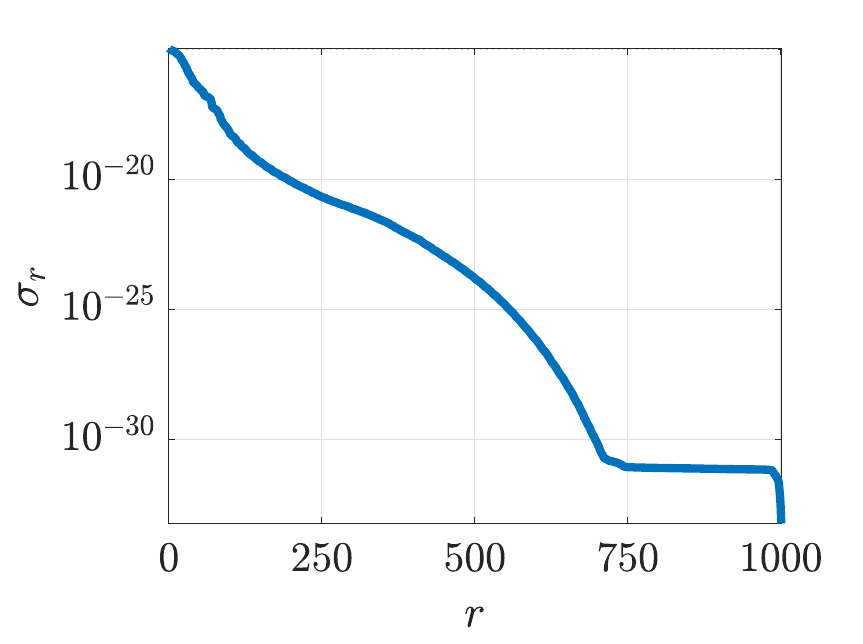}
        \includegraphics[width=0.33\linewidth]{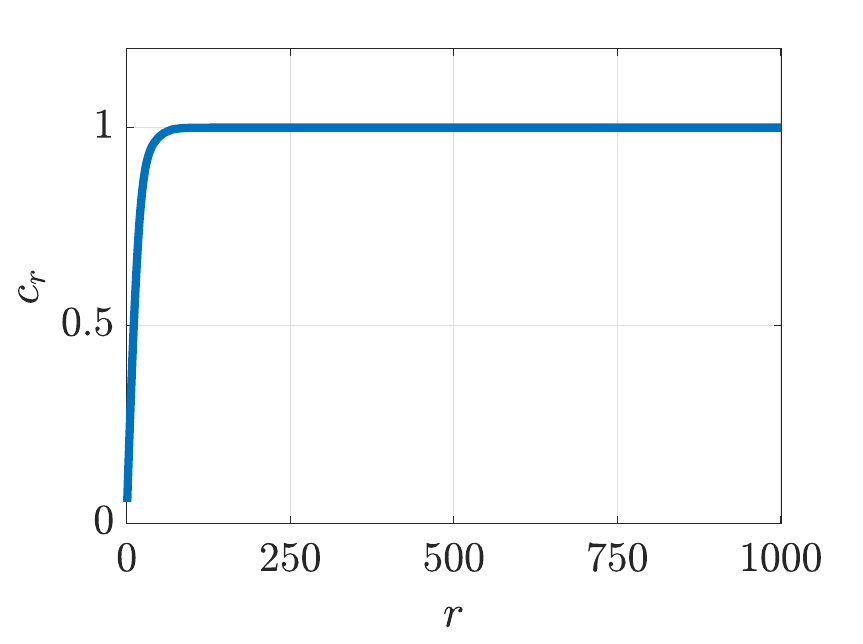}
  \caption{Decay of the singular values of the snapshot matrix of the displacement (left) and cumulative sum of the singular values (right) for the mechanical model of the damaged FML.}
\label{fig:sv_decay_damagedFML} 
\end{figure}

Next, we apply time continuous operator inference approach to this model. We recall that in order to apply operator inference one has to obtain the approximation of the second derivative data. In our case, the COMSOL software provides only the data for the solution and its first time derivative. Therefore, we employ an eighth order finite difference scheme, see \cite{Sharma2024}, in order to approximate accurately the second time derivatives of the solution. Then we proceed with obtaining the reduced basis by carrying out a SVD of the snapshot matrix and thereafter we constructed reduced versions of the displacement and its derivatives. The decay of the singular values of the snapshot matrix of the displacement is shown in Fig.~\ref{fig:sv_decay_damagedFML} (left). The cumulative sum of the singular values is shown in Fig.~\ref{fig:sv_decay_damagedFML} (right). We observe that at least the first $50$ singular values are needed for recovering $98$ $\%$ of the energy embedded in the snapshot data.\par 

The operator inference approach used for this example is the time continuous unconstrained version which does not assume full knowledge about the external loads, see Section~\ref{sec:Opinf2}. In this approach, we assume only the presence of the input signal (the excitation signal) mentioned earlier. Therefore, the dimension of the control matrix $\widehat{\mathbf{B}}_M$ on the right hand side of the reduced system is $r \times 1$.\par 

In this example, we have different scales of the snapshot matrices and the input signal. Therefore, we apply the scaling pre-processing step as explained in Section~\ref{sec:Opinf3}. The data matrix of the OpInf problem $\mathbf{\widehat{D}}$ is full rank and therefore one can solve the system of normal equations by inverting the matrix $\mathbf{\widehat{D}}\mathbf{\widehat{D}}^T$ to determine the unknown operators. The solution of the linear systems is performed using Matlab 2023b. We infer the reduced operators for different sizes of the reduced model in order to study the effect of the reduced dimension on the quality of the approximation. In particular, we simulate the OpInf problem for $r=10,20,\ldots,100$. We obtained the reduced operators of $\widehat{\mathbf{K}}_M$ and $\widehat{\mathbf{B}}_M$ for each setting. These reduced operators are then used to integrate the second order mechanical system using the Newmark method. We would like to remark that the results of the reduction given by any set of the learned operators are affected by the choice of the Newmark parameters, which are called $\beta$ and $\gamma$. We refer the reader to \cite[Section 5.3]{Chopra_Dynamics}, where two different types of the Newmark method are introduced, namely, the constant average acceleration method ($\beta = \frac{1}{4}, \, \gamma = \frac{1}{2}$) and the linear acceleration method ($\beta = \frac{1}{6}, \, \gamma = \frac{1}{2}$). We will show the results of the choice which has shown more accuracy which is the linear acceleration method, however, we will also show a comparison between the OpInf-ROM results for both choices for comparison sake.\par 

\begin{figure}[htbp]
\centering
\includegraphics[width=0.49\linewidth]{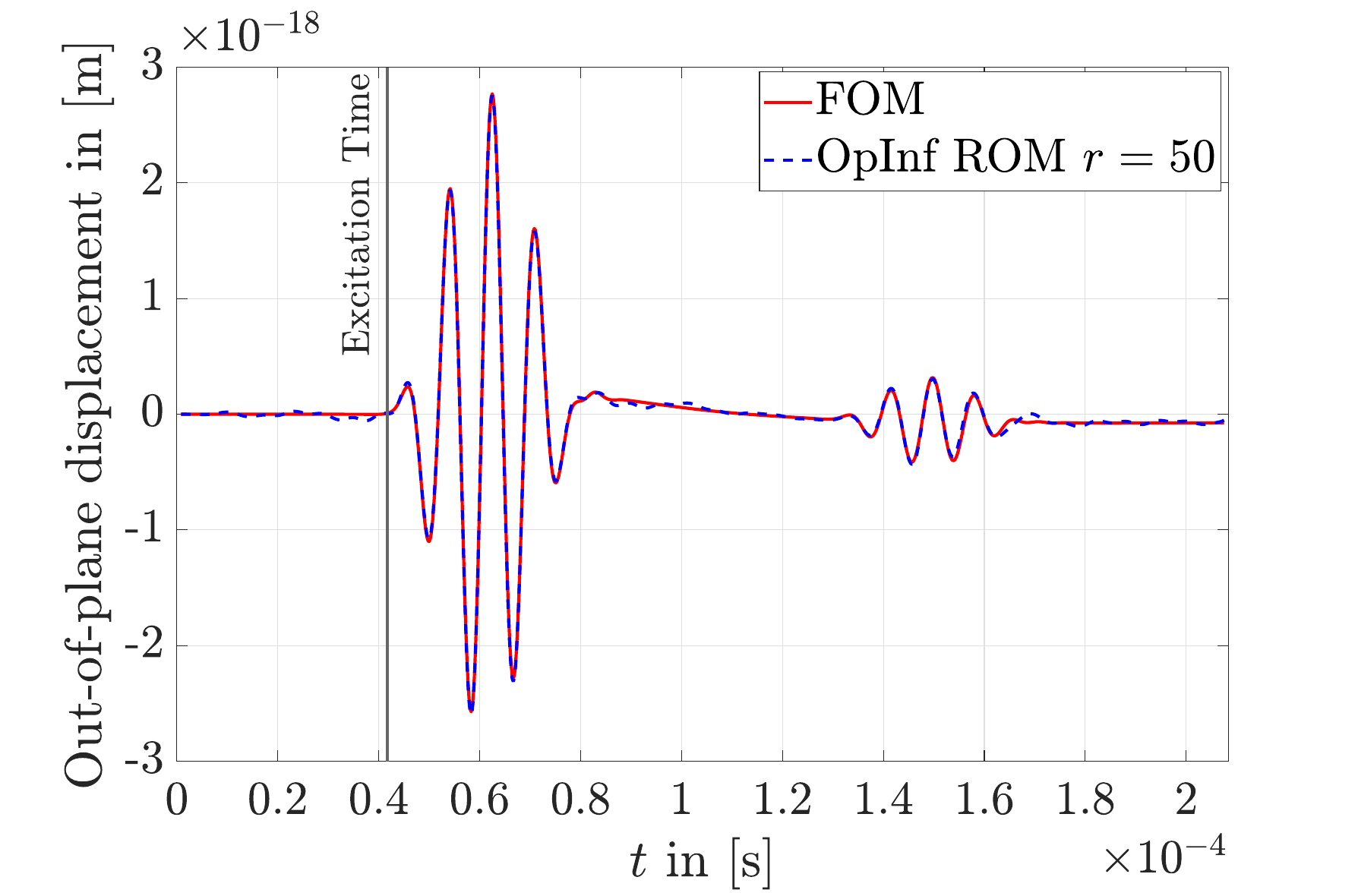}
\includegraphics[width=0.49\linewidth]{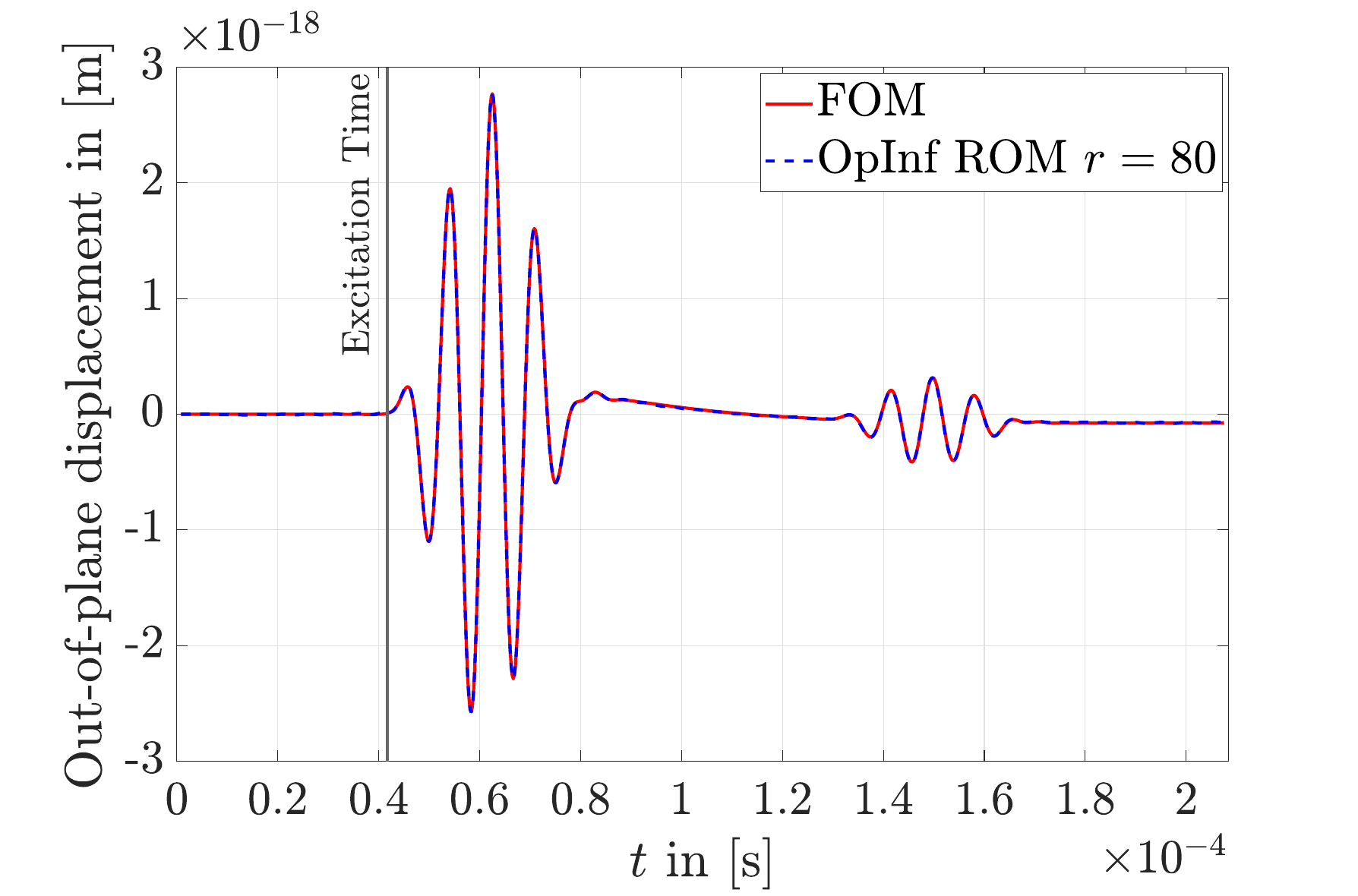}

  \caption{Reduced order model solution obtained via unconstrained OpInf for $r=50$ (left) and $r=80$ (right) versus the reference solution of COMSOL at the sensor for the mechanical model of the damaged FML.}\label{fig:opinf_rom_vs_comsol} 
\end{figure}

\begin{figure}[htbp]
  \centering
    \includegraphics[width=0.49\linewidth]{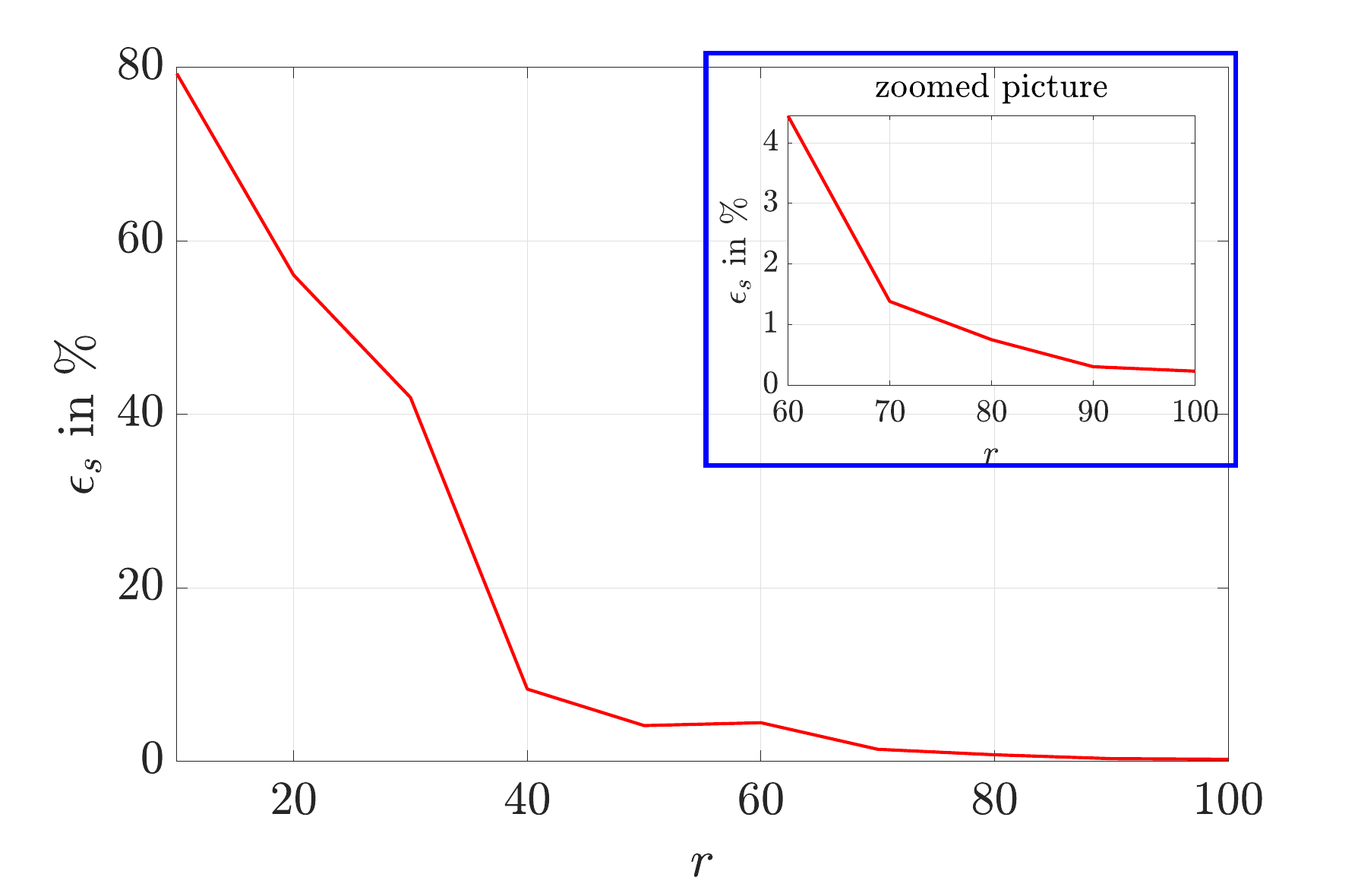} 
    \includegraphics[width=0.49\linewidth]{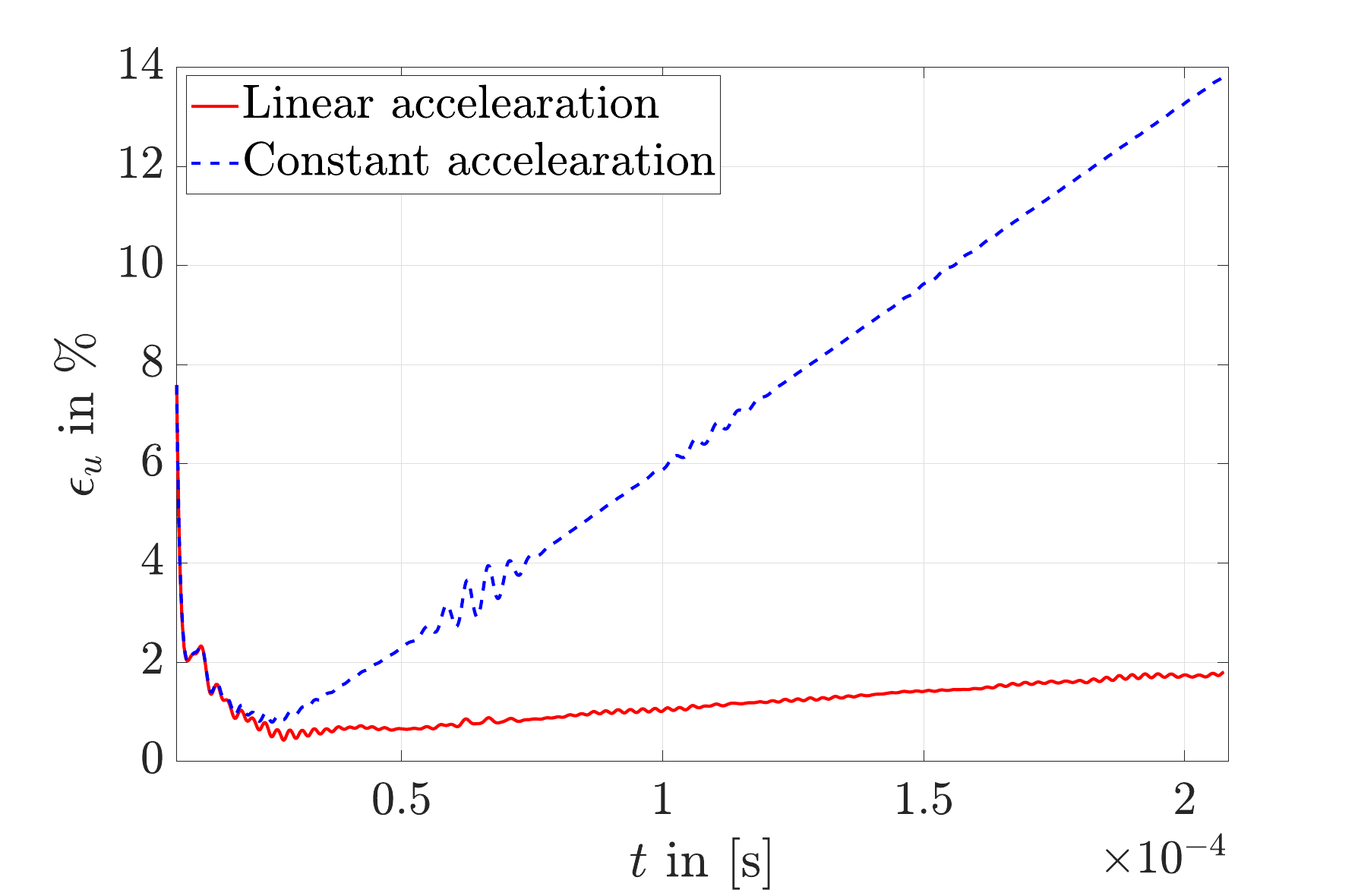}
  \caption{Relative $L^2$-error $\epsilon_s$ at the sensor location as function of the reduced dimension $r$ (left) and relative $L^2$-error $\epsilon_u(t)$ over the entire spatial domain as a function of time for the reduced dimension $r=80$ (right) of OpInf-ROM for the damaged FML model. The values are in percentages.}
\label{fig:epsilon_s_VS_r_FML} 
\end{figure}

The solution obtained by the OpInf based ROM and the reference solution of COMSOL are plotted at the sensor location in Fig.~\ref{fig:opinf_rom_vs_comsol}. The plots indicate that the OpInf approach allows for qualitatively accurate approximations of the reference solution for $r=50$ and it can be seen that the reduced approximation is improved when considering $r=80$. A more quantitative assessment of the approach is provided by computing the relative error at the sensor $\epsilon_s$ as defined in \eqref{eq:l2_error}. To analyze the convergence of the OpInf reduced approximation with respect to the reduced dimension, we compute the values of $\epsilon_s$ for increasing numbers of reduced basis functions. The values of the sensor error are plotted versus the reduced dimension $r$ in Fig.~\ref{fig:epsilon_s_VS_r_FML} (left) showing that the error reaches smaller values when increasing the reduced dimension $r$.\par 

To further evaluate the quality of the reduced approximation, we also examine the time-dependent relative error over the spatial domain $\epsilon_u (t)$, defined in \eqref{eq:l2_error_t}. The error $\epsilon_u (t)$ for $r=80$ is presented in Fig.~\ref{fig:epsilon_s_VS_r_FML} (right) for two different settings of the Newmark time integrator. The two settings which are described above provide different trends of the error as it can be seen from the plot. This shows that OpInf-ROM accuracy could vary significantly on the choice of the time discretization scheme. The spatial error (of the linear acceleration case) $\epsilon_u (t)$ remains under the threshold of $2$ $\%$ for the whole part of the simulation which shows that the OpInf-ROM is able to contain the error growth.

\subsection{Numerical results for the nonlinear hyperelastic model}\label{sec:nonlhypnum}
\begin{figure}[htbp]
  \centering
  \begin{minipage}[b]{1.0\linewidth}
    \centering
    \includegraphics[width=0.6\linewidth]{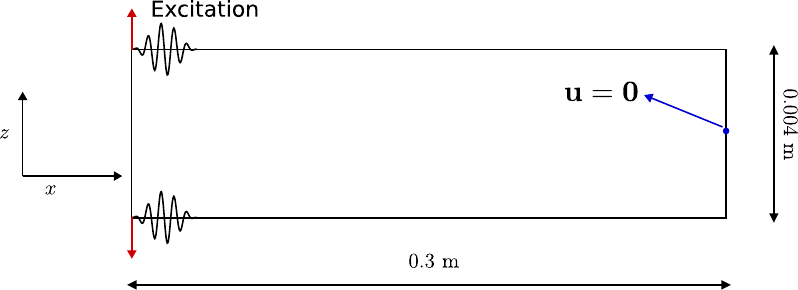} 
  \end{minipage}
  \caption{Aluminium model setup for the GUW simulation. The red arrows are used to indicate the sign of the excitation signal at each of the two nodes. A symmetric excitation is used for this numerical example through the displacement.}
\label{fig:Alu_sketch} 
\end{figure}

The third numerical example investigates a nonlinear elastic model of an aluminium plate. The model consists of a two-dimensional, intact rectangular aluminium plate with dimensions $0.3$ \si{m} (length) × $0.004$ \si{m} (thickness). The material parameters of aluminium are $\lambda = 5.1 \times 10^{10} \, \si{Pa}$ and $\mu = 2.6 \times 10^{10} \, \si{Pa}$. A sketch of the model is shown in Fig.~\ref{fig:Alu_sketch}. We note that this model has been investigated in our article \cite{GH24}, where we considered antisymmetric lamb wave A0. In contrast, here we consider a symmetric lamb wave mode S0. Actuators are positioned at the top and bottom left corners of the plate, applying a five-cycle Hanning window sinusoidal burst with a central frequency of $100$ kHz and a magnitude of excitation which is equal to $10^{-9}$. The material is modeled as hyperelastic, with nonlinearity captured using the Neo-Hookean model \cite{Rivlin1948}. The nonlinear effect of the Neo-Hookean material model can be seen e.g.\ in \cite[Figure 5.1]{RauterPhD}. The model employs displacement-based excitation, implemented as nonhomogeneous Dirichlet boundary conditions for the vertical displacement at the two left ends of the plate. A homogeneous Dirichlet boundary condition is set for the mid-point on the right side of the plate.  The simulation is run for $2.5 \times 10^{-4}$ \si{s}, with excitation active during the initial $5 \times 10^{-5}$ \si{s}. The simulation of the GUW is obtained by COMSOL-Multiphysics software\textcopyright.\par

\begin{figure}[htbp]
    \centering
         \includegraphics[width=0.33\linewidth]{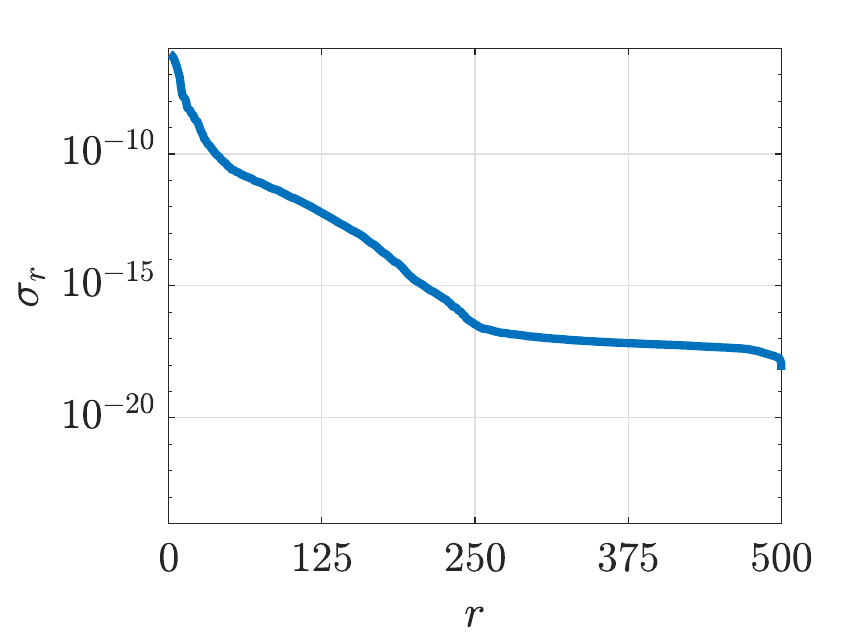}
     \includegraphics[width=0.33\linewidth]{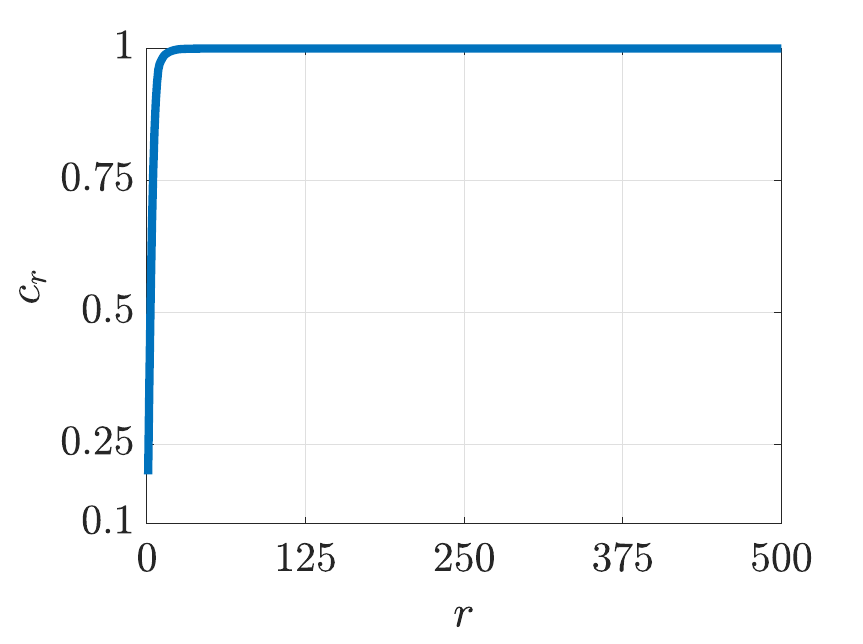}
  \caption{Decay of the singular values of the snapshot matrix of the displacement (left) and cumulative sum of the singular values (right) for the nonlinear hyperelastic model of aluminium.}
\label{fig:sv_decay_nonlModel} 
\end{figure}

In this example, we apply the DMD and mrDMD reduction techniques. A total of $M = 501$ time snapshots are collected, covering the full simulation window. The dimension of the FOM is $N=40834$. Figure~\ref{fig:sv_decay_nonlModel} (left) shows the decay of the singular values of the snapshot matrix of the displacement. The cumulative sum of the singular values is shown in Fig.~\ref{fig:sv_decay_nonlModel} (right). We can observe that with around $30$ singular values, the cumulative sum exceeds the threshold of $99.9$ $\%$ of the energy.\par 

\begin{figure}[htp]
\centering
\includegraphics[width=0.55\linewidth]{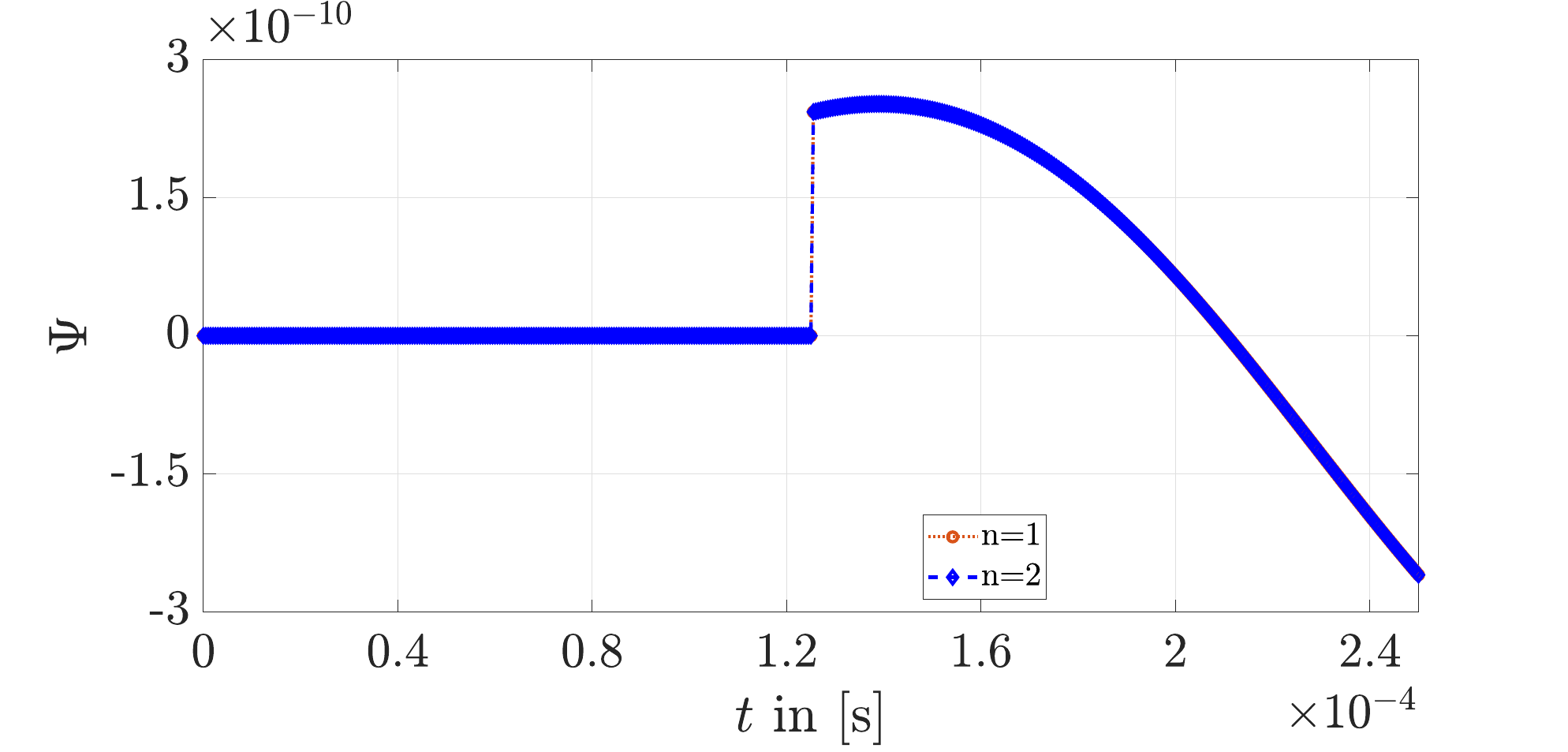}
\caption{The first-level time dynamics obtained by mrDMD for the nonlinear aluminium model.}
\label{fig:mrdmd_1stlevel_dynamics}
\end{figure}

The maximum number of mrDMD time levels is set to $6$. The first level of the mrDMD decomposition yields $2$ modes, whose time dynamics $\mathbf{\Psi} = \mathbf{b} \, \exp(\Omega t)$ are shown in Fig.~\ref{fig:mrdmd_1stlevel_dynamics}.

\begin{figure}[htp]
\centering
\includegraphics[width=0.495\linewidth]{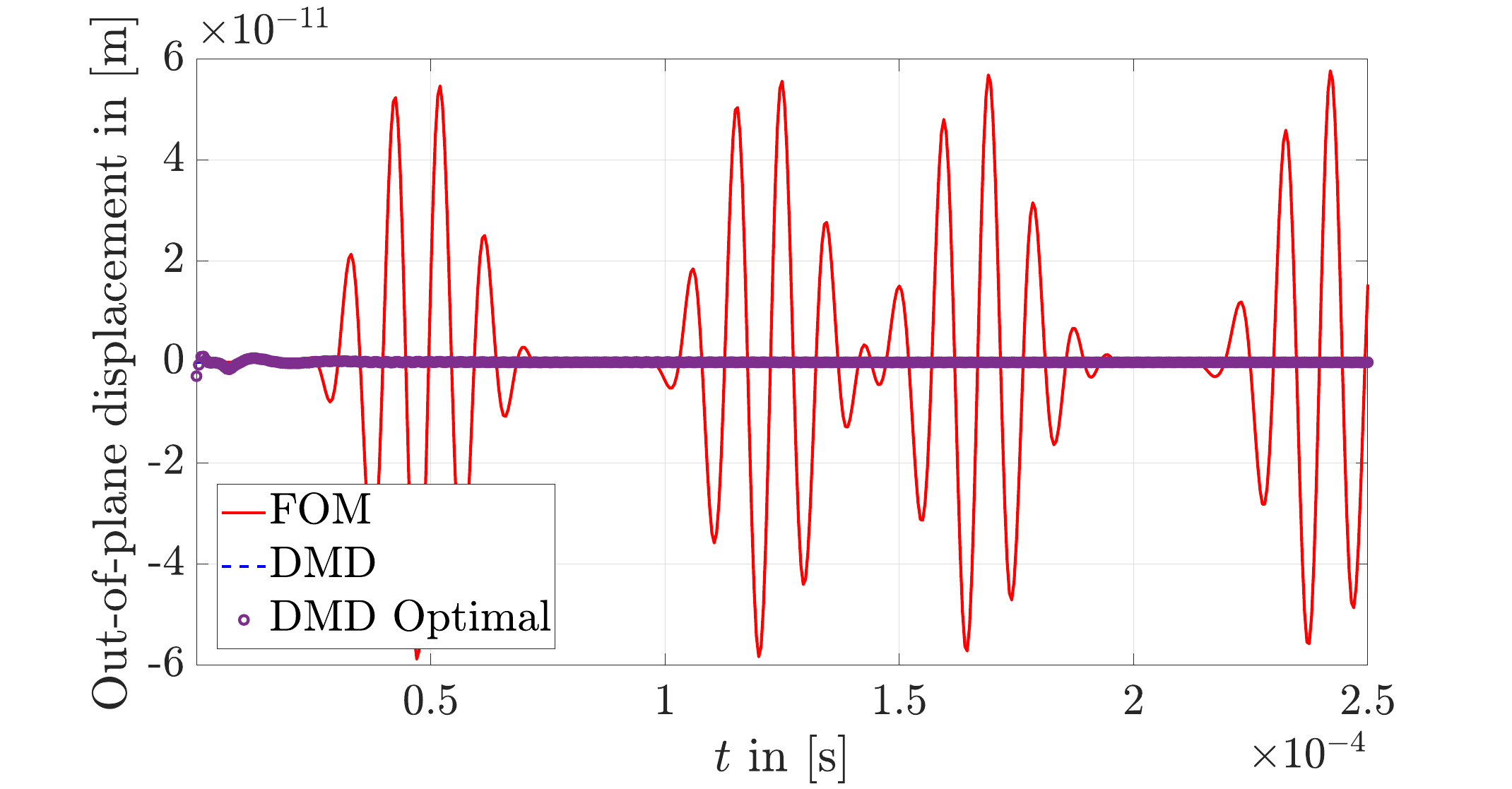}
\includegraphics[width=0.495\linewidth]{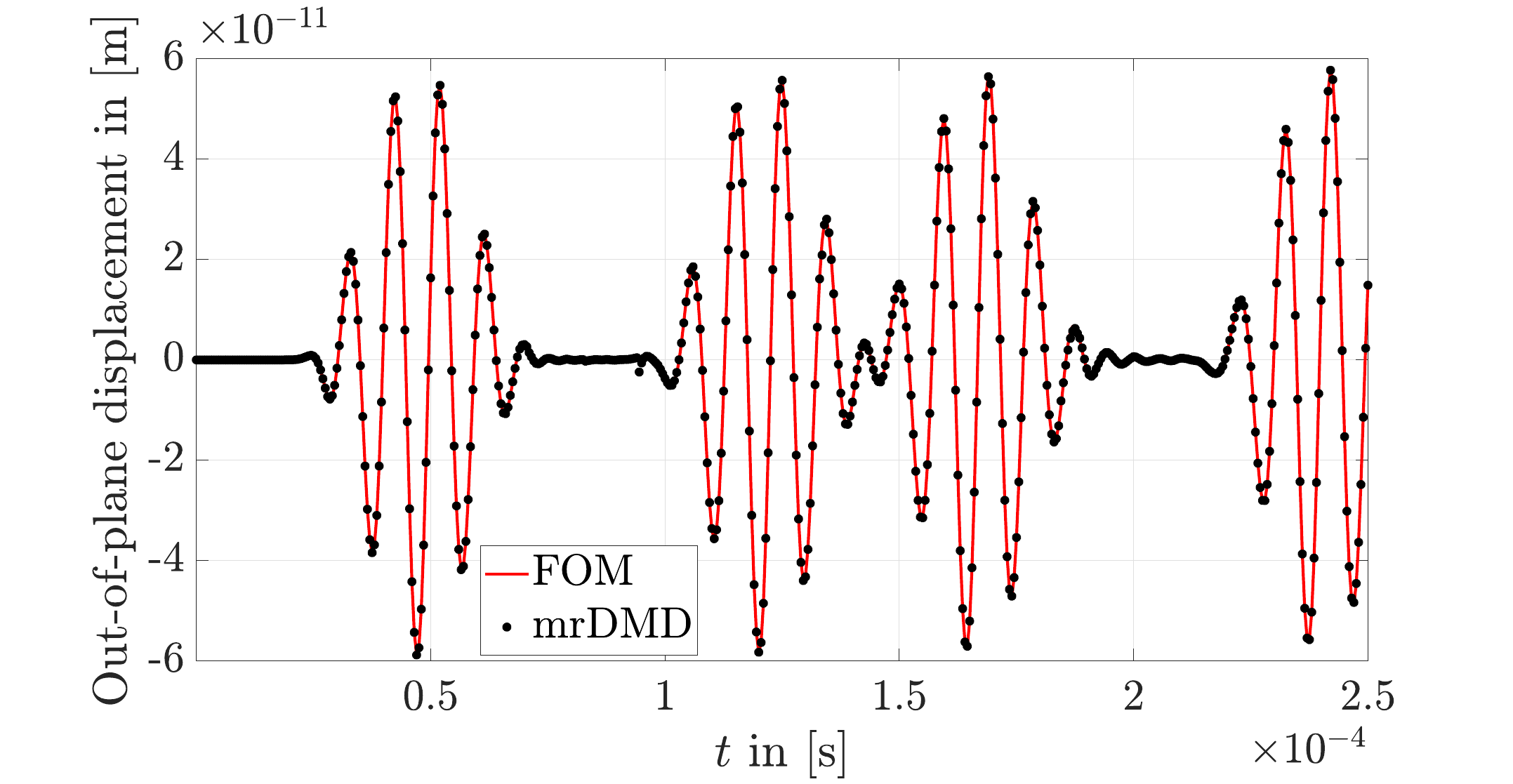}
\caption{Comparison between the FOM and DMD predicted displacement (left) and the mrDMD approximation (right) at the sensor for the nonlinear aluminium model. The reduced dimension used for DMD and DMD optimal is $r=80$.}
\label{fig:FOM_vs_DMD_AluModel_80modes}
\end{figure}

Standard DMD with a reduced dimension of $r=80$ is applied to the same snapshot set, with initial amplitudes computed using equations~\eqref{eq:amplitudes_using_projectedPOD} and \eqref{eq:amplitudes_optimal}. However, the method fails in both cases, as demonstrated in Fig.~\ref{fig:FOM_vs_DMD_AluModel_80modes} (left). This figure compares the out-of-plane displacement measured by the sensor for the reference COMSOL solution and the DMD approximations; standard DMD does not provide even qualitatively accurate results. This does not change as we increase the reduced dimension $r$. In contrast, Fig.~\ref{fig:FOM_vs_DMD_AluModel_80modes} (right) shows a similar comparison using mrDMD, where the alignment between the signals is significantly better. The relative error $\epsilon_s$ for the mrDMD approximation is $1.023$ $\%$, which is within an acceptable range.

These results highlight a key limitation of classical DMD: it assumes linear dynamics and often struggles to capture strongly nonlinear or transient behavior. In contrast, the multiresolution structure of mrDMD enables better isolation of localized temporal features, making it more suitable for modeling nonlinear systems such as the one considered here.

\begin{figure}[htp]
\centering
\includegraphics[width=0.55\linewidth]{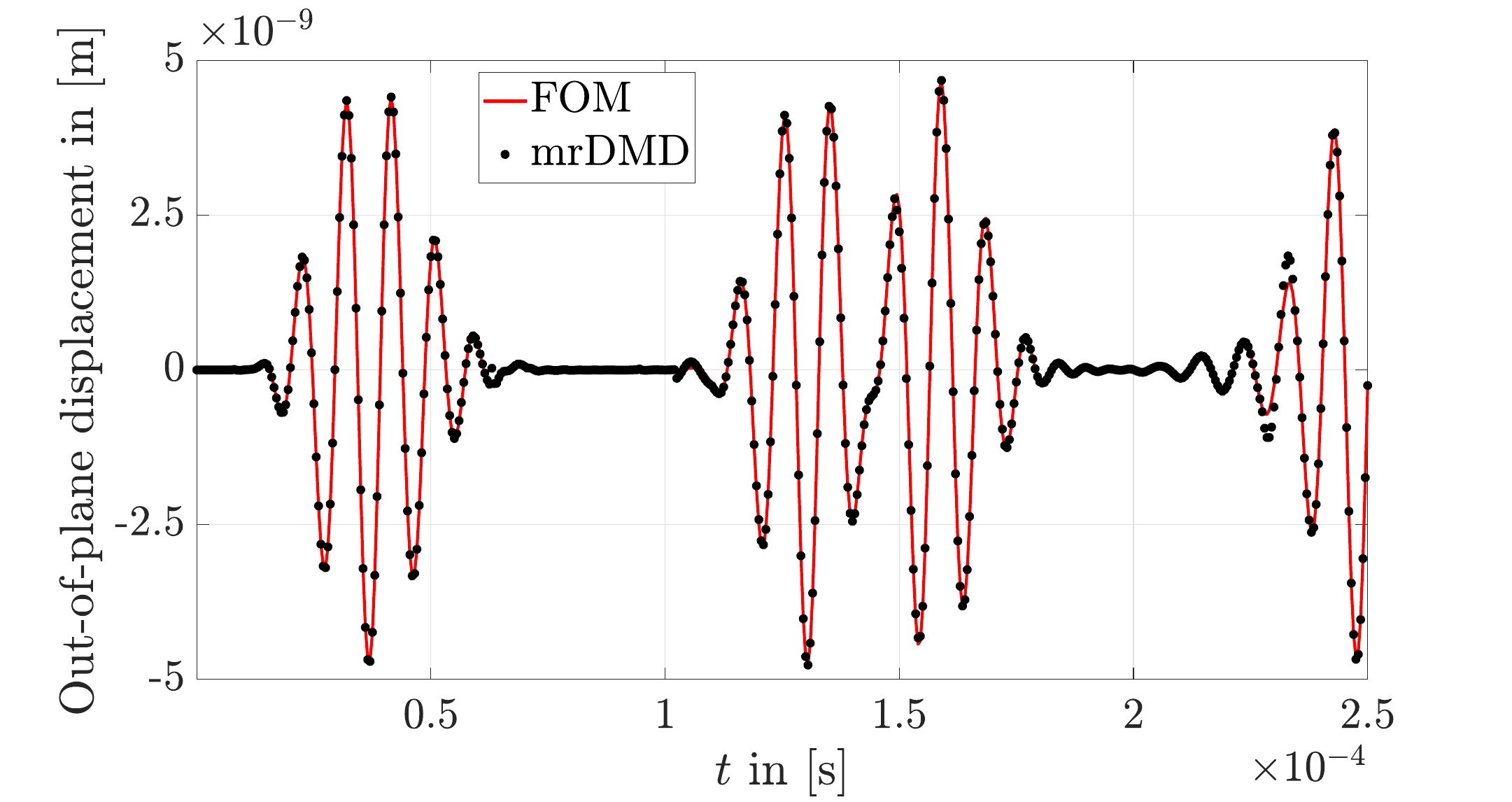}
\caption{Comparison between the FOM and mrDMD predicted displacement at the sensor for the nonlinear aluminium model with damage.}
\label{fig:FOM_vs_DMD_AluDamagedModel_80modes}
\end{figure}

In the final numerical test presented in this section, we consider a damage scenario in the aluminium plate, where the damage is modelled as a reduction of the Young’s modulus of elasticity $E$ within a rectangular subdomain. The Young’s modulus of the undamaged aluminum is $E = 72 \times 10^9$ \si{Pa}. In the damaged region, the reduced modulus $E_d$ is set to $10^4$ \si{Pa}. The damaged zone is defined by the subdomain $[0.06, 0.07] \times [0.00124, 0.00136]$.
The model is then simulated, and displacement snapshots are extracted in order to apply the mrDMD method. In this example, we employ $8$ levels of mrDMD. The results of the full order model and the mrDMD reconstruction at the sensor location are shown in Fig.~\ref{fig:FOM_vs_DMD_AluDamagedModel_80modes}. The figure demonstrates that mrDMD is able to qualitatively reconstruct the FOM solution. Quantitatively, the sensor error $\epsilon_s$ is equal to $4.388$ $\%$.

\section{Outlook}\label{sec:outlook}
In this work, we apply both intrusive and non-intrusive model reduction approaches to the wave equation posed on a domain containing damage. Subsequently, we apply non-intrusive reduction techniques to discretized second order mechanical systems obtained from guided ultrasonic wave (GUW) simulations. The test cases considered include models with damage and / or nonlinear behavior.\\
Damage identification and characterization — determining both the location and severity of defects — are central challenges in the field of structural health monitoring (SHM). This motivation underlies our current research  to develop a reduced order model (ROM) capable of generalizing across different damage configurations. To achieve this, we parameterize the damage and address the associated inverse problem, which involves estimating unknown damage parameters from observed system outputs. In this case, the damage parameters influence the full order stiffness matrix and, consequently, the displacement field. This dependency is carried over to the projected reduced operator of the stiffness matrix, i.e.\ $\mathbf{K}_r = \mathbf{K}_r(\bm{\theta}) $. Since the reduced stiffness matrix is not parameter separable in this case, for each new parameter $\bm{\theta}$ the construction of the reduced stiffness matrix $\mathbf{K}_r(\bm{\theta})$ requires the assembly of the full order stiffness matrix $\mathbf{K}(\bm{\theta})$.  Thus, the reduced model remains coupled to the high dimension $N$. This bottleneck for multi-query scenarios is well-known such that hyper reduction techniques are developed like matrix gappy POD \cite{matrixgappypod} or matrix DEIM \cite{mdeim}. However, applying such methods requires that the full order solver or computational code allows direct computation of specific entries of the stiffness matrix without assembling the entire matrix. This is not feasible in our case, since we employ COMSOL Multiphysics as the simulator for the GUW problem, and direct access to the full order matrices or the time discretization scheme is restricted.
Consequently, we extend the existing OpInf-based ROM framework to handle parameterized systems. The proposed extension must be designed to accurately capture the dependence of the system outputs — in this case, the displacement field — on the damage-related input parameters.

\section{Acknowledgement}
The authors would like to thank Peter Benner for helpful discussions on the topic. The authors would like to thank our student assistants for their contributions. Moreover, the authors acknowledge funding by the German Research Foundation (Deutsche Forschungsgesellschaft, DFG), project number 418311604 (Research Unit FOR3022: Ultrasonic Monitoring of Fibre Metal Laminates Using Integrated Sensors).

\section{Research data for this article}
The numerical source code and datasets generated during the current study are publicly available to support the reproducibility of the presented results:

\begin{itemize}
    \item \textbf{Dataset:} The numerical data used to generate the figures and benchmarks are deposited on Zenodo and accessible via DOI: \href{https://doi.org/10.5281/zenodo.23194087}{10.5281/zenodo.23194087} \cite{hijazi_wave_problems_2026}.
    \item \textbf{Source Code:} The MATLAB and Python implementation files associated with all numerical examples in this manuscript are hosted on GitHub at \url{https://github.com/saddamhijazi/Codes_Manuscript_DataDrivenMOR_for_WavePbs/} \cite{hijazi_code_2026}.
\end{itemize}

\bibliographystyle{elsarticle-num}
\bibliography{bibfile} 
\end{document}